%% file: main.tex
\RequirePackage{fix-cm}
\documentclass[smallextended]{svjour3}       
\smartqed  
\AtBeginDocument{%
  \paperwidth=\dimexpr
    1in + \oddsidemargin
    + \textwidth
    + 1in + \oddsidemargin
  \relax
  \paperheight=\dimexpr
    1in + \topmargin
    + \headheight + \headsep
    + \textheight
    + 1in + \topmargin
  \relax
  \usepackage[pass]{geometry}\relax
}

\usepackage{graphicx}
\usepackage{amsmath}
\usepackage{amssymb}
\usepackage{xcolor}
\usepackage[caption=false]{subfig}
\usepackage{mathrsfs}
\usepackage{pgfplots}
\pgfplotsset{compat = 1.18}
\usepackage{csvsimple}
\usepackage{longtable}
\usepackage{siunitx}

\usepackage{hyphenat}

\usepackage{todonotes}

\usepackage[%
    colorlinks=true,
    pdfborder={0 0 0},
    linkcolor=red,
    citecolor=blue
]{hyperref}

\newtheorem{rule_def}{Consistency Rule}
\begin{document}

\title{Global Optimization Algorithm for Mixed-Integer Nonlinear Programs with Trigonometric Functions
}
\subtitle{}


\author{Christopher Montez \and Sujeevraja Sanjeevi \and Kaarthik Sundar 
}


\institute{C. Montez \at
    Department of Mechanical Engineering \\
    Texas A\&M University, College Station, TX \\
    \email{christopher.martin.montez@gmail.com}         
    \and
    S. Sanjeevi \at 
    OpsLab, Austin, TX \\ 
    \email{sujeev.sanjeevi@gmail.com}
    \and 
    K. Sundar (corresponding author) \at
    Information Systems \& Modeling Group, \\ 
    Los Alamos National Laboratory, Los Alamos, NM \\
    \email{kaarthik@lanl.gov}
}

\date{Received: date / Accepted: date}

\titlerunning{Global Optimization Algorithm for MINLPs with Trigonometric Functions}

\maketitle

\begin{abstract}
This article presents the first mixed-integer linear programming (MILP)-based iterative algorithm to solve factorable mixed-integer nonlinear programs (MINLPs) with bounded, differentiable periodic functions to global optimality with an emphasis on trigonometric functions.
At each iteration, the algorithm solves a MILP relaxation of the original MINLP to obtain a bound on the optimal objective value. The relaxations are constructed using 
partitions of variables involved in each nonlinear term and across successive iterations, the solution of the relaxations is used to refine these partitions further leading to tighter relaxations. Also, at each iteration, a heuristic/local solve on the MINLP is used to obtain a feasible solution to the MINLP. The iterative algorithm terminates till the optimality gap is sufficiently small. This article proposes novel refinement strategies that first choose a subset of variables whose domain is refined, refinement schemes that specify the manner in which the variable domains are refined, and MILP relaxations that exploit the principal domain of the periodic functions. We also show how solving the resulting MILP relaxation may be accelerated when two or more periodic functions are related by a linking constraint. This is especially useful as any periodic function may be approximated to arbitrary precision by a Fourier series.  Finally, we examine the effectiveness of the proposed approach by solving a path planning problem for a single fixed-wing aerial vehicle and present extensive numerical results comparing the various refinement schemes and techniques.

\end{abstract}

\section{Introduction}
\input{Sections/intro}
\section{Initial Overview of the Algorithm}
\input{Sections/overview}

\section{Preliminaries}
\input{Sections/preliminaries}

\section{MILP Relaxations}
\input{Sections/relaxations}

\section{Partitions}
\input{Sections/tighten}

\section{Principal Domains for Periodic Functions}
\input{Sections/domains}

\section{Motivating Example - Markov-Dubins Path Planning Problem}
\input{Sections/dubins}

\section{Computational Results}
\input{Sections/results_updated}

\section{Conclusion}
\input{Sections/conclusion}

\appendix 
\section{Appendix} 
\input{Sections/appendix_new}

\bibliographystyle{spmpsci}
\bibliography{refs}

\end{document}

%% file: Sections/intro.tex
Optimization problems in many applications such as chemical process networks \cite{lee1996mixed,jia2003mixed,yadav2012short,rubio2013global}, energy systems \cite{catalao2011hydro}, and wastewater treatment \cite{galan1998optimal} to name a few, are traditionally modeled as Mixed Integer Nonlinear Programs (MINLPs). MINLPs are mathematical programs that include both continuous and discrete decision variables, where the objective function and/or the constraints may be nonlinear and possibly non-convex. Algorithms to solve this class of optimization problems to global optimality have garnered extensive attention from both academia and industry, resulting in the development of both open-source (Couenne \cite{Couenne} and SCIP \cite{SCIP}) and commercial solvers (BARON \cite{BARON}, LINDOGlobal \cite{LINDO}, Gurobi \cite{gurobi}, and ANTIGONE \cite{ANTIGONE}) for the same. Each of these solvers is specialized and solves a subset of factorable \cite{mccormick1976computability} MINLPs where nonlinearities arise due to a certain class of functions such as multilinear, logarithmic, exponential, etc. Since the focus of this article is on MINLPs with periodic functions with an emphasis on trigonometric functions, we remark that among all these global optimization solvers, to the best of our knowledge, only Couenne, LINDOGlobal, SCIP, Gurobi, and FICO XPress Global implement global optimization algorithms that are equipped to solve MINLPs with trigonometric terms.

Algorithmic approaches to solving MINLPs to global optimality typically contain two main features: convex relaxations and search. For instance, any of the aforementioned solvers for MINLPs first isolate each non-convex term in the problem and construct convex over- and under-estimators for these terms. Doing so for every non-convex term in the MINLP yields a convex relaxation of the MINLP, which in turn is solved to optimality to provide a bound to the optimal objective value of the MINLP. This process is typically combined with a standard search procedure such as spatial branch-and-bound (sBB) \cite{horst2013global,Smith1997} to explore the full space of feasible solutions to the MINLP. Furthermore, at each node of the sBB tree, a local solve or heuristic is applied to keep generating feasible solutions to the MINLPs and to keep track of the optimality gap. Mathematical properties of sBB such as consistency and exhaustiveness \cite{horst2013global,Floudas2014} in turn provide a guarantee of convergence of the algorithm to global optimality. Each solver differs in how the convex relaxations are constructed, the strength of the respective convex relaxations, and the support for the different types of non-convex terms it provides, resulting in disparate computational performance for the same MINLP.

More recently, global optimization algorithms for MINLPs that rely on solving a Mixed Integer Linear Program (MILP) at each iteration have gained considerable interest \cite{Wicaksono2008,Misener2011,Teles2013,Castillo2018,Nagarajan2019}. The main motivation for developing MILP-based algorithms for MINLPs has been the meteoric improvement in the speed of off-the-shelf MILP solvers \cite{GurobiBenchmarks}. While the sBB-based algorithms rely on constructing and solving convex relaxations of non-convex structures, MILP-based algorithms rely on constructing and solving piecewise convex or polyhedral relaxations of non-convex structures \cite{sundar2021sequence}. The search procedure for MILP-based methods are domain partitioning schemes which partition the domains of the variables involved in the non-convex terms \cite{Nagarajan2019}. Similar to the sBB search procedure, the domain partitioning schemes require the mathematical properties of exhaustiveness and consistency to ensure convergence. It has been shown in recent work \cite{Nagarajan2019,Castro2016} that the MILP-based methods perform comparably with sBB-based methods and sometimes even outperform the sBB counterpart for MINLPs with multilinear and quadratic functions. The work in \cite{Nagarajan2019} that deals with MINLPs with multilinear and quadratic functions has resulted in a solver by the name of Alpine.jl, written using the Julia programming language.

Several notable contributions exist in this space. For example, \cite{lundell2013reformulation} developed a method based on transforming non-convex MINLPs into convex MINLPs using $\alpha$BB and spline-$\alpha$BB under-estimators, subsequently solved by an outer approximation algorithm. While their method is effective, it ultimately solves convex MINLPs at each iteration. In contrast, our work aims to directly leverage advancements in MILP technology by entirely avoiding the need to solve non-linear problems. Other works such as \cite{lundell2012global,westerlund2006some} focus on MINLPs with signomial functions and use transformations involving single-variable powers and exponentials, which are then approximated using piecewise linear functions to construct valid lower bounds. These works also introduced various refinement strategies, such as breakpoint refinement at specific locations (e.g., midpoints, most deviating points) -- strategies which we incorporate and extend. Our work goes beyond by employing non-uniform partitioning approaches informed by recent successes in the field \cite{Nagarajan2019}, offering flexibility and better approximation in challenging regions. An extension of these ideas was proposed in \cite{lundell2018solving}, where a unified framework combining $\alpha$BB and signomial transformations was used to extend applicability to more general MINLPs.

In this article, we develop piecewise polyhedral relaxations for non-convexities arising from bounded, differentiable periodic terms and embed these relaxations in an MILP-based algorithmic framework for global optimization of underlying MINLPs. To the best of our knowledge, this is the first article to develop a MILP-based global optimization algorithm for MINLPs with bounded, differentiable periodic functions. Existing open-source and commercial solvers that support trigonometric functions like Gurobi, SCIP, LINDOGlobal use the sBB framework to solve the MINLPs to global optimality. The framework we develop in this article can be easily integrated into existing MILP-based MINLP solvers like Alpine to handle periodic terms, such as trigonometric terms, which originally could not be considered. 

For the remainder of this article, we restrict our focus to trigonometric functions, but the presented work can be easily extended to other periodic functions of interest. We have chosen to do this as (i) trigonometric functions, namely sine and cosine, are prevalent in many models and (ii) any bounded, periodic function can be approximated to arbitrary precision using a truncated Fourier series. Despite this restricted focus, the work presented in this article can be easily applied to other periodic functions should it be more desirable to retain the original periodic function rather than replace it with an approximation.

The fundamental contributions of this article are (i) development of polyhedral and piecewise polyhedral relaxations for univariate trigonometric functions (which may be extended to general univariate, differentiable periodic functions), (iii) techniques to reformulate these relaxations when the principle domain of the periodic function is greater than its period and its impact in the computation times, (iii) integration of these relaxations into an MILP-based algorithmic framework for global optimization of MINLPs, (iv) development of novel variable selection for partitioning strategies and non-uniform partitioning schemes for variable domain partitioning, and
(v) illustration of the effectiveness of the overall algorithm to solve a MINLP that models a path planning problem for a single fixed-wing vehicle. The main reason for choosing the path planning problem as a use case is that this is one of the most common application in robotics which leads to an MINLP with sine and cosine function and with domain of these functions being greater than its period, $2\pi$. This feature presents a substantial challenge for global optimization algorithms and the proposed global optimization algorithm with all the enhancements is shown to have substantial computational gains for this problem.

\subsection{Problem Statement} \label{subsec:problem}
In this paper we consider the problem of finding the globally optimal solution\footnote{Global optimality is defined numerically by a specified tolerance, $\epsilon$.} of an MINLP composed of linear, trigonometric, and bilinear terms that has been factored into the form

\begin{subequations} \label{P}
    \begin{align}
        (\mathcal{F}) \qquad \text{minimize } \quad &c^Tx + d_1^Ty + d_2^Tz \\
        \text{subject to} \quad &A_1 x + A_2 y + A_3 z\leq 0 \label{P inequality}\\
        &y_i = f_i(x) \,\,\, \qquad i = 1, \hdots, n_t \label{univariate term}\\
        &z_j = g_j(x, y) \quad j = 1, \hdots, n_b \label{bilinear term} \\
        &x \in X \subseteq \mathbb{Z}^{n_i} \times \mathbb{R}^{n - n_i}
    \end{align}
\end{subequations}
In \eqref{P}, the vector $x$ consists of continuous and integer-valued variables and represents the variables of the original MINLP. The vectors $c$, $d_1$, and $d_2$ and the matrices $A_1$, $A_2$, and $A_3$ are of appropriate size and are the result of the factoring procedure used to convert the original MINLP into the form of \eqref{P}.
The functions $f_i \colon I_i \subset \mathbb R \mapsto \mathbb R$ for $i = 1, \hdots n_t$, are univariate, trigonometric functions that are differentiable and bounded over a closed interval $I_i$, where $n_t$ is the number of unique trigonometric terms after factoring. The functions $g_j \colon R_j \subset \mathbb{R}^2 \mapsto \mathbb R$ for $j = 1, \hdots, n_b$ are bilinear terms. Apart from bilinear terms involving two variables in $x$, \eqref{bilinear term} may also consist of the one original variable in $x$ and one auxiliary variable in $y$ or two auxiliary variables in $y$. The functions $g_j$ are each defined over a closed rectangle $R_j$ and $n_b$ represents the number of unique bilinear terms after factoring. We assume the feasible space of \eqref{P} is non-empty. We note that multilinear terms may also be considered by recursively defining bilinear variables. In doing so, $g_j$ would also be a function of the introduced $z$ variables.

We approach finding the globally optimal solution to $\mathcal{F}$ by introducing MILP relaxations for the trigonometric terms \eqref{univariate term} and the bilinear terms \eqref{bilinear term}. For the univariate trigonometric functions $f_i$ with domain $I_i$, we partition $I_i$ using information of the convexity of $f_i$ over sub-intervals of $I_i$. An MILP relaxation of each $f_i$ is then constructed using its corresponding partition where the MILP relaxation is the disjunctive union of triangles containing $f_i$. For the bilinear terms $g_j$ with domain $R_j$, we partition $R_j$ by dividing the domain of one variable and leaving the second variable's domain untouched. An MILP relaxation of each $g_j$ is then constructed using its corresponding partition where the MILP relaxation is the disjunctive union of tetrahedrons containing $g_j$. All the relaxations are formulated using a so-called ``incremental formulation'' \cite{yildiz2013incremental}. Once all trigonometric and bilinear terms have been relaxed, the resulting MILP is solved to find a lower bound. This lower bound is successively tightened by refining some or all of the partitions used to construct the MILP relaxations of the trigonometric and bilinear terms. This process is repeated until the gap is sufficiently small. It is assumed a method exists for determining an upper bound given a lower bound. If no such method exists, a standard local nonlinear programming (NLP) solver can be used to find a feasible solution and hence, an upper bound to the optimal objective value. This algorithm we present can be integrated into existing solvers to handle a larger class of MINLPs than those considered in this paper.

\subsection{Structure of the Paper}
The remainder of this paper is organized as follows. In Section \ref{Overview Section}, we present a simplified flowchart of the algorithm presented in this paper to serve as a guide for the reader. In Section \ref{Preliminaries}, we present preliminary definitions which will be used throughout the remainder of the paper. In Section \ref{MILP Relaxations Section}, we present the MILP relaxations for trigonometric and bilinear terms. In Section \ref{Partitions Section}, we present an adaptive partition refinement scheme that will be used to successively tighten the MILP relaxation until global optimality is reached. In Section \ref{Principal Domains Section}, we introduce the concept of principle domains and discuss their impact on solving MINLPs with periodic functions. In Section \ref{Motivating Example Section}, we present a motivating example in the form of the path planning of a vehicle moving in the plane through a specified sequence of points subject to kinematic constraints. In Section \ref{Computational Results Section}, we present computational results. In Section \ref{Conclusion}, we give concluding remarks and discuss potential future work. 

%% file: Sections/overview.tex
\label{Overview Section}

\begin{figure}
    \centering
    \includegraphics[width = 0.9\textwidth]{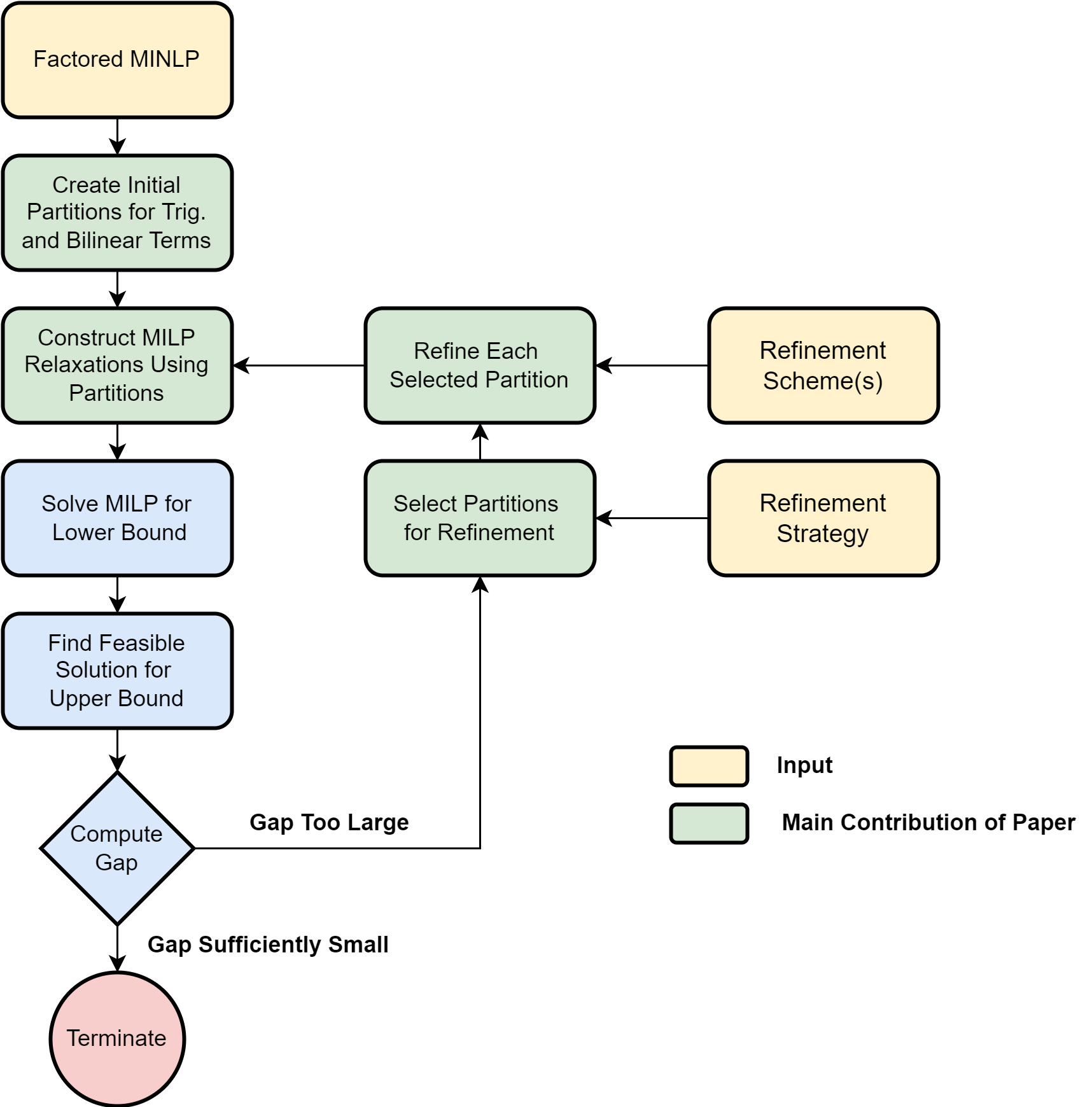}
    \caption{Simplified flowchart for the proposed algorithm.}
    \label{fig: flowchart}
\end{figure}

A simplified flowchart of the algorithm is shown in Figure \ref{fig: flowchart}. The algorithm first takes an MINLP composed of linear, trigonometric, and bilinear terms in the factored form \eqref{P}. An initial partition is created for each trigonometric and bilinear term in the factored MINLP. 
The procedure for constructing these initial partitions is the subject of Sections \ref{Preliminaries}, \ref{MILP Relaxations Section}, and \ref{Partitions Section}. Using these initial partitions, an MILP relaxation for each trigonometric and bilinear term in the factored MINLP is constructed. The procedure for constructing MILP relaxations of trigonometric and bilinear terms given a set of partitions is the subject of Section \ref{MILP Relaxations Section}. The MILP relaxation is then solved using a standard MILP solver (such as CPLEX \cite{cplex2009user} or Gurobi \cite{gurobi}) to obtain a lower bound on the original MINLP's optimal solution. It is assumed a method is available to find a feasible solution to the MINLP.  
This method is then used to find a feasible solution, the best known upper bound is updated, and the resulting relative gap is computed. If the gap is sufficiently small (using a user-specified gap tolerance $\varepsilon$), the algorithm terminates. If the gap is too large, some or all partitions are selected for refinement to tighten the relaxation. The procedure used for selecting partitions for further refinement is referred to as a refinement strategy in this paper. Refinement strategies are the subject of Section \ref{Refinement Strategy}. Each of the selected partitions are then refined. The procedure used for refining a given partition is referred to as a refinement scheme in this paper. Refinement schemes are the subject of Section \ref{Refinement Scheme Section}. Once each of the selected partitions are further refined, a new, tighter MILP relaxation is constructed using these partitions. The tighter MILP relaxation is solved and this process is repeated until the gap is less than $\varepsilon$. 



%% file: Sections/preliminaries.tex
\label{Preliminaries}

We first introduce some terminology that will be used throughout the remainder of the paper. The terminology we use closely follows the work of \cite{sundar2021sequence}.

\begin{definition} \label{def: partition}
    Given a closed interval $[l, u] \subset \mathbb R$, a \textit{partition} $p$ of $[l, u]$ is an ordered sequence of distinct rational numbers 
    $(x_0, x_1, \hdots, x_m)$ such that $l = x_0 < x_1 < \hdots < x_m = u$. We denote the set of partition points of $p$ by $\mathcal P(p) = \{x_0, x_1, \hdots, x_m\}$.
\end{definition}

\begin{definition} \label{def: break point}
    Given a nonlinear, bounded, and differentiable function \\ $f \colon [l, u] \mapsto \mathbb R$, we say $b \in [l, u]$ is a \textit{break point} of $f$ if at $b$ the function changes from being convex to concave or vice-versa. 
\end{definition}

\begin{definition} \label{def: admissible}
    Given $f \colon [l, u] \mapsto \mathbb R$, we say a partition $p$ is \textit{admissible} if $p$ contains all break points of $f$ in $[l, u]$.
\end{definition}


\begin{definition} \label{def: base partition}
    Given a function $f \colon [l, u] \mapsto \mathbb{R}$, the partition $p^0$ of $[l, u]$ is referred to as a \textit{base partition} of $[l, u]$ for $f$ if (i) it is admissible, (ii) for any sub-interval $[x_i, x_{i + 1}]$ defined by the partition, $f'(x_i) \neq f'(x_{i + 1})$ (slope condition), and (iii) $|\mathcal{P}(p^0)|$ is minimum.  
\end{definition}

\begin{definition} \label{def: valid refinement}
    Given $f \colon [l, u] \mapsto \mathbb R$ and admissible partition $p$ of $[l, u]$ for $f$, we say the partition $q$ is a \textit{valid refinement} of $p$ if $\mathcal{P}(p) \subset \mathcal{P}(q)$ and for each sub-interval $[x_i, x_{i + 1}]$ defined by $q$ we have $f'(x_i) \neq f'(x_{i + 1})$. 
\end{definition}

As an example, consider the function $f(x) = \sin x$ for $x \in [0, 2\pi]$ and the partition $p = (0, \frac{\pi}{2}, \pi, \frac{3\pi}{2}, 2\pi)$. The break points of $f$ in the domain $[0, 2\pi]$ are 0, $\pi$, and $2\pi$ and so $p$ is admissible. We also see the slope conditions (ii) in Definition \ref{def: base partition} are satisfied. However, $\mathcal{P}(p)$ satisfying conditions (i) and (ii) in Definition \ref{def: base partition} is not of minimum size and so $p$ is not a base partition. The partition $p^0 = (0, \pi, 2\pi)$ is a base partition of $[0, 2\pi]$ for $f$ and so $p$ is a valid refinement of $p^0$. In general, to construct a base partition we start with all break points and only add additional points to the partition when the slope condition (ii) in Definition \ref{def: base partition} for a sub-interval $[x_i, x_{i + 1}]$ is not satisfied. As a result, a base partition is not unique in general.

%% file: Sections/relaxations.tex
\label{MILP Relaxations Section}

In this section we present a method for constructing piecewise polyhedral  relaxations of univariate trigonometric terms \eqref{univariate term} and bilinear terms \eqref{bilinear term} in $\mathcal F$. Once these nonlinearities have been relaxed, the resulting MILP can be solved using a standard MILP solver. 

\subsection{Trigonometric Terms} \label{Trigonometric Terms Subsection}

Consider a univariate, trigonometric function $f \colon [x^L, x^U] \mapsto \mathbb{R}$, where $f$ is differentiable and bounded over $[x^L, x^U] \subset \mathbb R$. We would like to construct a MILP relaxation of the constraint $y = f(x)$ with $x \in [x^L, x^U]$, which can then be used for the constraints \eqref{univariate term} in $\mathcal F$. Define an admissible partition $p = (x_0, x_1, \hdots, x_m)$ of $[x^L, x^U]$ for $f$ where $p$ is a valid refinement of a base partition $p^0$ or is $p^0$ itself. Let $[x_i, x_{i + 1}]$ be a sub-interval defined by $p$. Define the following:
\begin{equation} \label{tangent 1}
    h_i(x) = f(x_i) + f'(x_i) \cdot (x - x_i)
\end{equation}
\begin{equation} \label{tangent 2}
    h_{i + 1}(x) = f(x_{i + 1}) + f'(x_{i + 1}) \cdot (x - x_{i + 1})
\end{equation}
\begin{equation} \label{nonsecant}
    t_{i + 1}(x) =
    \begin{cases}
        \max\{h_i(x), \, h_{i + 1}(x) \}, \quad \text{if $f$ is convex in $[x_i, x_{i + 1}]$} \\
        \min\{h_i(x), \, h_{i + 1}(x) \}, \quad \text{if $f$ is concave in $[x_i, x_{i + 1}]$}
    \end{cases}
\end{equation}
\begin{equation} \label{secant}
    s_{i + 1}(x) = f(x_i) + \frac{f(x_{i + 1}) - f(x_i)}{x_{i + 1} - x_i} \cdot (x - x_i)
\end{equation}
Equations \eqref{tangent 1} and \eqref{tangent 2} define tangent lines at $f(x_i)$ and $f(x_{i + 1})$, respectively. When $f$ is convex (resp. concave) the tangent lines lie below (resp. above) the function and so a simple underestimate (resp. overestimate) of $f$ over the sub-interval $[x_i, x_{i + 1}]$ is obtained by taking the maximum (resp. minimum) of these two tangents, corresponding to \eqref{nonsecant}. Equation \eqref{secant} defines a secant line from $f(x_i)$ to $f(x_{i + 1})$, which is above (resp. below) $f$ when it is convex (resp. concave) in the sub-interval $[x_i, x_{i + 1}]$. For each partition point $x_i \in \mathcal{P}(p)$, let $v_i = f(x_i)$ be the corresponding point on the curve. For each sub-interval $[x_i, x_{i + 1}]$, let $v_{i, i + 1}$ be the point of intersection for the tangent lines defined by \eqref{tangent 1} and \eqref{tangent 2}.

\begin{figure}    
    \centering
    \subfloat[Using base partition.]{
    \includegraphics[width = 0.75\textwidth]{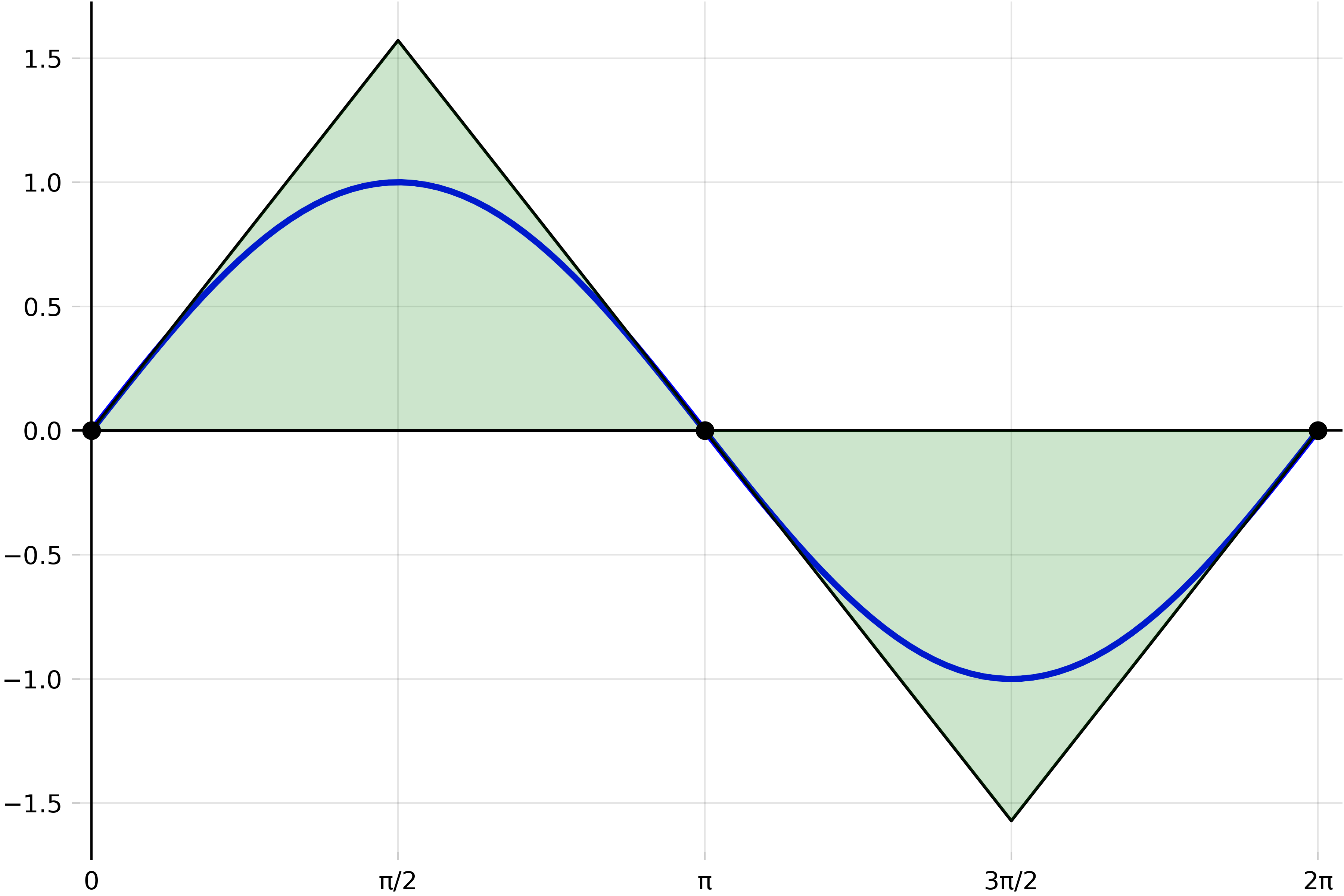}
    }
    \hfill
    \subfloat[Using valid refinement.]{
    \includegraphics[width = 0.75\textwidth]{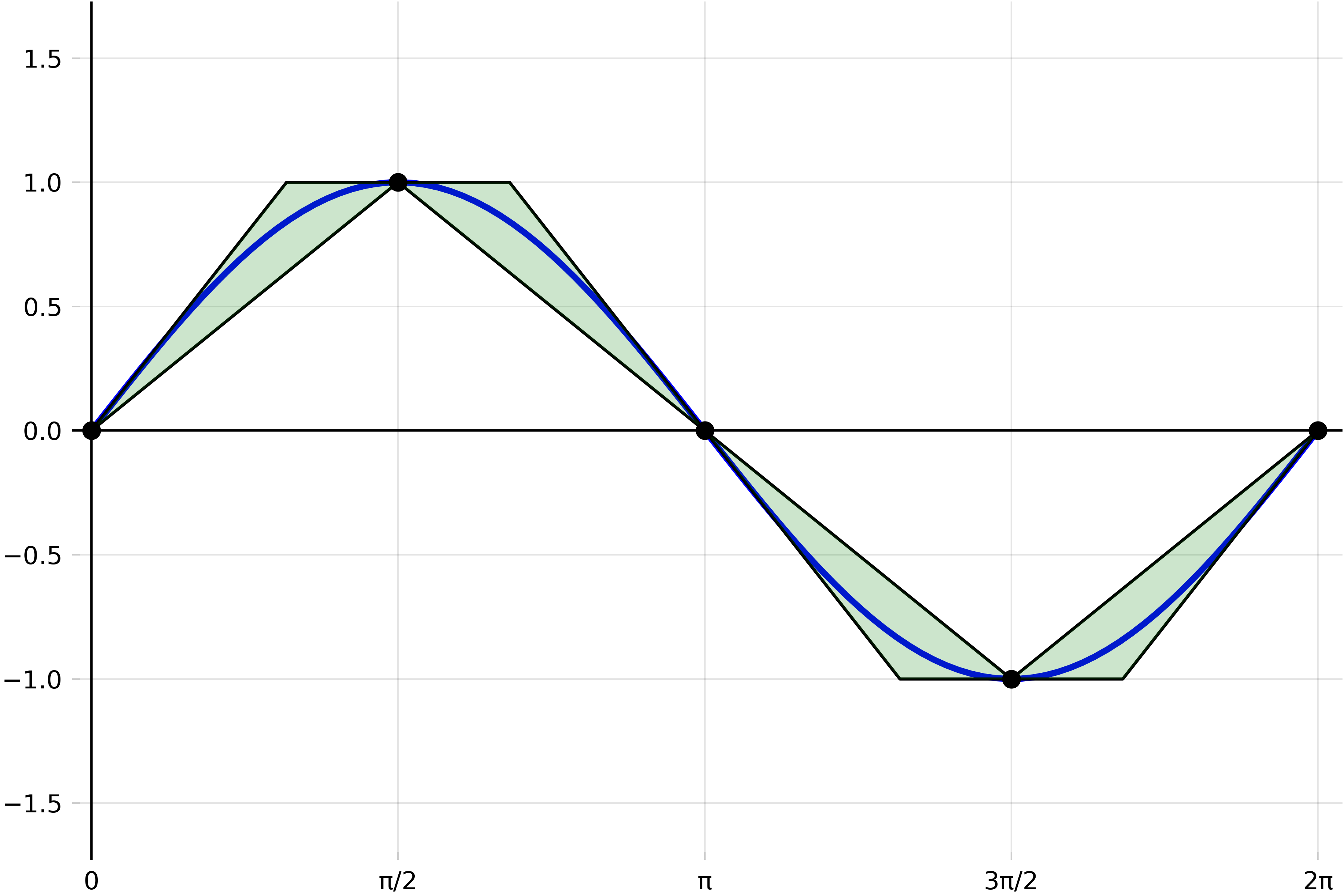}
    }
    \caption{Example of triangles formed by overestimates, underestimates, and secant lines of $f(x) = \sin x$ with $x \in [0, 2\pi]$ and admissible partition $p$. The points on the curve corresponding to the partition points are shown with black markers. The green triangles correspond to the points $v_i$, $v_{i, i + 1}$, and $v_{i + 1}$ for each sub-interval $[x_i, x_{i + 1}]$ of $p$. (a) Triangles formed when using base partition $p^0 = (0, \pi, 2\pi)$. (b) Triangles formed when using the admissible partition $p = (0, \frac{\pi}{2}, \pi, \frac{3\pi}{2}, 2\pi)$, which is a valid refinement of the base partition $p^0$. }
    \label{fig: sin example}
\end{figure}

We now make a few key observations. For any sub-interval $[x_i, x_{i + 1}]$ defined by $p$, the points $v_i$, $v_{i, i + 1}$, and $v_{i + 1}$ form a triangle which contains the curve $f$ over this sub-interval. This can be seen by noting the partition $p$ is admissible and so $f$ is either convex or concave in each sub-interval. We are also guaranteed the intersection point $v_{i, i + 1}$ exists for each sub-interval as $p$ satisfies the slope condition in Definitions \ref{def: base partition} and \ref{def: valid refinement}. The number of triangles formed by $p$ is $m = |\mathcal{P}(p)|- 1$. If $q$ is a valid refinement of $p$, the number of triangles will increase as $\mathcal{P}(p) \subset \mathcal{P}(q)$. Additionally, by the nature of the construction of the triangles, the triangles themselves will decrease in size and approach the curve itself. This can be seen by the simple example $f(x) = \sin x$ with $x \in [0, 2\pi]$ (see Figure \ref{fig: sin example}). We have a unique base partition of $p^0 = (0, \pi, 2\pi)$. Define $p = (0, \frac{\pi}{2}, \pi, \frac{3\pi}{2}, 2\pi)$, which is simply the result of bisecting each sub-interval of $p^0$. We can see $p$ is a valid refinement of $p^0$ and so once again we can form triangles containing the curve over each sub-interval. This is shown in Figure \ref{fig: sin example}. It can be clearly seen the triangles approach the curve itself as the number of partition points increases with valid refinements. 

The main idea behind the MILP relaxation is then as follows: For each sub-interval $[x_i, x_{i + 1}]$ of $p$, any point in the triangle defined by the points $v_i$, $v_{i, i + 1}$, and $v_{i + 1}$ provides a relaxation of the constraint $y = f(x)$ when $x \in [x_i, x_{i + 1}]$. Therefore, the disjunctive union of triangles formed by $p$ using \eqref{tangent 1}-\eqref{secant} provides a relaxation of $y = f(x)$ with $x \in [x^L, x^U]$. As valid refinements are iteratively constructed, this relaxation is tightened and will ultimately converge to the original function in the limit.

\subsubsection{Incremental Formulation for Trigonometric Terms} \label{Trigonometric Incremental Formulation Subsubsection}

We are now prepared to construct the MILP relaxation for the constraint $y = f(x)$ with $x \in [x^L, x^U] \subset \mathbb R$ where $f$ is a univariate, trigonometric function that is differentiable and bounded over $[x^L, x^U]$. Let $p = (x_0, x_1, \hdots, x_m)$ be a valid refinement of a base partition $p^0$ of $[x^L, x^U]$ for $f$ or be $p^0$ itself. As before, let $v_i = f(x_i)$ and let $v_{i, i + 1}$ be the intersection of the two tangent lines defined by \eqref{tangent 1} and \eqref{tangent 2} for each sub-interval $[x_i, x_{i + 1}]$. Let $v_{i}^x$ and $v^y_{i}$ be the $x$- and $y$-coordinate of $v_i$ and similarly for the vertex $v_{i, i + 1}$. Define binary variables $u_i$ for $i = 1, \hdots, m - 1$ and non-negative continuous variables $\delta^i_1$ and $\delta^i_2$ for $i = 1, \hdots, m$. The disjunctive union of triangles formed by $v_{i}$, $v_{i, i + 1}$, and $v_{i + 1}$ for each sub-interval $[x_i, x_{i + 1}]$ provides a relaxation of $y = f(x)$ with $x \in [x^L, x^U]$ and can be expressed using a standard incremental formulation \cite{yildiz2013incremental,sundar2021sequence} as
\begin{subequations} \label{MILP relaxation : trig}
    \begin{equation} \label{x coordinate MILP relaxation}
        x = v_0^x + \sum_{i = 1}^m \left\{\delta^i_1(v^x_{i - 1, i} - v^x_{i - 1}) + \delta^i_2(v^x_i - v^x_{i - 1}) \right\}
    \end{equation}
    \begin{equation} \label{y coordinate MILP relaxation}
        y = v_0^y + \sum_{i = 1}^m \left\{\delta^i_1(v^y_{i - 1, i} - v^y_{i - 1}) + \delta^i_2(v^y_i - v^y_{i - 1}) \right\}
    \end{equation}
    \begin{equation}
        \delta^1_1 + \delta^1_2 \leq 1
    \end{equation}
    \begin{equation} \label{linking constraint MILP relaxation}
        \delta^i_1 + \delta^i_2 \leq u_{i - 1} \leq \delta^{i - 1}_2 \quad \forall i \in \{2, \hdots, m\}
    \end{equation}
    \begin{equation} \label{delta bounds}
        0 \leq \delta^i_1, \, \delta^i_2 \leq 1 \quad \forall i \in \{1, \hdots, m\}
    \end{equation}
    \begin{equation}
        u_i \in \{0, 1\} \quad \forall i \in \{1, \hdots, m - 1\}
    \end{equation}
\end{subequations}
Though there are many ways to formulate the disjunctive union of triangles with theoretical properties similar to that of the incremental formulation as shown in \cite{yildiz2013incremental,vielma2015mixed}, we utilize the incremental formulation given its computational superiority to other ways of formulating the disjunctive union \cite{huchette2022nonconvex}. Nevertheless, we remark that comparing the different formulations for relaxing trigonometric functions is an interesting problem in its own right and we delegate this to future work.

The relaxation \eqref{MILP relaxation : trig} can be understood as follows. The terms $v_0^x$ and $v_0^y$ represent a starting point for $x$ and $y$. For partition $p = (x_0, x_1, \hdots, x_m)$, we will refer to the triangle corresponding to the sub-interval $[x_i, x_{i + 1}]$ as the $(i + 1)$-th triangle and so the triangles are ordered. Let $(x, y)$ denote the relaxation value of $(x, f(x))$. The binary variable $u_i$ takes value 1 if the $i$-th triangle must be passed to reach the triangle containing the point $(x, y)$ and takes value 0 otherwise. It can be shown \cite{yildiz2013incremental} that \eqref{MILP relaxation : trig} has the ordering property $u_1 \geq u_2 \geq \hdots \geq u_{m - 1}$. Therefore, if $u_{i^*} = 1$ then we must have $u_{i} = 1$ for all $i < i^*$, which simply states we must pass the triangles in order before reaching the triangle containing $(x, y)$. The variables $\delta^i_1$ and $\delta^i_2$ for the $i$-th triangle have two purposes which are most easily understood through example. Suppose $(x, y)$ lies in the first triangle. In this case, we must have $u_i = 0$ for all $i \in \{1, \hdots, m-1\}$, which implies $\delta^i_1 = 0$ and $\delta^i_2 = 0$ for all $i \in \{2, \hdots, m\}$ by \eqref{linking constraint MILP relaxation} and \eqref{delta bounds}. Relaxation \eqref{MILP relaxation : trig} then becomes
\begin{subequations} \label{first triangle example}
    \begin{equation}
        x = v_0^x + \delta^1_1(v_{0, 1}^x - v^x_0) + \delta^1_2(v_1^x - v_0^x)
    \end{equation}
    \begin{equation}
        y = v_0^y + \delta^1_1(v_{0, 1}^y - v^y_0) + \delta^1_2(v_1^y - v_0^y)
    \end{equation}
    \begin{equation}
        \delta^1_1 + \delta^1_2 \leq 1
    \end{equation}
    \begin{equation}
        0 \leq \delta^1_1, \delta^1_2 \leq 1
    \end{equation}  
\end{subequations}
From \eqref{first triangle example}, we see $\delta^1_1$ and $\delta^1_2$ are used to capture any $(x, y)$ in the first triangle defined by $v_0$, $v_{0, 1}$, and $v_1$. In the case where $(x, y)$ lies in the $i$-th triangle, variables $\delta^i_1$ and $\delta^i_2$ hold similar meaning. Next, suppose $(x, y)$ lies in the second triangle. In this case, we must have $u_1 = 1$ and $u_i = 0$ for all $i \in \{2, \hdots, m - 1\}$. From \eqref{linking constraint MILP relaxation} we then must have $\delta^1_1 = 0$ and $\delta^1_2 = 1$. Similarly, we must also have $\delta^i_1 = 0$ and $\delta^i_2 = 0$ for all $i \in \{3, \hdots, m\}$. Relaxation \eqref{MILP relaxation : trig} then becomes
\begin{subequations} \label{second triangle example}
    \begin{equation}
        x = v_0^x + (v_1^x - v_0^x) + \delta^2_1(v_{1, 2}^x - v^x_1) + \delta^2_2(v_2^x - v_1^x)
    \end{equation}
    \begin{equation}
        y = v_0^y + (v_1^y - v_0^y) + \delta^2_1(v_{1, 2}^y - v^y_1) + \delta^2_2(v_2^y - v_1^y)
    \end{equation}
    \begin{equation}
        \delta^2_1 + \delta^2_2 \leq 1
    \end{equation}
    \begin{equation}
        0 \leq \delta^2_1, \delta^2_2 \leq 1
    \end{equation}
\end{subequations}
We now note that $\delta^1_2 = 1$ is the coefficient of $(v_1^x - v_0^x)$ and $(v_1^y-v_0^y)$, which corresponds to the secant line connecting $v_0$ and $v_1$. Therefore, $\delta^1_2$ simply shifts the starting point from $v_0$ to $v_1$ or, equivalently, shifts the search for $(x,y)$ from the first triangle to the second triangle. After this shift, we again see $\delta^2_1$ and $\delta^2_2$ are used to capture any $(x, y)$ in the second triangle, as expected. The variable $\delta^1_1$ is zero as we do not need it to shift from $v_0$ to $v_1$ by construction. In the case where $(x, y)$ lies in the $i^*$-th triangle, variables $\delta^i_1$ and $\delta^i_2$ for all $i < i^*$ hold similar meaning. 

\subsubsection{Convergence Guarantee}

We refer the reader to \cite{sundar2021sequence} for a proof the previously described MILP relaxation converges to the original function as all sub-intervals are refined in the limit. We note this convergences relies on the refinement strategy used on the partition. This is addressed in Section \ref{Refinement Scheme Section}.

\subsection{Bilinear Terms}
This section presents an incremental formulation to obtain a MILP relaxation of the constraint $z = xy$ where $x \in [x^L, x^U] \subset \mathbb R$ and $y \in [y^L, y^U] \subset \mathbb R$. Such a relaxation can then be used for the bilinear terms \eqref{bilinear term} in $\mathcal{F}$ found after factoring. The relaxation we present has a similar construction as the relaxation for trigonometric terms. Using MILP relaxations for bilinear terms is not new and has been dealt with extensively in the literature \cite{bergamini2005logic,Nagarajan2019}. Typically the MILP relaxations use either a big-M reformulation or a convex hull representation of a disjunctive union \cite{Nagarajan2019}. We instead use an incremental formulation similar to the relaxation used for trigonometric terms. Furthermore, we assume that only one of the variables involved in the bilinear term, i.e., either $x$ or $y$ is partitioned. To the best of our knowledge, the following incremental formulation is not explicitly presented in the literature and so, for the sake of completeness we present it here. 

\begin{figure}
    \centering
    \includegraphics[width = 0.8\textwidth]{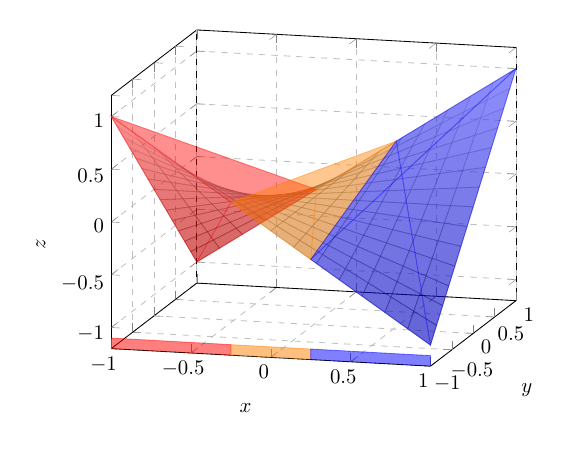}
    \caption{Disjunctive union of three tetrahedrons containing $z = xy$ over the domain $[-1, 1] \times [-1, 1]$. The variable $x$ has been partitioned using $p_x = (-1, -0.25, 0.25, 1)$, shown by the colored bars along the $x$-axis with tetrahedrons being colored accordingly.}
    \label{fig: disjunctive union bilinear}
\end{figure}

Consider the standard convex envelope \cite{mccormick1976computability} commonly used for relaxing bilinear terms. 

\begin{subequations} \label{McCormick}
    \begin{equation}
        z = w 
    \end{equation}
    \begin{equation}
        w \geq x^L y + x y^L - x^L y^L
    \end{equation}
    \begin{equation}
        w \geq x^U y + x y^U - x^Uy^U
    \end{equation}
    \begin{equation}
        w \leq x^U y + x y^L - x^U y^L
    \end{equation}
    \begin{equation}
        w \leq x y^U + x^L y - x^L y^U
    \end{equation}
\end{subequations}
The convex envelope defined by \eqref{McCormick} is precisely the convex hull of $z = xy$ over the domain $[x^L, x^U] \times [y^L, y^U]$.
%
\begin{figure}
    \centering
    \includegraphics[width = 0.8\textwidth]{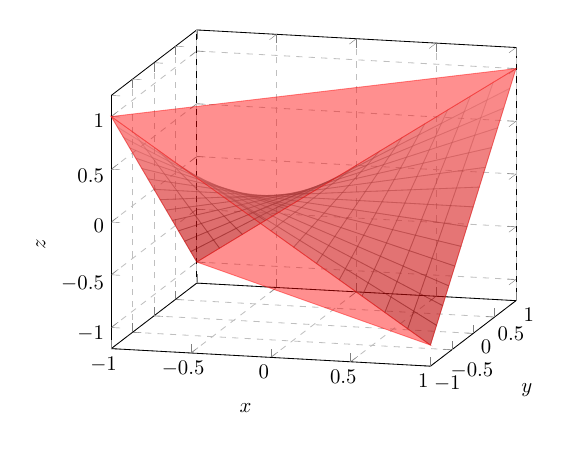}
    \caption{Convex hull of $z = xy$ with domain $[-1, 1] \times [-1, 1]$.  }
    \label{fig: basic mccormick}
\end{figure}
This convex hull is the tightest tetrahedron containing the graph of the bilinear term (see Figure \ref{fig: basic mccormick}). Consider a partition $p_x = (x_0, x_1, \hdots, x_m)$ of $[x^L, x^U]$. Let $R_{i + 1} = [x_i, x_{i + 1}] \times [y^L, y^U]$ be the $(i + 1)$-th rectangle defined by sub-interval $[x_i, x_{i + 1}]$ of $p_x$. We can describe the surface over $R_{i + 1}$ by
\begin{subequations} \label{bilinear over subdomain}
    \begin{equation}
        z = xy
    \end{equation}
    \begin{equation}
        x_i \leq x \leq x_{i + 1}
    \end{equation}
    \begin{equation}
        y^L \leq y \leq y^U
    \end{equation}
\end{subequations}
We see \eqref{bilinear over subdomain} is once again a bilinear term with box constraints on $x$ and $y$, so we can again construct a convex envelope of \eqref{bilinear over subdomain} to get a tetrahedron containing the graph $z = xy$ over sub-domain $R_{i + 1}$. Figure \ref{fig: disjunctive union bilinear} shows an example of this procedure.

The main idea behind the MILP relaxation is then as follows: For each sub-domain $R_{i + 1}$ defined by partition $p_x$, construct the tightest tetrahedron (i.e., convex envelope) containing $z = xy$ over $R_{i + 1}$. Any point in the corresponding tetrahedron provides a relaxation of $z = xy$ with $(x, y) \in R_{i + 1}$. Therefore, the disjunctive union of tetrahedrons formed by $p_x$ provides a relaxation of $z = xy$ over the original domain $[x^L, x^U] \times [y^L, y^U]$.

\subsubsection{Incremental Formulation for Bilinear Terms}

We are now prepared to construct the MILP relaxation for the bilinear term $z = xy$ with $x \in [x^L, x^U] \subset \mathbb R$ and $y \in [y^L, y^U] \subset \mathbb R$. Let $p_x = (x_0, x_1, \hdots, x_m)$ be a partition of $[x^L, x^U]$. Define binary variables $u_i$ for $i = 1, \hdots, m - 1$ and non-negative continuous variables $\delta^i_1$, $\delta^i_2$, and $\delta^i_3$ for $i = 1, \hdots, m$. The MILP relaxation is the disjunctive union of tetrahedrons constructed from $p_x$ and can be expressed \cite{yildiz2013incremental} as

\begin{subequations} \label{MILP relaxation : bilinear}
    \begin{equation} \label{x coordinate bilinear}
        x = x^L + \sum_{i = 1}^{m} \left[\delta^i_2(x_{i} - x_{i - 1}) + \delta^i_3(x_{i} - x_{i - 1}) \right]
    \end{equation}
    \begin{equation} \label{y coordinate bilinear}
        y = y^L + \sum_{i = 1}^{m}\left[\delta^i_1(y^U-y^L) + \delta^i_2(y^U - y^L) \right]
    \end{equation}
    \begin{equation} \label{z coordinate bilinear}
        z = x^Ly^L + \sum_{i = 1}^{m} 
        \begin{pmatrix}
        \delta^i_1 \\ \delta_2^i \\ \delta_3^i
        \end{pmatrix}^T 
        \begin{pmatrix}
            x_{i - 1}y^U - x_{i - 1}y^L \\
            x_{i}y^U - x_{i - 1} y^L \\
            x_{i}y^L - x_{i - 1} y^L
        \end{pmatrix}
    \end{equation}
    \begin{equation} \label{bilinear linking}
        \delta^i_1 + \delta^i_2 + \delta^i_3 \leq u_{i - 1} \leq \delta^{i -1}_3 \quad \forall i \in \{2, \hdots, m\}
    \end{equation}
    \begin{equation}
        \delta^1_1 + \delta^1_2 + \delta^1_2 \leq 1
    \end{equation}
    \begin{equation} \label{bilinear delta bounds}
        0 \leq \delta^i_1, \, \delta^i_2 \, \delta^i_3 \leq 1, \quad \forall i \in \{1, \hdots, m\}
    \end{equation}
    \begin{equation}
        u_i \in \{0, 1\}, \quad \forall i \in \{1, \hdots, m - 1\}
    \end{equation}
\end{subequations}

Relaxation \eqref{MILP relaxation : bilinear} can be understood as follows. For the given partition $p_x$, we will refer to the tetrahedron corresponding to sub-domain $R_{i + 1}$ as the $(i + 1)$-th tetrahedron and so the tetrahedrons are ordered. Note the $(i+1)$-th tetrahedron over $R_{i + 1}$ has the four extreme points
\begin{subequations} \label{Extreme points of tetrahedron}
    \begin{align}
        &v_{i + 1}^0 = (x_i, y^L, x_iy^L) \\
        &v_{i + 1}^1 = (x_i, y^U, x_iy^U)\\
        &v_{i + 1}^2 = (x_{i + 1}, y^U, x_{i + 1}y^U) \\
        &v_{i + 1}^3 = (x_{i + 1}, y^L, x_{i + 1}y^L)
    \end{align}
\end{subequations}
The terms $x^L$, $y^L$, and $x^Ly^L$ represent starting points for $x$, $y$, and $z$, respectively. In particular, we see $(x^L, y^L, x^Ly^L)$ is the extreme point $v_1^0$. Let $(x, y, z)$ denote the relaxation value of $(x, y, xy)$. The binary variable $u_i$ takes value 1 if the $i$-th tetrahedron must be passed to reach the tetrahedron containing the point $(x, y, z)$ and takes value 0 otherwise. It can be shown \cite{yildiz2013incremental} that \eqref{MILP relaxation : bilinear} has the ordering property $u_1 \geq u_2 \geq \hdots \geq u_{m - 1}$. Therefore, if $u_{i^*} = 1$ then we must have $u_i = 1$ for all $i < i^*$, which simply states we must pass the tetrahedrons in order before reaching the tetrahedron containing $(x, y, z)$. The variables $\delta^i_1$, $\delta^i_2$, and $\delta^i_3$ for the $i$-th tetrahedron have similar meaning as the $\delta$ variables in the trigonometric case.

\subsubsection{Convergence Guarantee}

It can be shown \cite{androulakis1995alphabb} the largest gap between the convex envelope \eqref{McCormick} over a sub-domain $R_{i + 1} = [x_i, x_{i + 1}] \times [y^L, y^U]$ is at most
\begin{equation}
    \max_{(x, y) \in R_{i + 1}} |w_{i + 1} - xy| =  \frac{(x_{i + 1} - x_i)(y^U - y^L)}{4}
\end{equation}
where $w_{i + 1}$ denotes the relaxation value. This bound holds for each sub-domain and so as all sub-domains are further refined the MILP relaxation \eqref{MILP relaxation : bilinear} converges to the original bilinear function in the limit. We note this convergence relies on all sub-domains being further refined. This will be addressed in Section \ref{Refinement Scheme Section}.

%% file: Sections/tighten.tex
\label{Partitions Section}

In this section we focus on the schemes used for partitioning selected variable domains which in turn are used to construct the MILP relaxations previously described. We first show how to reduce the total number of variables needed across MILP relaxations of different nonlinear functions by sharing partitions between functions when multiple functions use the same variable. This is especially important when using trigonometric functions, as many applications will have several functions sharing the same variable (for example, $\sin x$ and $\cos x$ may appear in pairs when looking at signals or geometric constraints).  We next introduce several refinement schemes to refine a given partition after finding an optimal solution to the MILP relaxation of $\mathcal{F}$ to further tighten the relaxation. The refinement scheme chosen will directly impact the convergence rate of the overall algorithm and so we present three options - bisection, direct, and non-uniform. Our list of refinement schemes is not exhaustive and it is entirely possible to use alternative refinement schemes not discussed in this paper. When refining partitions for the MILP relaxation of $\mathcal{F}$, we may either choose to refine all available partitions or elect to refine a subset of partitions subject to some criteria (under certain conditions so convergence to a globally optimal solution is ensured). We briefly discuss these two possibilities and introduce a simple criteria for refining a subset of all partitions as a motivating example. 

\subsection{Sharing Partitions}

For ease of discussion, we will consider two univariate trigonometric functions $f_1(x)$ and $f_2(x)$ which share the same variable $x$ and have been relaxed using the MILP relaxation presented in Section \ref{MILP Relaxations Section}. The same procedure presented applies to bilinear terms that have been relaxed as well and may be extended to more than two functions. The variable $x$ has a box constraint $[x^L, x^U]$.

Let $p_1^0$ and $p_2^0$ be base partitions of $[x^L, x^U]$ for $f_1$ and $f_2$, respectively. Rather than constructing two separate partitions, we can use information from both base partitions to construct a single partition that can be used for both functions simultaneously in the MILP relaxation. Define a partition as the union of points in $p_1^0$ and $p_2^0$. If the slope condition in Definition \ref{def: base partition} is not initially satisfied (applied to both functions), add additional partition points so that the slope condition is satisfied and the resulting partition, $p^0$, is of minimum size. The resulting partition is in a sense a base partition of the collection of functions. We note that for any sub-interval $[x_i, x_{i + 1}]$ defined by $p^0$, the functions $f_1$ and $f_2$ are either convex or concave over the sub-interval (the functions need not have the same convexity). This allows the use of the incremental formulation presented in Section \ref{MILP Relaxations Section} and so we can refine $p^0$ to tighten all parts of the formulation containing $f_1$ and $f_2$ rather than refining two separate partitions. In doing so, we reduce the number of binary variables needed for $f_1$ and $f_2$ in the incremental formulation by half (or a factor of $M$ when $M$ functions are considered instead of two). As shown later in the computational results, this will significantly reduce computation time. 

\subsubsection{Sharing Partitions Example}

\begin{figure}
    \centering
    \includegraphics[width = 0.95\textwidth]{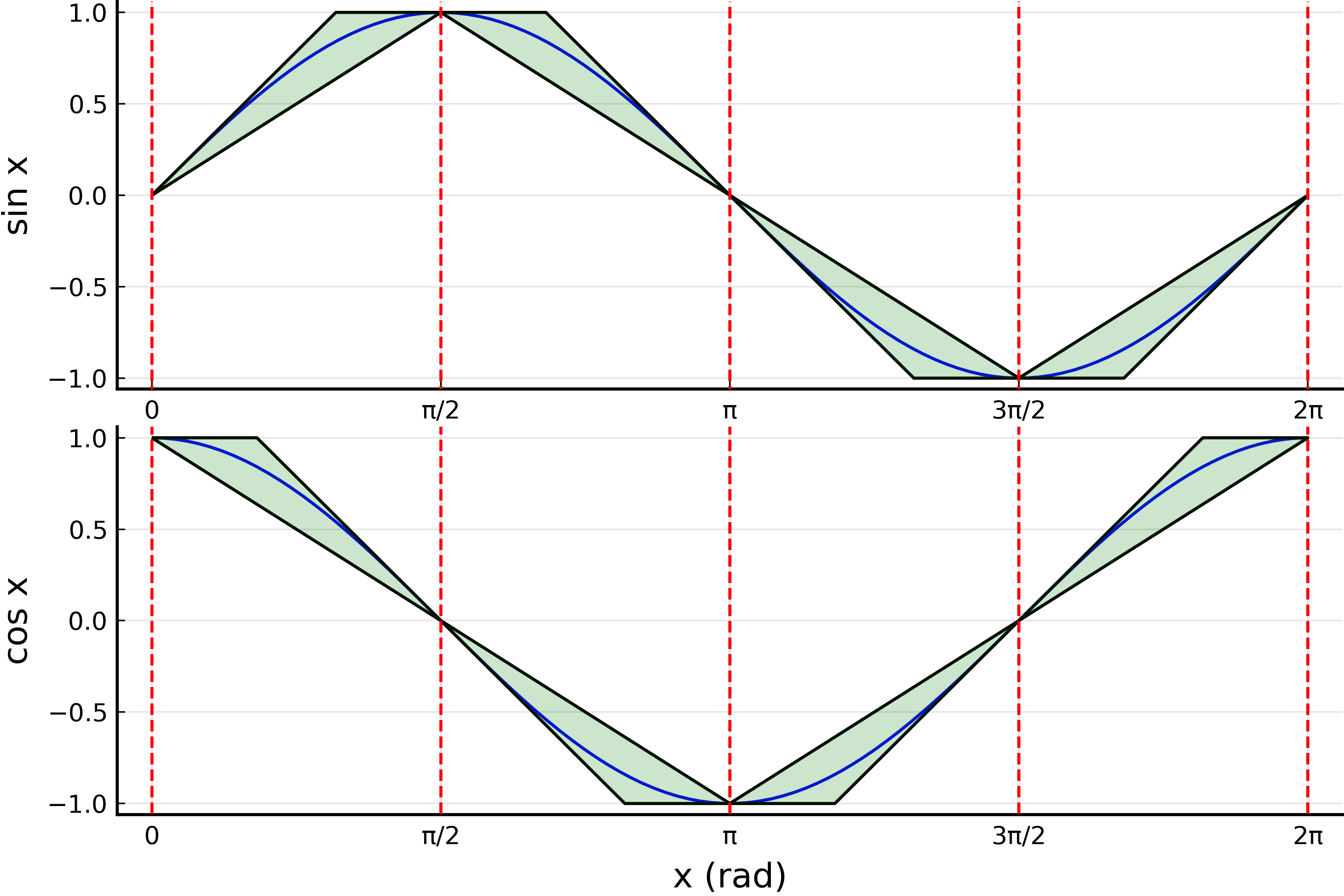}
    \caption{Example of using a shared partition, $p^0 = (0, \frac{\pi}{2}, \pi, \frac{3\pi}{2}, 2\pi)$, for $f_1(x) = \sin x$ and $f_2(x) = \cos x$ over the closed interval $[0, 2\pi]$. It can be seen the $i$-th triangle for $f_1$ and the $i$-th triangle for $f_2$ are defined over the same sub-interval, so they can be linked by sharing the same binary variables in the MILP relaxations. The binary variables then indicate which sub-interval the solution is in, rather than which triangle.}
    \label{fig:sharing partitions example}
\end{figure}

As an example, consider the two functions $f_1(x) = \sin x$ and $f_2(x) = \cos x$ where $x \in [0,2\pi]$. The (unique) base partitions for $f_1$ and $f_2$ are $p_1^0 = (0, \pi, 2\pi )$ and $p_2^0 = (0, \frac{\pi}{2}, \frac{3\pi}{2}, 2\pi)$, respectively. We can then define a new partition $p^0 = (0, \frac{\pi}{2}, \pi, \frac{3\pi}{2}, 2\pi)$. We see the slope condition in Definition \ref{def: base partition} is satisfied for both $f_1$ and $f_2$ with $p^0$ and so $p^0$ is effectively a base partition for the set of functions. We can then use relaxation \eqref{MILP relaxation : trig} for both $f_1$ and $f_2$ using the same partition. The $i$-th triangle for $f_1$ and the $i$-th triangle for $f_2$ are linked since they are defined over the same sub-interval. This is readily seen in Figure \ref{fig:sharing partitions example}. With this in mind, the binary variables in relaxation \eqref{MILP relaxation : trig} for both functions may be shared.

\subsection{Refinement Schemes} \label{Refinement Scheme Section}

We first introduce the following definition. 

\begin{definition} \label{Refinement Scheme}
    Given an admissible partition $p$ of $[x^L, x^U]$ for a function $f$, the procedure used to generate a valid refinement of $p$ is referred to as a \textit{refinement scheme}.
\end{definition}
We note that a refinement scheme is for a \textit{single} partition and it is entirely possible to use various refinement schemes depending on the partition considered. As a simple example, consider the base partition $p^0 = (0, \pi, 2\pi)$ of $[0, 2\pi]$ for $f(x) = \sin x$ (see Figure \ref{fig: sin example}). We can construct the valid refinement $q = (0, \frac{\pi}{2}, \pi, \frac{3\pi}{2}, 2\pi)$, where the refinement scheme used is: For each sub-interval defined by $p^0$ we add a point to bisect that sub-interval. We note that in this particular example the refinement scheme used is indifferent of any known solution to the MILP relaxation. This may be considered inefficient, as some sub-intervals are further refined despite them potentially not being useful to the solution of the original problem. We soon consider more targeted refinement schemes by using information of the optimal solution to the MILP relaxation.

For the remainder of this section, we use the following notation. Suppose $\mathcal{F}$ has been relaxed using the previously described MILP relaxations. For a given term (trigonometric or bilinear) that has been relaxed, let $x$ be the variable that has been partitioned with $p = (x_0, \hdots, x_m)$ and let $x^*$ denote its value in the optimal solution of the MILP. We note that only a single partition $p$ needs to be associated with $x$ because multiple relaxed terms using $x$ may share the same partition as previously discussed. We will assume $x^*$ is strictly in the interior of a sub-interval $[x_i, x_{i + 1}]$ defined by $p$. In the case where $x^* = x_i$ for some $x_i \in \mathcal{P}(p)$, we do not refine partition $p$ \textit{for that particular iteration of the overall algorithm}. In a later iteration we may have $x^*$ strictly in the interior of a sub-interval defined by $p$ and so we may refine $p$ for that iteration. We now discuss three proposed refinement schemes --- bisection, direct, and non-uniform. 

\subsubsection{Bisection Refinement Scheme}

The first refinement scheme we consider is bisection. Given $x^* \in (x_i, x_{i + 1})$, the bisection refinement scheme is to simply add a new partition point $x' = \frac{1}{2}(x_i + x_{i + 1})$ to $p$, giving the valid refinement $q = (x_0, \hdots, x_i, x', x_{i + 1}, \hdots, x_m)$. 

\begin{figure}
    \centering
    \includegraphics[width=0.75\textwidth]{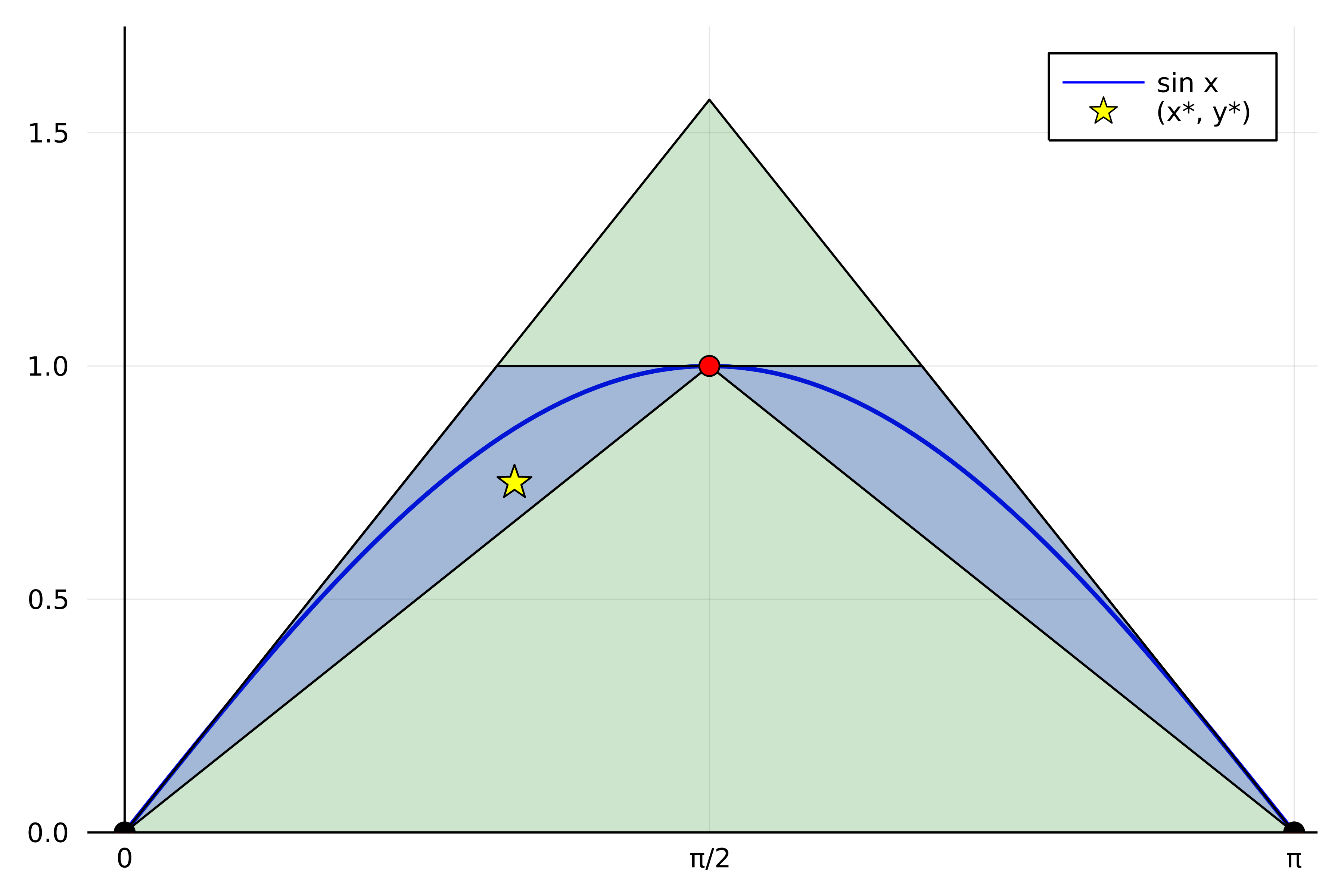}
    \caption{MILP relaxation of $f(x) = \sin x$ over $[0, \pi]$ using $p = (0, \pi)$ (green) and the valid refinement $q = (0, \frac{\pi}{2}, \pi)$ (blue), which was constructed using a bisection refinement scheme. The optimal solution (yellow star) of the resulting MILP, $(x^*, y^*)$, when using $p$ may still lie in the triangles constructed using $q$. If this happens, $(x^*, y^*)$ remains a possible solution to the MILP when using $q$ and so the lower bound may not improve for that particular iteration of the algorithm.}
    \label{fig: bisection visual}
\end{figure}

It is important to note that while this refinement scheme is simple, it may not \textit{always} improve the lower bound to $\mathcal{F}$. This can be seen by considering $f(x) = \sin x$ with domain $[0, \pi]$. We may first start with the base partition $p = (0, \pi)$ and use relaxation \eqref{MILP relaxation : trig} for $f$ using $p$, where we introduce a variable $y$ for the relaxation value of $f(x)$. This will result in a single triangle containing the curve (green triangle in Figure \ref{fig: bisection visual}). After relaxing $f$ (and other constraints), we may solve the resulting MILP, yielding values $(x^*, y^*)$. It is possible for $(x^*, y^*)$ to lie in the \textit{interior} of the triangle defined by $p$. Since $x^*$ is in the interior of $[0, \pi]$, we may refine $p$ using bisection to get the valid refinement $q = (0, \frac{\pi}{2}, \pi)$. Using relaxation \eqref{MILP relaxation : trig} with $q$ results in two triangles containing $f$ (blue triangles in Figure \ref{fig: bisection visual}). Because $(x^*, y^*)$ was in the interior of the original triangle, it is possible this point is also in one of the two newly created triangles. This case is shown in Figure \ref{fig: bisection visual}. In the event this happens, solving the new MILP relaxation \textit{may not result in a tighter lower bound}, though this will largely depend on any other constraints present in the problem. We also note a similar situation can occur when $(x^*, y^*)$ is on an \textit{edge} of a triangle, which can be seen by noting the overlapping edges in Figure \ref{fig: bisection visual}.

It is important to note that while the lower bound may not improve for a particular iteration using bisection, the MILP's solution will still converge to the optimal solution of $\mathcal{F}$ in the limit. Once again consider the example shown in Figure \ref{fig: bisection visual} and suppose $(x^*, y^*)$ remains optimal when relaxing using $q$. We then apply the bisection refinement scheme to $q$, resulting in $r = (0, \frac{\pi}{4}, \frac{\pi}{2}, \pi)$ and so two triangles are created over $[0, \frac{\pi}{2}]$. Visually, it should be clear from Figure \ref{fig: bisection visual} that the two triangles will likely not contain $(x^*, y^*)$. When this happens, the lower bound will continue to improve. This same idea holds in the case of bilinear terms as well. We are guaranteed this will always eventually happen because the MILP relaxations approach the original functions in the limit. As such, any relaxation solution that is not on the original function will eventually become infeasible as partitions are further refined. 

\subsubsection{Direct Refinement Scheme}

We next consider the direct refinement scheme, which is another natural choice. Given $x^* \in (x_i, x_{i + 1})$, the direct refinement scheme is to simply add $x^*$ to $p$, giving the valid refinement $q = (x_0, \hdots, x_i, x^*, x_{i + 1}, \hdots, x_m)$. 

\begin{figure}
    \centering
    \includegraphics[width = 0.75\textwidth]{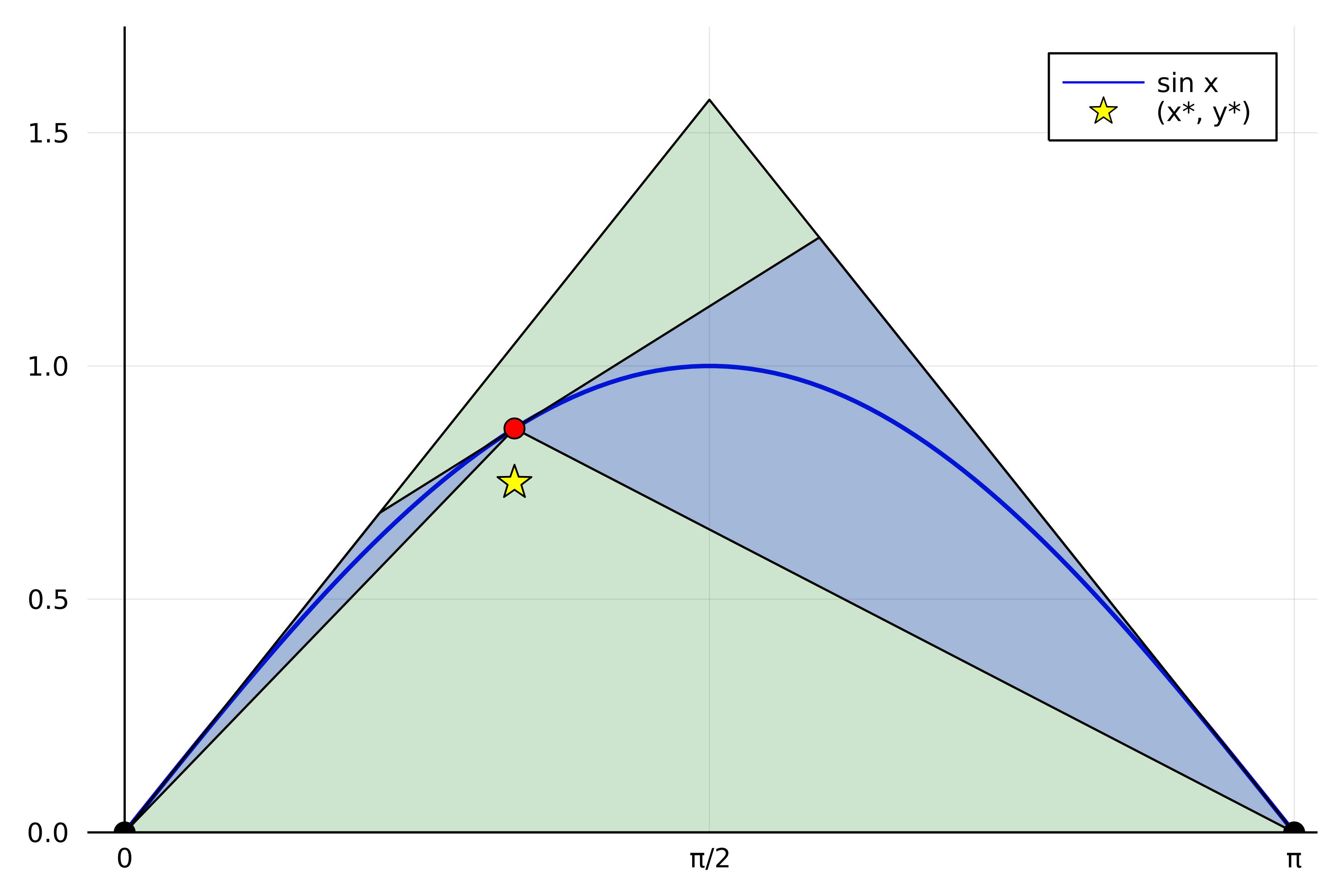}
    \caption{MILP relaxation of $f(x) = \sin x$ over $[0, \pi]$ using $p = (0, \pi)$ (green) and the valid refinement $q = (0, x^*, \pi)$ (blue), which was constructed using a direct refinement scheme. The previous optimal solution to the MILP using $p$ (yellow star) is no longer feasible when constructing the relaxation of $f$ using $q$. }
    \label{fig:direct}
\end{figure}

Unlike bisection, the direct refinement scheme is 
guaranteed to remove the previous MILP optimal solution from the relaxation's feasible space except in the extremely rare (and favorable) case the optimal solution lies on the original curve. Once again consider $f(x) = \sin x$ with domain $[0, \pi]$. As before, we start with the base partition $p = (0, \pi)$, relax $f$ using relaxation \eqref{MILP relaxation : trig} (introducing variable $y$ for the relaxation value of $f(x)$), and solve the resulting MILP giving optimal solution $(x^*, y^*)$. Suppose $(x^*, y^*)$ lies in the interior of the triangle constructed using $p$; this can happen since the optimal solution to the MILP, when projected onto a smaller dimensional space of variables, can lie in the interior in the projected space. Since $x^*$ is in the interior of $[0, \pi]$, we may refine $p$ using the direct refinement scheme to get the valid refinement $q = (0, x^*, \pi)$. Suppose we now use relaxation \eqref{MILP relaxation : trig} using $q$. Since $y^*$ lies on the vertical line $x = x^*$, we see the tangent lines and secant line defined in \eqref{MILP relaxation : trig} will only contain $(x^*, y^*)$ if either (i) $y^* = f(x^*)$ or (ii) a tangent line is vertical, which cannot happen as we require the function to be differentiable. This can be seen in Figure \ref{fig:direct}. A similar discussion holds for bilinear terms that have been relaxed using \eqref{MILP relaxation : bilinear}. 

\subsubsection{Non-Uniform Refinement Scheme}

Another refinement scheme to consider is the non-uniform refinement scheme. In bisection, the location of the added partition point only considered which sub-interval contained $x^*$ and was indifferent to the exact location of $x^*$. In contrast, the direct refinement scheme only considered the exact location of $x^*$. The non-uniform refinement scheme attempts to strike a balance between these two refinement schemes and is as follows. Let $\Delta_1, \Delta_2 > 1$ be user-defined constants and let $x^* \in (x_i, x_{i + 1})$ as before. Define
\begin{equation}
    x'_1 = x^* - \frac{x^* - x_i}{\Delta_1}
\end{equation}
\begin{equation}
    x'_2 = x^* + \frac{x_{i + 1} - x^*}{\Delta_2}
\end{equation}
In the \textit{two-point non-uniform refinement scheme}, points $x'_1$ and $x'_2$ are added to $p$ to yield the valid refinement $q = (x_0, \hdots, x_i, x'_1, x'_2, x_{i + 1}, \hdots, x_m)$. In a sense, the two-point refinement scheme attempts to remove large amounts of the relaxation space, like in bisection, while still being influenced by the exact location of $x^*$. Larger values of $\Delta_1$ and $\Delta_2$ will produce a smaller polyhedron (triangle for trigonometric or tetrahedron for bilinear) containing the curve or surface near $x^*$. We remark that this refinement scheme has been proposed previously in the literature with $\Delta_1 = \Delta_2$ in the context of adaptive partitioning schemes \cite{Nagarajan2019}. In the \textit{three-point non-uniform refinement scheme}, we also add the point $x^*$ to yield the valid refinement $q = (x_0, \hdots, x_i, x'_1, x^*, x'_2, \hdots, x_m)$. Adding point $x^*$ will ensure the previous MILP's optimal solution will be removed from the feasible space, just like in the direct refinement scheme. In the two-point refinement scheme, two binary variables are added for each partition refined using this scheme. In the three-point refinement scheme, three binary variables are added. 

\begin{figure}
    \centering
    \subfloat[]{
    \includegraphics[width = 0.75\textwidth]{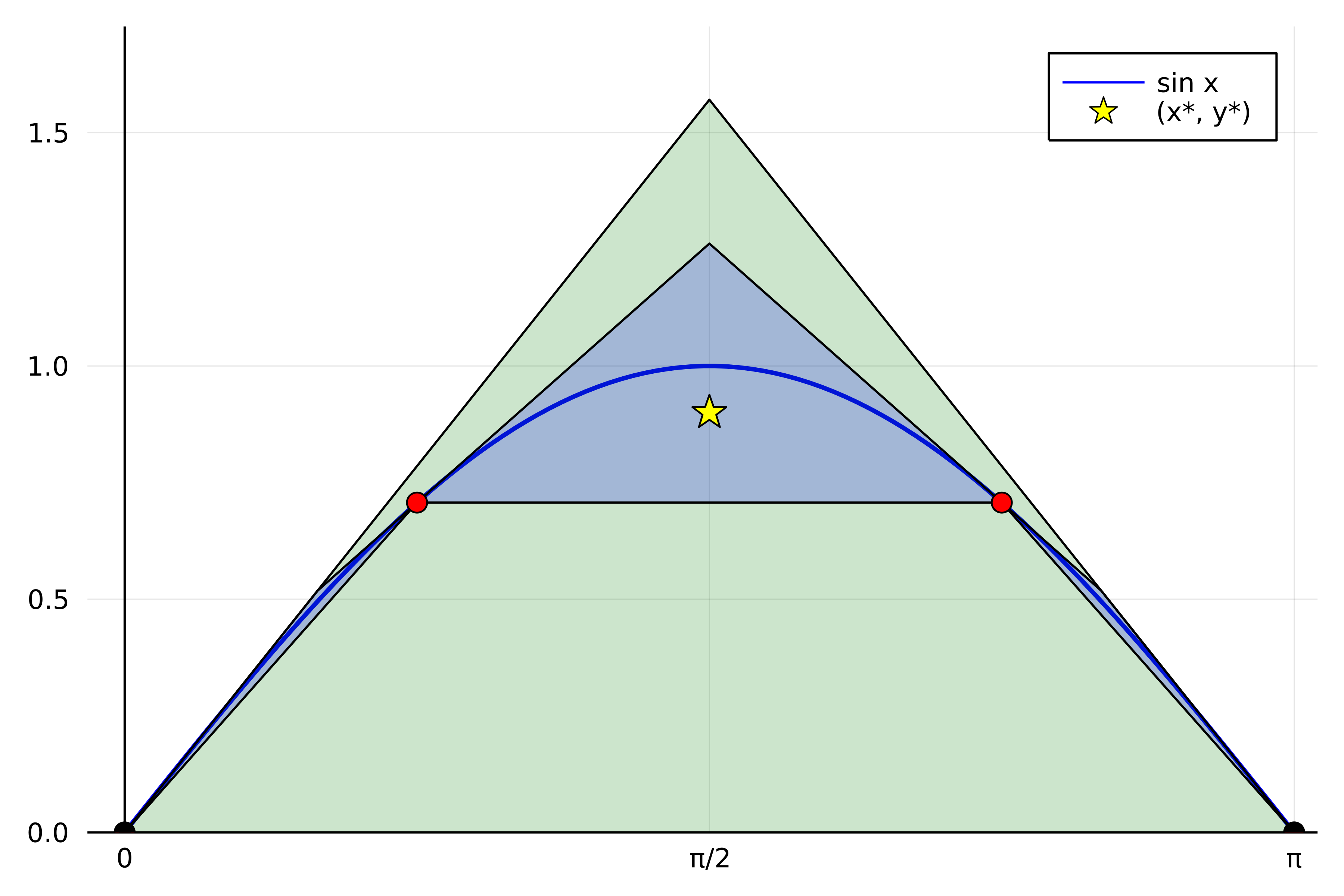}
    \label{fig: two point non uniform}
    }
    \hfill
    \subfloat[]{
    \includegraphics[width = 0.75\textwidth]{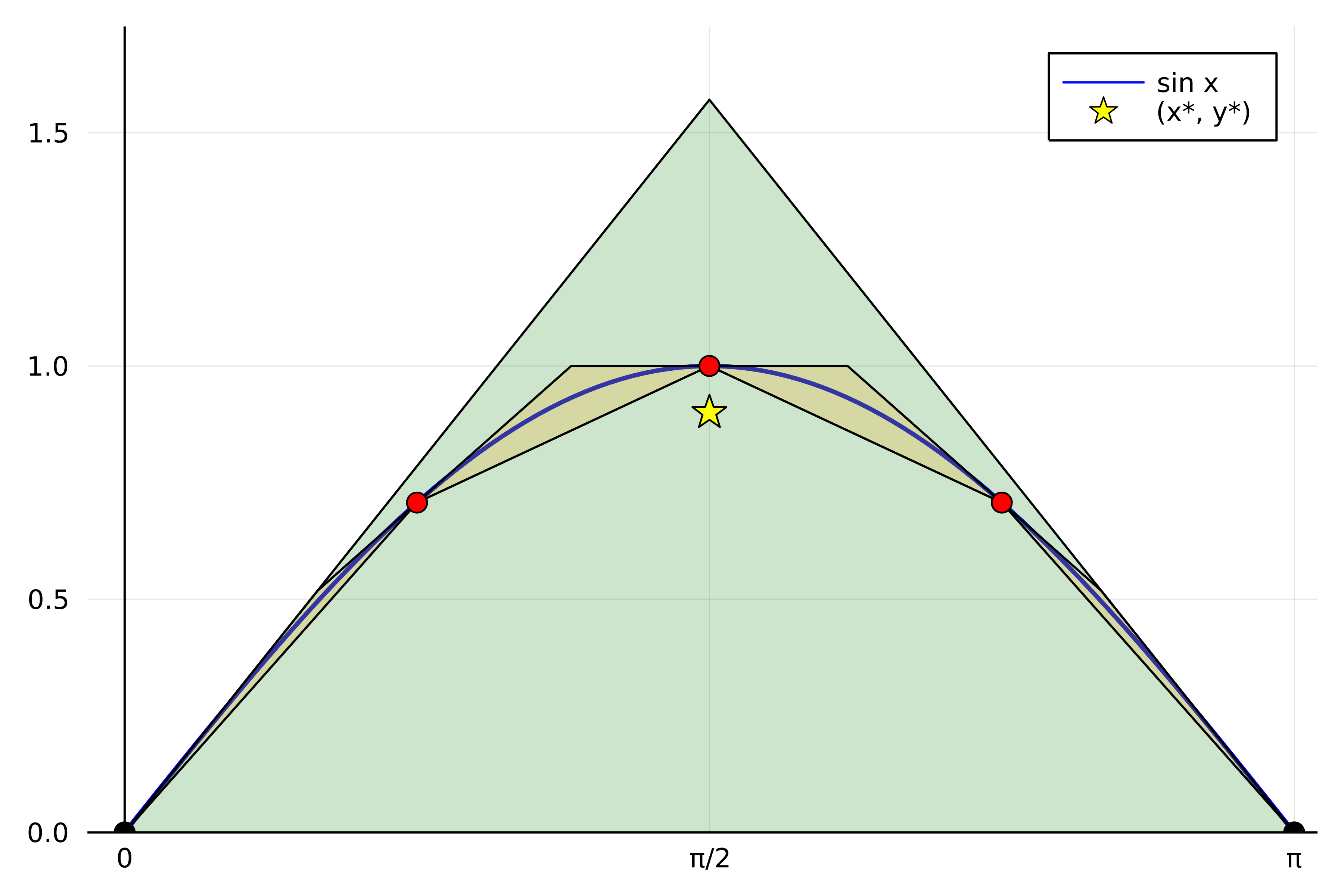}
    \label{fig: three point non uniform}
    }
    \caption{Example showing partition $p = (0, \pi)$ being refined using (a) the two-point non-uniform refinement scheme and (b) the three-point non-uniform refinement scheme. The initial MILP relaxation using $p$ is shown in green. In (a) the optimal solution (yellow star) for the previous iteration remains feasible, while in (b) it is no longer feasible after refinement. In both cases we have set $\Delta_1 = \Delta_2 = 2$.}
    \label{fig: non uniform figures}
\end{figure}

To better visualize both the two-point and three-point non-uniform refinement schemes, once again consider $f(x) = \sin x$ with domain $[0, \pi]$. As before, we start with the initial base partition $p = (0, \pi)$, apply relaxation \eqref{MILP relaxation : trig} for $f$ using $p$ (introducing variable $y$ for the relaxation value of $f(x)$), and solve the resulting MILP giving optimal solution $(x^*, y^*)$. Suppose $(x^*, y^*)$ lies in the interior of the triangle constructed using $p$. As in the bisection refinement scheme, the two-point refinement scheme may not remove $(x^*, y^*)$ from the feasible space after refinement. However, as in the direct refinement scheme, the three-point refinement scheme will always remove $(x^*, y^*)$ from the feasible space except for in the rare case $y^* = f(x^*)$. These observations directly follow from the discussion of the two previous refinement schemes. This can be seen in Figure \ref{fig: non uniform figures}. 

Note that in Figure \ref{fig: non uniform figures}, the three-point refinement scheme results in a significantly tighter relaxation, but has four triangles instead of three in the two-point refinement. This appears to make the three-point refinement scheme the better choice of the two. However, suppose the domain was larger than $[0, \pi]$. The base partition would necessarily consist of more points, leading to more initial triangles. It is possible both refinement schemes may tighten the relaxation over $[0, \pi]$ enough for the optimal solution to move to a different sub-interval and never return to $[0, \pi]$. If this happens, the three-point refinement scheme will have generated an additional binary variable that no longer contributes to the overall convergence of the algorithm. With this in mind, the choice between using two points or three points is not clear in general and would require computational experiments to make an informed decision.

\subsubsection{Consistent Refinement Scheme} \label{Ensuring Consistency}

Consider the function $f$ (trigonometric or bilinear) that has been relaxed using one of the previously described MILP relaxations using partition $p$ for variable $x$. As $p$ is refined using a chosen refinement scheme, the MILP relaxation of $f$ is tightened. When all sub-intervals defined by $p$ are further refined, the MILP relaxation will approach $f$ in the limit. However, in the refinement schemes we have presented we only further refine a single sub-interval to limit the number of additional variables generated at each iteration of the algorithm. We require an additional mechanism so the MILP relaxation approaches $f$ in the limit. Following the standard language used in the literature, we introduce the following definition.
\begin{definition} \label{Consistent Refinement Scheme}
    We say a refinement scheme is \textit{consistent} if at every step any sub-interval defined by the current partition is capable of further refinement and the length of every sub-interval approaches zero (or a small $\varepsilon > 0$) in the limit.
\end{definition}
We note the use of the word ``capable'' in Definition \ref{Consistent Refinement Scheme}. We do not \textit{require} every sub-interval be further refined. This concept is similar to the concept of a consistent branch-and-bound algorithm, where the feasible space may be subdivided (branching) and all sub-domains may be explored unless it is known a sub-domain cannot contain the optimal solution (pruning). 

The previous three refinement schemes are not consistent. In order to make these refinement schemes consistent, we add the following simple rule. 

\begin{rule_def} \label{Consistency rule}
    Let $p$ be a given partition and let $[x_i, x_{i + 1}]$ be a sub-interval of $p$ to be further refined by adding one or more partition points according to a refinement scheme. If $|x_i - x_{i + 1}| < \varepsilon$ for small $\varepsilon > 0$, bisect the largest sub-interval defined by $p$ instead. 
\end{rule_def}
This rule ensures the algorithm does not get stuck refining the same initial sub-interval when the optimal solution lies in a different sub-interval. Using the largest sub-interval ensures all sub-intervals initially defined by a base partition are \textit{capable} of further refinement. The choice of bisection is for simplicity and alternative refinement schemes could be used. Note that if the optimal solution to the MILP relaxation results in $x$ being equal to some $x_i \in \mathcal{P}(p)$, we do not refine $p$ for that iteration and so we would not need to invoke this rule.

Using any of the previously described refinement schemes along with Consistency Rule \ref{Consistency rule}, we are guaranteed the MILP relaxation of a function (trigonometric or bilinear) will approach the function in the limit \cite{Nagarajan2019}. This only concerns a single variable in $\mathcal{F}$ and so we need additional criteria to ensure \textit{all} relaxations of functions approach their original functions in the limit. This is the subject of Section \ref{Refinement Strategy}.

\subsection{Refinement Strategies (Partition Selection)} \label{Refinement Strategy}

We first introduce the following definition.

\begin{definition}
    Let $\mathscr{P}$ be the set of partitions used to construct the MILP relaxation of $\mathcal{F}$.  The procedure used to choose a subset $\mathscr{P}_1 \subseteq \mathscr{P}$ of partitions to be refined is referred to as a \textit{refinement strategy}.
\end{definition}
Let $\mathcal{R}(\mathcal{F}, \mathscr{P})$ denote the relaxation of $\mathcal{F}$ using the previously described MILP relaxations using a set of partitions $\mathscr{P}$ (after sharing partitions). Consider a partition $p \in \mathscr{P}$ corresponding to function $f$ (a similar discussion holds for when $p$ is shared among a set of functions) with partitioned variable $x$. In Section \ref{Refinement Scheme Section}, it was noted the MILP relaxation of $y = f(x)$ will approach the original curve if a consistent refinement scheme is used for $p$. This only considers a single partition and so the purpose of this section is extend this to the set $\mathscr{P}$. In doing so, we will ensure the optimal solution of $\mathcal{R}(\mathcal{F}, \mathscr{P})$ approaches the optimal solution of $\mathcal{F}$ in the limit as further refinements are performed. Similar to the concept of a consistent refinement scheme, we will introduce the concept of a consistent refinement strategy to accomplish this.

The remainder of this section is organized as follows. We introduce the concept of a consistent refinement strategy to guarantee the optimal solution of $\mathcal{R}(\mathcal{F}, \mathscr{P})$ approaches the optimal solution of $\mathcal{F}$ in the limit.  We then provide two consistent refinement strategies. The first refinement strategy is the most natural refinement strategy where all partitions are refined at each iteration (i.e., $\mathscr{P}_1 = \mathscr{P}$). The second refinement strategy considers a subset of partitions based on the quality of the relaxations associated with those partitions. 

\subsubsection{Consistent Refinement Strategy}

We now introduce the concept of a consistent refinement strategy.

\begin{definition} \label{Consistent Refinement Strategy}
    We say a refinement strategy is \textit{consistent} if at every iteration of the algorithm every partition used for the MILP is capable of being selected for refinement and in the limit every partition will be refined such that the corresponding relaxations converge to the original functions.
\end{definition}
We note the use of the word "capable" in Definition \ref{Consistent Refinement Strategy}. Similar to consistent refinement scheme, we do not \textit{require} every partition to be refined. In fact, it is possible for one or many initial partitions to never be refined and still have the optimal value of the relaxation of $\mathcal{F}$ approach the optimal solution to $\mathcal{F}$. This can happen when the optimal solution has a point that corresponds to one or many initial partition points. This is similar to the case in branch-and-bound when a feasible solution is found at a node in the branch-and-bound tree and so no branching is needed for that node. The condition a refinement strategy be consistent is equivalent to requiring every relaxation (trigonometric or bilinear) be capable of being tightened at every iteration of the overall algorithm until the gap is sufficiently small. 

\subsubsection{Complete Refinement Strategy}

The most obvious consistent refinement strategy is to refine all partitions $\mathscr{P}$ at each iteration of the overall algorithm. In doing so, all relaxations (trigonometric and bilinear) are tightened at each iteration until the gap is sufficiently small. We will refer to this refinement strategy as the complete refinement strategy.

The complete refinement strategy has advantages and disadvantages. It is easy to implement and avoids additional computations and sorting that may be used in an alternative refinement strategy that selects a subset of partitions. However, a major disadvantage of this strategy is the number of variables (with emphasis on binary variables) will grow at a faster rate than any other refinement strategy that only uses a subset of partitions. As a result, the complete refinement strategy's effectiveness will depend on the number of iterations needed in the overall algorithm (i.e., the number of times the refinement strategy is invoked) before the gap is sufficiently small. If the number of iterations to reach optimality is small, then the growth in the number of variables will be manageable. 

\subsubsection{k-Worst Refinement Strategy}

We next present an alternative consistent refinement strategy that selects a subset of partitions $\mathscr{P}_1 \subseteq \mathscr{P}$ for further refinement at each iteration of the overall algorithm. To do so, we introduce the following definition. 

\begin{definition} \label{def: Measure of Partition}
    Let $p_x \in \mathscr{P}$ be a partition for $x \in [x^L, x^U]$ corresponding to functions $\{f_1(x), \hdots, f_M(x)\}$ with relaxation values $y_1, \hdots, y_M$ in $\mathcal{R}(\mathcal{F}, \mathscr{P})$. Let $(x^*, y_i^*)$ denote the value of $x$ and $y_i$ in the optimal solution to $\mathcal{R}(\mathcal{F}, \mathscr{P})$. We say the \textit{measure of $p_x$}, denoted by $\mu(p_x)$, is then
    
    \begin{equation} \label{eq: measure of partition}
        \mu(p_x) = \max_{1 \leq i \leq M} |f_i(x^*) - y_i^*|.
    \end{equation}
\end{definition}
We note that if $p_x \in \mathscr{P}$ only corresponds to a single function $f$ ($f$ includes both the \eqref{univariate term} and \eqref{bilinear term}), we simply have $\mu(p_x) = |f(x^*) - y^*|$ in Definition \ref{def: Measure of Partition}. The measure of a partition $p_x$ provides a means of quantifying the quality of the MILP relaxations of the corresponding functions. Suppose $\mathcal{R}(\mathcal{F}, \mathscr{P})$ has been solved to optimality and the measure of all partitions in $\mathscr{P}$ have been computed. The proposed alternative consistent refinement strategy is then to choose the $k$ partitions with the $k$ highest measures (with ties broken arbitrarily) for further refinement. We refer to this refinement strategy as the $k$-worst refinement strategy. The integer $k$ is set in advance by the user and requires knowledge of the number of functions to be relaxed in $\mathcal{F}$. The user may also choose to modify $k$ as the algorithm is running, either by increasing or decreasing $k$ to a valid non-zero integer. For the purposes of this paper, we will keep $k$ constant when using this refinement strategy. We note $k = |\mathscr{P}|$ results in the complete refinement strategy.

The $k$-worst refinement strategy has the advantage of reducing the number of variables (with emphasis on binary variables) added each time partitions are refined to tighten the relaxation. However, a major disadvantage of this refinement strategy is the number of iterations needed to get a sufficiently small gap may increase as a result. This may lead to a longer overall computational time and will largely depend on the maximum gap allowed before terminating the algorithm. 

%% file: Sections/domains.tex
\label{Principal Domains Section}
In this section we discuss a simple reformulation for bounded, periodic functions with period $T$ where the corresponding variable has a domain that does not align with a known, more convenient domain of width $T$. We will refer to this more convenient domain as a \textit{principal domain}. A few simple examples of such principal domains for $T = 2\pi$ are $[0, 2\pi]$, $[-\pi, \pi]$, and more generally $[\theta, \theta + 2\pi]$ for any $\theta \in \mathbb R$. For sine and cosine, these example principal domains capture all possible values and so if a variable's domain contains one of these principal domains (i.e., the domain's width is at least $T$), then the points outside a principal domain of interest are in some sense redundant. A similar observation holds for when the variable's original domain can be shifted to contain a principal domain of interest. Noting this, we aim to reduce the number of binary variables introduced in the MILP relaxation of trigonometric functions (or more generally, any bounded and differentiable $T$-periodic function), as these binary variables will likely have a significant impact on the overall computational time needed to solve the original MINLP. We also show the presented reformulation gives the ability to effectively tighten multiple regions of the original domain simultaneously. 

\subsection{Reformulation Using Principal Domains} \label{Principal Domain Reformulation Subsection}

Let $f(x)$ be a periodic, bounded, univariate function with period $T$ and with $x \in [x^L, x^U] \subset \mathbb R$ where $x^U - x^L \geq T$. Let $[\hat x^L, \hat x^U] \subset \mathbb R$ be a principal domain of interest for $f$ where $\hat x^U - \hat x^L = T$. For brevity, we will denote $[x^L, x^U]$ and $[\hat x^L, \hat x^U]$ by $I$ and $\hat I$, respectively. We will focus our attention to the case $\hat I \subseteq I$. The procedure for $\hat I \not \subset I$ can be easily deduced based on the following discussion and will only differ from $\hat I \subseteq I$ by considering where the endpoints of $I$ are relative to the endpoints of $\hat I$. If $I = \hat I$, clearly no additional work needs to be done. In the case $\hat I \subset I$, define
\begin{equation}
    \alpha^L = \left \lfloor \frac{x^L - \hat x^L}{T}\right \rfloor
\end{equation}
\begin{equation}
    \alpha^U = \left \lceil \frac{x^U - \hat x^U}{T} \right \rceil
\end{equation}
where $\lfloor \cdot \rfloor$ and $\lceil \cdot \rceil$ denote the floor and ceil function, respectively. We note $\alpha^L$ and $\alpha^U$ are integers by construction and may take on negative values. We may then reformulate $y = f(x)$ with $x \in [x^L, x^U]$ by
\begin{subequations} \label{Principal Domain Reformulation}
    \begin{equation} \label{Principal Domain Function New}
        y = f(\hat x)
    \end{equation}
    \begin{equation} \label{Principal Domain Map}
        \hat x = x - \alpha T
    \end{equation}
    \begin{equation}
        \hat x \in [\hat x^L, \hat x^U]
    \end{equation}
    \begin{equation}
        x \in [x^L, x^U]
    \end{equation}
    \begin{equation}
        \alpha \in \{\alpha^L, \hdots, \alpha^U\}
    \end{equation}
\end{subequations}
We have introduced a new variable $\hat x$ and made it so the periodic function in \eqref{Principal Domain Function New} takes in $\hat x$ as its argument where $\hat x$ is related to the original variable $x$ by the mapping \eqref{Principal Domain Map}. 
The trade-off for doing this is the introduction of the integer variable $\alpha$ which takes on an integer value between $\alpha^L$ and $\alpha^U$. The variable $\alpha$ represents the number of integer steps of width $T$ a point $x \in [x^L, x^U]$ is from an equivalent point $\hat x \in [\hat x^L, \hat x^U]$. Depending on the original domain $I$, this reformulation may not be more efficient than using the original formulation. However, as we will see in Section \ref{Computational Results Section}, this reformulation may lead to noticeable improvements in the computational time needed to solve the original MINLP. 

\begin{figure}
    \centering
    \includegraphics[width=\textwidth]{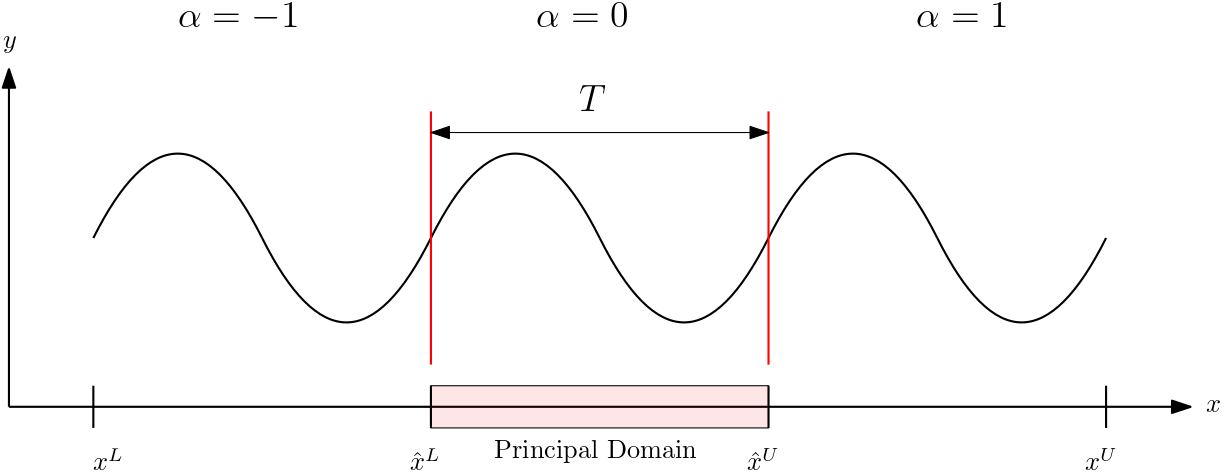}
    \caption{Simple example of selecting a principle domain and relating the points within to the points in the original domain that are outside the chosen principle domain. All points to the right of $\hat x^U$ are $T$ units to the right from a point within $[\hat x^L, \hat x^U]$, indicated by $\alpha = 1$. Similarly, all points to the left of $\hat x^L$ are $T$ units to the left of a point within $[\hat x^L, \hat x^U]$, indicated by $\alpha = -1$.}
    \label{fig:principal domain figure}
\end{figure}

The reformulation \eqref{Principal Domain Reformulation} is illustrated visually in Figure \ref{fig:principal domain figure}. In Figure \ref{fig:principal domain figure}, the principal domain $\hat I$ has been selected so that all points outside of the principal domain are at most $T$ units away from either end of $\hat I$. Because of this, the entire domain $I$ can be captured by introducing $\alpha \in \{-1, 0, 1\}$. 

\subsection{Impact on MILP Relaxation}

In order to see how reformulation \eqref{Principal Domain Reformulation} impacts the MILP relaxations of trigonometric terms as described in Section \ref{Trigonometric Incremental Formulation Subsubsection}, consider the case where the MINLP has constraints $y_1 = \sin(\theta)$ and $y_2 = \cos(\theta)$ with $\theta \in [-4\pi, 4\pi]$. From Section \ref{Partitions Section}, a single partition $p$ may be used for relaxing $y_1$ and $y_2$ simultaneously. We will use the union of the base partitions of $y_1$ and $y_2$ as the shared partition, which will consist of 17 points, noting that we require a point for every multiple of $\frac{\pi}{2}$ in $[-4\pi, 4\pi]$. The corresponding initial MILP relaxation of $y_1$ and $y_2$ using this partition will then consist of 15 binary variables and 32 non-negative continuous variables. We will choose $[0, 2\pi]$ to be the principal domain of interest.  For this principal domain, we may rewrite these constraints as

\begin{subequations} \label{Principal Domain Example}
    \begin{equation}
        y_1 = \sin(\hat \theta)
    \end{equation}
    \begin{equation}
        y_2 = \cos(\hat \theta)
    \end{equation}
\begin{equation}
    \hat \theta = \theta - 2\pi \alpha
\end{equation}
\begin{equation}
    \hat \theta \in [0, 2\pi]
\end{equation}
\begin{equation}
    \theta \in [-4\pi, 4\pi]
\end{equation}
\begin{equation}
    \alpha \in \{-2, -1, 0, 1\}
\end{equation}
\end{subequations}

\noindent From \eqref{Principal Domain Example}, we have changed the domain of the trigonometric functions from $[-4\pi, 4\pi]$ to $[0, 2\pi]$. For this new domain, we can use the smaller shared partition $p = (0, \frac{\pi}{2}, \pi, \frac{3\pi}{2}, 2\pi)$, which only has 5 points (see Figure \ref{fig:sharing partitions example}) and so the corresponding MILP relaxation will consist of 3 binary variables and 8 non-negative continuous variables. In exchange, we have introduced an integer variable $\alpha$ which may take on 4 possible values. Despite the number of binary variables being reduced from 17 to 3, the introduction of $\alpha$ may not lead to improved computational times. This is because $\alpha$ may require multiple branching steps when solving the MILP, which may prove to be too computationally expensive. Even so, this example clearly illustrates the potential benefit of using this reformulation. 

In addition to reducing the number of variables in the resulting MILP relaxation, the presented reformulation effectively allows the relaxation over multiple sub-domains to be tightened simultaneously. For example, once again consider \eqref{Principal Domain Example}. In the original formulation, the functions are relaxed over the domain $[-4\pi, -2\pi] \cup [-2\pi, 0] \cup [0, 2\pi] \cup [2\pi, 4\pi]$. Conversely, in the principal domain reformulation the functions are relaxed over the single principle domain $[0, 2\pi]$. Now, consider the shared (effective base) partition $p = (0, \frac{\pi}{2}, \pi, \frac{3\pi}{2}, 2\pi)$ of $[0, 2\pi]$. If we are now to refine $p$ using any of the previously described refinement schemes, the refinement will be equivalent to adding a partition point to each of the previous sub-domains simultaneously at the expense of a single partition point and the use of the integer variable $\alpha$ in \eqref{Principal Domain Example}. More generally, for $T$-periodic functions, adding a partition point in a principle domain is equivalent to adding a partition point in every sub-domain of width $T$ that has been shifted by $\alpha T$ from the principal domain. In the original formulation, adding multiple partition points to any of the sub-domains may prove to be expensive if the optimal solution lies in a different sub-domain. However, in the principal domain reformulation any added partition point will be beneficial as it will impact \textit{all} of the sub-domains simultaneously. 

Even with the benefits of reformulating the problem to use principal domains, computational studies need to be conducted to determine if this formulation will lead to reduced computational times for the user-specific application at hand. The need for the integer variables $\alpha$ may ultimately lead to longer solve times despite these benefits due to the branching needed to identify the optimal value of $\alpha$.

\subsection{Choice of Principal Domain}

As previously mentioned, there are many choices for the principal domain to be used for the reformulation in general. This is because we may always shift an interval of width $T$ by any amount to produce another principal domain (assuming the function remains properly defined). The choice of principal domain has a direct impact on the performance of the reformulation. To see this, consider $y = \sin(\theta)$ with $\theta \in [-\pi, 3\pi]$. For the principal domain $[0, 2\pi]$, we find $\alpha^L = -1$ and $\alpha^U = 1$, and so $\alpha \in \{-1, 0, 1\}$. However, for the principal domain $[-\pi, \pi]$, we find $\alpha^L = 0$ and $\alpha^U = 1$ and so $\alpha \in \{0, 1\}$. Clearly the principal domain $[-\pi, \pi]$ is superior in this case since $\alpha$ only has two possible values, which will likely reduce the computational effort in solving the resulting MILP relaxation. Therefore, care should be taken when choosing a principal domain for reformulation. 

\subsection{Relating Principal Domain Variables}

Consider two $T$-periodic functions $f_1(\theta_1)$ and $f_2(\theta_2)$ with the same principal domain $[\hat \theta^L, \hat \theta^U]$ and with original domains $I_1 = [\theta_1^L, \theta_1^U]$ and $I_2 = [\theta_2^L, \theta_2^U]$. We will denote by $\alpha_1$ and $\alpha_2$ the integer variables introduced for $\theta_1$ and $\theta_2$ when implementing the principal domain reformulations for $f_1$ and $f_2$. Suppose $\theta_1$ and $\theta_2$ are related by a linking constraint $\theta_2 = h(\theta_1)$. We will assume $h$ is bounded with lower and upper bounds $K_1$ and $K_2$, respectively. We will assume $I_1 \subseteq I_2$, i.e., $I_2$ is the result of expanding $I_1$ a finite amount. Such a case occurs when $\theta_1$ and $\theta_2$ both represent the accumulation of some quantity (time, angle, cost, etc.) in one stage of a process for $\theta_1$ and a subsequent stage for $\theta_2$. We will also assume $K_1 \in I_2$ and $K_2 \in I_2$. In the case $K_1 \notin I_2$ (resp. $K_2 \notin I_2$), we will be unable to improve the lower bound (resp. upper bound) of $\alpha_2$. Define 
\begin{equation} \label{Above Deviation}
    \gamma^U = K_2 - \theta_1^U
\end{equation}
\begin{equation} \label{Below Deviation}
    \gamma^L = \theta_1^L - K_1
\end{equation}
where $\gamma^U$ (resp. $\gamma^L$) is the maximum amount $\theta_2$ can deviate from the maximum (resp. minimum) of $I_1$. Recall in the principal domain reformulation the introduced $\alpha$ integer variable represents how many steps of width $T$ are needed to reach the original variable's value from the chosen principal domain. With this in mind, $\gamma^U$ and $\gamma^L$  give a tighter bound on the number of steps of width $T$ needed to capture all possible values of $\theta_2$ with respect to any value of $\theta_1$ as a result of the linking constraint $\theta_2 = h(\theta_1)$. More specifically, we have
\begin{equation} \label{Alpha Link Inequality}
    -\left\lceil \frac{\gamma^L}{T} \right\rceil \leq \alpha_2 - \alpha_1 \leq \left \lceil \frac{\gamma^U}{T} \right \rceil
\end{equation}

As a simple example, consider $y_1 = \sin(\theta_1)$ and $y_2 = \sin(\theta_2)$ with $\theta_1 \in [-2\pi, 2\pi]$ and $\theta_2 \in [-4\pi, 4\pi]$. We will suppose $\theta_1$ and $\theta_2$ are related by the linking constraint $\theta_2 = \theta_1 + b$ with $0 \leq b \leq 2\pi$. This linking constraint may represent a physical constraint on a vehicle i.e., it may not turn more than an angle of $2\pi$ in addition to the initial angle $\theta_1$ in the stage corresponding to $\theta_2$. Let $[0, 2\pi]$ be the chosen principal domain for both $\theta_1$ and $\theta_2$. Without any modifications, the principal domain reformulations would give
\begin{subequations} \label{Example Relating Alpha}
    \begin{equation}
        y_1 = \sin(\hat\theta_1)
    \end{equation}
    \begin{equation}
        y_2 = \sin(\hat\theta_2)
    \end{equation}
    \begin{equation}
        \hat\theta_1, \hat\theta_2 \in [0, 2\pi]
    \end{equation}
    \begin{equation}
        \theta_1 \in [-2\pi, 2\pi]
    \end{equation}
    \begin{equation}
        \theta_2 \in [-4\pi, 4\pi]
    \end{equation}
    \begin{equation}
        \alpha_1 \in \{-1, 0\}
    \end{equation}
    \begin{equation} \label{alpha 2 example}
        \alpha_2 \in \{-2, -1, 0, 1\}
    \end{equation}
\end{subequations}
We now note by the linking constraint $K_1 = -2\pi$ and $K_2 = 4\pi$. From this we find $\gamma^L = 0$ and $\gamma^U = 1$. Therefore, we may add the additional constraint
\begin{equation} \label{Example Inequality Alpha}
    0 \leq \alpha_2 - \alpha_1 \leq 1
\end{equation}
In words, \eqref{Example Inequality Alpha} states the value of $\theta_2$ is either the same number of $2\pi$ steps from $[0, 2\pi]$ as $\theta_1$ or is one step further in the increasing direction. For this particular example, this reduces \eqref{alpha 2 example} to $\alpha_2 \in \{-1, 0, 1\}$. In general, a direct reduction like this may not always happen. However, the addition of \eqref{Example Inequality Alpha} introduces information to a solver when the value of $\alpha_1$ has been determined. For example, suppose \eqref{Example Relating Alpha} is relaxed using the previously described polyhedral relaxations and is being solved by a MILP solver using branch-and-bound. In one of the branches, we may have $\alpha_1 = -1$. If we include the constraint \eqref{Example Inequality Alpha}, we then know $- 1 \leq \alpha_2 \leq 0$ holds in this branch, halving the number of possible values of $\alpha_2$ from the \eqref{Example Relating Alpha}. As a result, if the MILP solver must branch once more, only two values of $\alpha_2$ need to be considered (in this branch). This may potentially lead to reduced computation times for solving the MILP. Attempting to take advantage of the branching decisions when solving the MILP corresponding to the principal domain reformulation is the subject of the discussion that immediately follows. 

\subsubsection{Branching Decisions}

Consider univariate, bounded trigonometric functions $f_1(x_1), \hdots, f_m(x_m)$ each with period $T$ that have been relaxed using the polyhedral relaxations of Section \ref{MILP Relaxations Section}. The following discussion will hold for general $T$-periodic functions that are relaxed and solved using a branch-and-bound procedure. Furthermore, suppose $x_1, \hdots, x_m$ are linked sequentially, i.e., $x_2 = h_1(x_1), x_3 = h_2(x_2), \hdots, x_m = h_{m - 1}(x_{m - 1})$ where each $h_i$ is bounded and continuous. As previously described, these linking constraints result in inequalities of the form \eqref{Alpha Link Inequality} for consecutive variables $\alpha_i$ and $\alpha_{i + 1}$ when using the same principal domain. When solving the corresponding MILP with these added inequalities, we may elect to branch on each $\alpha_i$ (corresponding to $x_i$) in sequential order. That is, when branching to solve the MILP we first consider $\alpha_1$, then $\alpha_2$, and so on until $\alpha_m$. Consider a branch where we have set $\alpha_i$, $i < m$, to be equal to some value. By the linking constraint $x_{i + 1} = h_i(x_i)$, the resulting inequality of form \eqref{Alpha Link Inequality} for $\alpha_i$ and $\alpha_{i + 1}$ may reduce the number of values to consider for $\alpha_{i + 1}$. Consequently, if we require to branch once more we may do so on $\alpha_{i + 1}$ and require fewer children subproblems for this branch as a result. This will ultimately reduce the number of subproblems needed to be solved in the extreme case of an exhaustive search (by effectively avoiding generating many infeasible subproblems due to the various values of $\alpha_i$). In the case where additional linking constraints between non-consecutive variables are present, a similar approach can be taken to identify a potentially stronger branching scheme. We do not handle these other cases in this paper. 

%% file: Sections/dubins.tex
\label{Motivating Example Section}

The Markov-Dubins path planning problem \cite{kaya2019markov} (MDPPP) is a natural extension of the classical two-point Markov-Dubins problem \cite{dubins1957curves}. The MDPPP may be formulated as follows. Consider a vehicle travelling in the plane starting at a point $p_1 = (x_1, y_1) \in \mathbb{R}^2$, passing through $n - 2$ intermediate points $p_i = (x_i, y_i) \in \mathbb{R}^2$, $i = 2, \hdots, n - 1$, and arriving at a final point $p_n = (x_n, y_n) \in \mathbb{R}^2$, in the sequence $(p_1, \hdots, p_n)$. The vehicle may or may not have a specified initial and final heading angle $\theta_0 \in [0, 2\pi]$ and $\theta_n \in [0, 2\pi]$ at $p_1$ and $p_n$, respectively. The vehicle is a Dubins vehicle \cite{dubins1957curves} and so it has a minimum turning radius $\rho > 0$ and travels at constant speed. The MDPPP is then to find a path of minimum length passing through the points in the specified sequence with radius of curvature at least $\rho$ everywhere from $p_1$ to $p_n$.

\begin{figure}

    \begin{minipage}{0.5\linewidth}
    \centering
    \subfloat[RSR]{
    \includegraphics[width=0.9\textwidth]{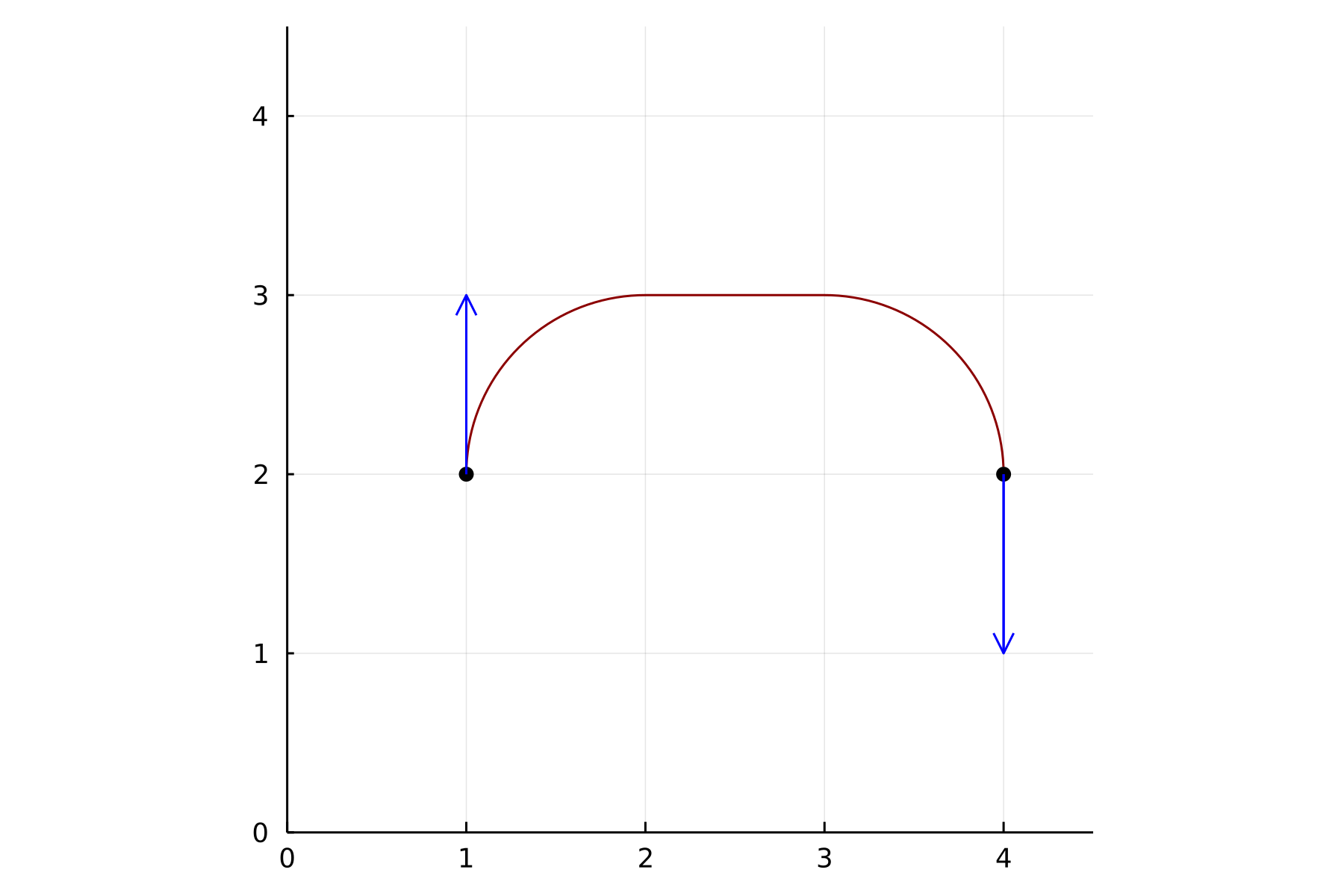}
    }
    \end{minipage}%
    \begin{minipage}{0.5\linewidth}%
    \subfloat[LRL]{
    \includegraphics[width=0.9\textwidth]{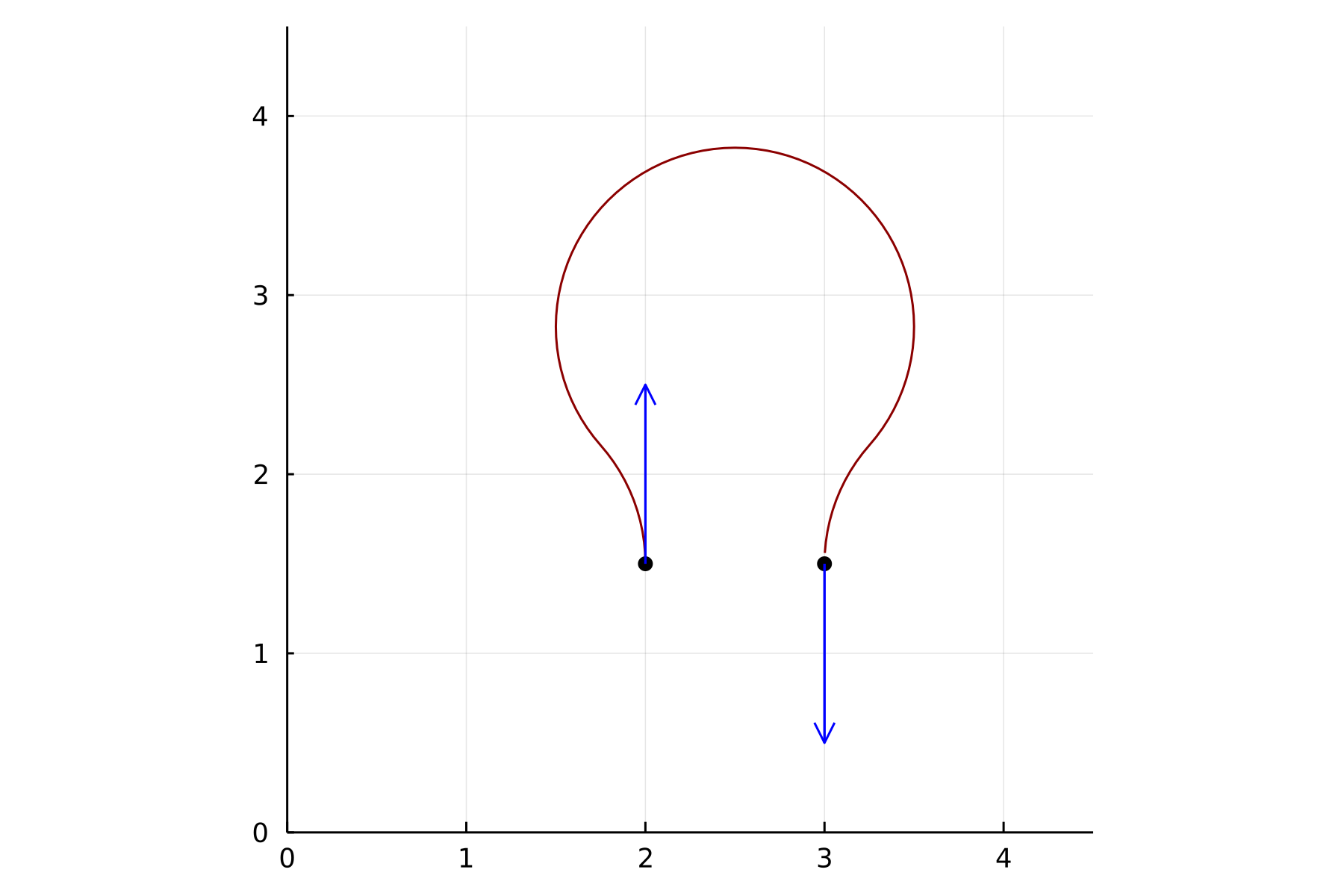}
    }
    \end{minipage}\par\medskip
    \centering
    \subfloat[LSR]{
    \includegraphics[width = 0.5\textwidth]{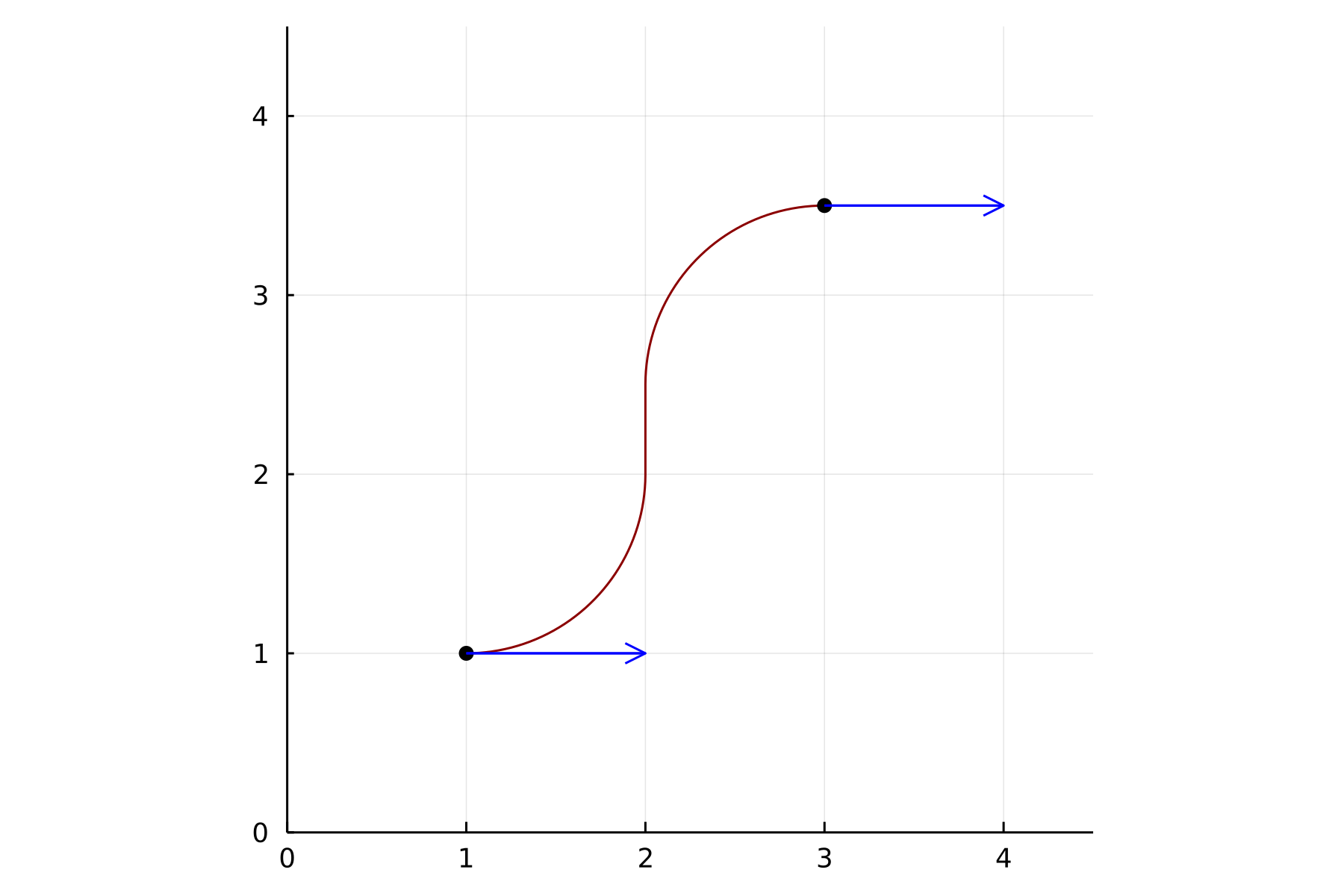}
    }
    \caption{Example of optimal Dubins paths when $n = 2$. Examples of the remaining optimal words are obtained by reversing the initial and final heading directions (blue arrows).}
    \label{Dubins path examples}
\end{figure}

\begin{figure}
    \centering
    \includegraphics[width = 0.8\textwidth]{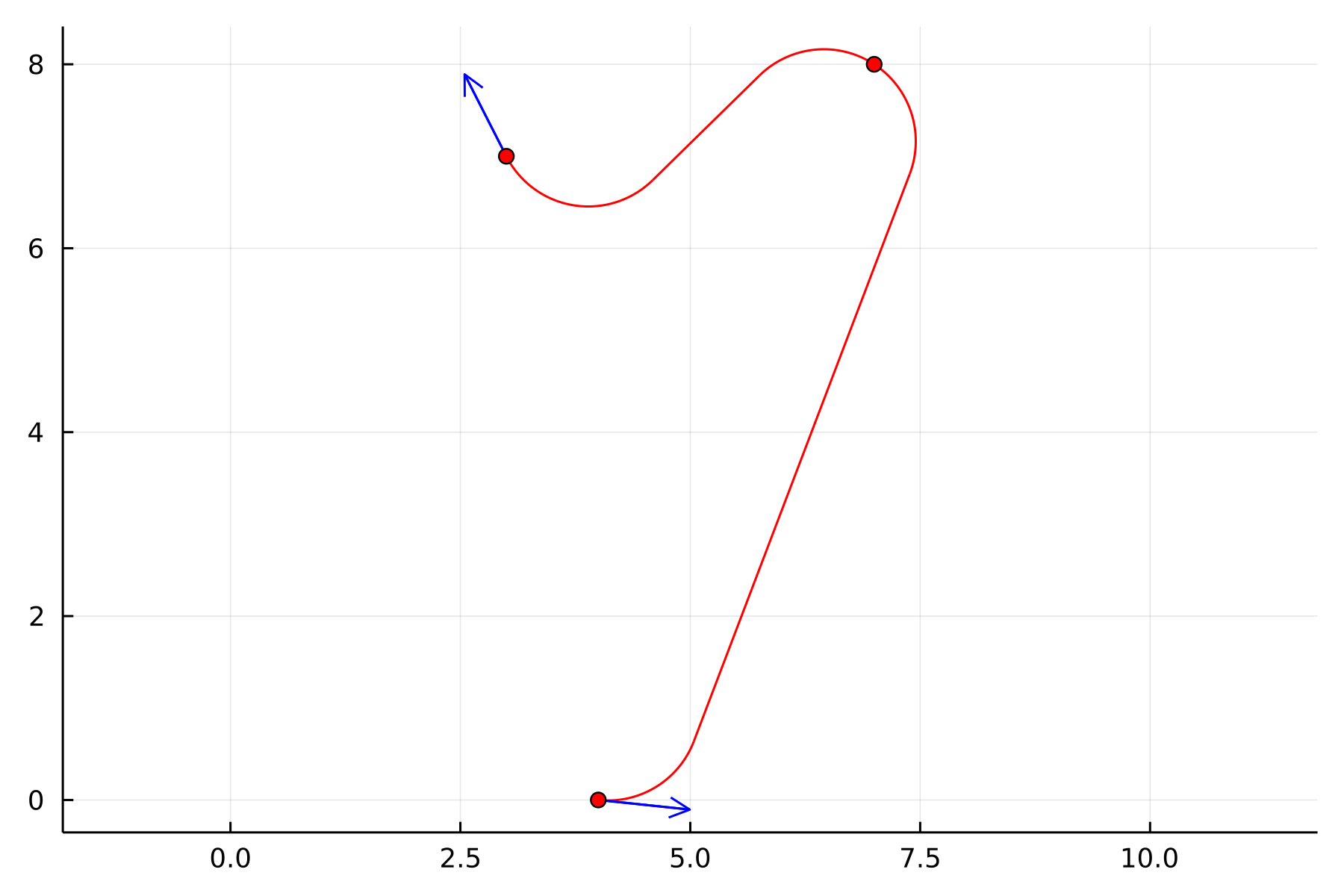}
    \caption{Optimal solution for MDPPP with $n = 3$. The path is LSL followed by LSR.}
    \label{fig: mdppp example}
\end{figure}

In Dubins' seminal work \cite{dubins1957curves}, it was shown the optimal path in the case where $n = 2$ is of type $CCC$ or $CSC$ where $C$ denotes a circular arc segment (left or right turn) of radius $\rho$ and $S$ denotes a straight line segment. For the circular arc segments of radius $\rho$, let a left turn (counter-clockwise) be denoted by $L$ and a right turn (clockwise) be denoted by $R$. For $n = 2$, the optimal solution is then one of the following: $LRL$, $RLR$, $LSL$, $LSR$, $RSR$, or $RSL$. Examples of $RSR$, $LRL$, and $LSR$ are shown in Figure \ref{Dubins path examples}. Examples of $LSL$, $RLR$, and $RSL$ are obtained by reversing the paths shown in Figure \ref{Dubins path examples}. We note that a segment in the optimal solution may have zero length. By Bellman's principle of optimality \cite{bellman1952dp}, the subpath connecting two consecutive points in the optimal solution to the MDPPP must then also be of type $CCC$ or $CSC$. An example solution to the MDPPP for $n = 3$ is shown in Figure \ref{fig: mdppp example}.

Following the work of \cite{kaya2019markov}, we can formulate the MDPPP as a NLP as follows. We refer to the subpath connecting consecutive points $p_i$ and $p_{i +1}$ as stage $i$. There are $N = n - 1$ stages. We note that the word describing the optimal subpath at stage $i$ is a subset of the sequence $(L, R, S, L, R)$. For each stage $i$, define $\xi^i_j \geq 0$, where $j = 1, \hdots, 5$ denotes which letter ($L$, $R$, or $S$) in the sequence $(L, R, S, L, R)$ and $\xi^i_j$ represents the length of the corresponding segment type. (For example, if the path taken at stage $i$ is $LRL$, then $\xi^i_3 = \xi^i_5 = 0$ and $\xi^i_1, \xi^i_2, \xi^i_4 \geq 0$.) For $j \in \{1, 2, 4, 5\}$, define $\theta^i_j$ to be heading angle of the vehicle \textit{at the end} of the corresponding segment ($L$ or $R$) at stage $i$. Let $\theta^i_0$ be the heading angle at the start of stage $i$, i.e., the heading angle when the vehicle reaches point $p_i$ or the initial heading angle in the case of $p_1$. Let $\mathcal{S} = \{1, \hdots, N\}$. The NLP formulation \cite{kaya2019markov} (after factoring) of the MDPPP is then:
\begin{subequations} \label{MDPPP NLP}
\allowdisplaybreaks 
    \begin{align}
        &(MDPPP) \qquad \text{minimize } \sum_{i = 1}^N \sum_{j = 1}^5 \xi^i_j \label{MDPPP: obj}\\
        \text{s.t.} \quad 
        &x_i - x_{i + 1} + \rho \left[-w_0^i + 2 w_1^i - 2w_2^i + 2w^i_4 - w^i_5 \right] + \mu^i = 0,  \quad \forall i \in \mathcal{S} \label{MDPPP: horiz}
        \\
        &y_i - y_{i + 1} + \rho \left[z^i_0 - 2 z^i_1 + 2 z^i_2 - 2z^i_4 + z^i_5 \right] + \nu^i = 0, \quad \,\qquad \forall i \in \mathcal{S} \label{MDPPP: vert}\\
        &w^i_j = \sin(\theta^i_j), \quad \forall i \in \mathcal{S}, \quad j \in \{0, 1, 2, 4, 5\} \label{MDPPP: w}\\
        &z^i_j = \cos(\theta^i_j), \,\quad \forall i \in \mathcal{S}, \quad j \in \{0, 1, 2, 4, 5\} \label{MDPPP: z}\\
        &\mu^i = \xi^i_3 z^i_2, \,\,\qquad \forall i \in \mathcal{S} \label{MDPPP: mu}\\
        &\nu^i = \xi^i_3 w^i_2, \,\qquad \forall i \in \mathcal{S} \label{MDPPP: nu}\\
        \begin{split}
        &\theta^i_1 = \theta^i_0 + \rho^{-1} \xi^i_1,
        \qquad \theta^i_2 = \theta^i_1 - \rho^{-1}\xi^i_2,
        \\
        &\theta^i_4 = \theta^i_2 + \rho^{-1}\xi^i_4,
        \qquad \theta^i_5 = \theta^i_4 - \rho^{-1}\xi^i_4, \quad \forall i \in \mathcal{S}
        \end{split} \label{MDPPP: turns}\\
        &w^{i + 1}_0 = w^i_5, \,\,\,\qquad i = 1, \hdots, N - 1 \label{MDPPP: w continuity} \\
        &z^{i + 1}_0 = z^i_5, \quad \qquad i = 1, \hdots, N - 1 \label{MDPPP: z continuity} \\
        \begin{split}
            &\theta^i_0 \in [0, 2\pi], \quad\quad\quad\,\,\,\,\, \theta^i_1 \in [0, 4\pi] \\
            &\theta^i_2, \theta^i_4 \in [-2\pi, 4\pi], \quad \theta^i_5 \in [-4\pi, 4\pi], \quad \forall i \in \mathcal{S} 
        \end{split} \label{MDPPP: box constraints}\\
        &\theta^1_0 = \theta_0, \quad\qquad
        w^N_5 = \sin(\theta_n), \quad\qquad
        z^N_5 = \cos(\theta_n) \label{MDPPP: boundary} \refstepcounter{equation} \tag{\theequation*} \\
        %
        &\xi^i_j \geq 0, \,\,\,\,\,\quad \qquad \forall i \in \mathcal{S}, \quad j = 1, \hdots, 5
    \end{align}
\end{subequations}
Objective \eqref{MDPPP: obj} says to minimize the total path length from $p_1$ to $p_n$. Constraints \eqref{MDPPP: horiz} and \eqref{MDPPP: vert} represent the horizontal and vertical displacement, respectively, after taking a $LRSLR$ path \cite{shkel2001classification,kaya2019markov} from $p_i$ to $p_{i + 1}$ for each stage $i$, where the trigonometric and bilinear terms have been factored out with constraints \eqref{MDPPP: w}-\eqref{MDPPP: z} and \eqref{MDPPP: mu}-\eqref{MDPPP: nu}, respectively. The left and right turns of radius $\rho$ are characterized by constraints \eqref{MDPPP: turns}, where an increase (resp. decrease) in angle from one segment to another corresponds to a left (resp. right) turn. Constraints \eqref{MDPPP: w continuity}-\eqref{MDPPP: z continuity} ensure the slope of the path at the end of stage $i$ and at the beginning of stage $i + 1$ are equal, i.e., continuity is maintained.

We have chosen to use this problem to illustrate the effectiveness of our proposed algorithm for several reasons. Firstly, all nonlinearities present in \eqref{MDPPP NLP} are trigonometric (sine and cosine) or bilinear. The MDPPP provides a natural example where sharing partitions for both trigonometric and bilinear terms is possible. The trigonometric terms are defined over closed intervals given by \eqref{MDPPP: box constraints}. The bilinear terms are also capable of being defined over closed intervals by using knowledge of the optimal solutions to the MDPPP (see Section \ref{Section Bounding Arc Lengths}). Secondly, by the nature of the problem the box constraints in \eqref{MDPPP: box constraints} are, in some sense, relatively loose. This motivates the idea of using principal domains to replace the original domains of the trigonometric terms. Finally, the complexity of the MDPPP will likely increase the influence of the chosen refinement schemes and refinement strategy on the overall time to solve. This will help in comparing these refinement schemes and refinement strategies.

\subsection{Additional Constraints} \label{Additional Constraints Subsection}

We now introduce some additional constraints that can be added to the MILP formulation of the MDPPP in order to reduce computation time. These are included to ensure the computational times stay within a reasonable range.

\subsubsection{Bounding Arc Lengths} \label{Section Bounding Arc Lengths}

In order to apply the MILP relaxations for the trigonometric terms \eqref{MDPPP: w}-\eqref{MDPPP: z} and the bilinear terms \eqref{MDPPP: mu}-\eqref{MDPPP: nu}, we require the arc lengths of the straight line segments, $\xi^i_3$, be bounded above (since we require a closed interval). Consider two consecutive points $p_i$ and $p_{i + 1}$ in the MDPPP. It can be shown that in an optimal solution of the MDPPP we must have
\begin{equation} \label{straight line segment upper bound}
    \xi^i_3 \leq \|p_i - p_{i + 1}\| + 4 \rho, \quad \forall i \in \mathcal{S}
\end{equation}
where $\| \cdot \|$ denotes the Euclidean norm. We may also bound the circular arc segments by noting that taking more than a full circle in any circular segment will never be optimal. Therefore,
\begin{equation} \label{circular arc segment upper bound}
    \xi^i_j \leq 2 \pi \rho, \quad \forall i \in \mathcal{S}, \quad j \in \{1, 2, 4, 5\}
\end{equation}

\subsubsection{Restricting the Number of Segments}

As previously mentioned, the sub-path for two consecutive points in the optimal solution to the MDPPP must be of type $CCC$ or $CSC$. In the current formulation, we are not strictly requiring at most three segments be used at each stage. This is because the optimal solution will automatically satisfy this property. In order to reduce the computation time needed to solve the MDPPP, we can add constraints enforcing at most three segments are used at each stage. For each stage $i$ and segment $j$, let $\beta^i_j$ be a binary variable where $\beta^i_j = 1$ if segment $j$ is used in stage $i$ and $\beta^i_j = 0$ otherwise. We can then add the constraints
\begin{equation} \label{binary arc segment}
    \sum_{j = 1}^5 \beta^i_j \leq 3, \quad \forall i \in \mathcal{S}
\end{equation}
and replace the upper bounds \eqref{straight line segment upper bound} and \eqref{circular arc segment upper bound} by
\begin{equation} \label{arc length bounded using binary}
    \xi^i_j \leq M^i_j \beta^i_j, \quad \forall i \in \mathcal{S}, \quad j = 1, \hdots, 5
\end{equation}
where $M^i_j$ represents an upper bound for the arc length $\xi^i_j$ as described in \eqref{straight line segment upper bound} and \eqref{circular arc segment upper bound}. If tighter bounds are available, they could be easily substituted into \eqref{arc length bounded using binary}. By including constraints \eqref{binary arc segment} and \eqref{arc length bounded using binary}, the formulation for the MDPPP becomes an MINLP.

\subsubsection{CSC Conditions} \label{CSC Conditions Section}

We next take advantage of two results that limit the type of sub-paths in an optimal solution of the MDPPP for two or more consecutive points. 

Consider two points $p_i$ and $p_{i + 1}$ and suppose the Euclidean distance between the two points is at least 4$\rho$. It can be shown the optimal Dubins path from $p_i$ and $p_{i + 1}$ cannot be of type $CCC$. Suppose points $p_i$ and $p_{i + 1}$ in the MDPPP, corresponding to stage $i$, are at least 4$\rho$ from each other. We may then introduce the following constraints.
\begin{align}
    \beta^i_1 + \beta^i_2 &\leq 1 \label{CSC First Turn} \\
    \beta^i_3 &= 1 \label{CSC Straight} \\
    \beta^i_4 + \beta^i_5 &\leq 1 \label{CSC Second Turn}
\end{align}
Recall $\beta^i_1, \beta^i_4$ and $\beta^i_2, \beta^i_5$ correspond to a left turn and a right turn, respectively, in the sequence $(L, R, S, L, R)$. Similarly, $\beta^i_3$ corresponds to a straight line segment. Constraints \eqref{CSC First Turn} and \eqref{CSC Second Turn} then simply says at most two turns are taken between points $p_i$ and $p_{i + 1}$. Additionally, constraint \eqref{CSC Straight} requires a straight line segment be used, though it may be a degenerate case with zero length (i.e., we still allow for $\xi^i_3 = 0$). 

Constraints \eqref{CSC First Turn}-\eqref{CSC Second Turn} impose additional constraints for a single stage due to optimality conditions from the physics of the problem. We can go further by also considering consecutive stages. In \cite{kaya2019markov}, the following result was derived using optimal control theory for the MDPPP. Consider two consecutive stages $i$ and $i + 1$ and suppose the Euclidean distance between consecutive points is at least 4$\rho$ for both stages. Therefore, each stages will only admit a $CSC$ path in the optimal solution to the MDPPP. Under these conditions, we can make use of two results in \cite{kaya2019markov} (see Theorem 4 and Proposition 2 therein) to get the following corollary.
\begin{corollary} \label{Consecutive CSC Corollary}
    Suppose two consecutive stages $i$ and $i + 1$ are such that both only admit $CSC$ paths in the optimal solution to the MDPPP. Then the type ($L$ or $R$) of the final turn in stage $i$ must be the same type ($L$ or $R$) as the first turn in stage $i + 1$. Furthermore, the length of the final turn in stage $i$ and the length of the first turn in stage $i + 1$ must be equal. Additionally, the first turn ($L$ or $R$) in stage $i + 1$ has length less than $\pi \rho$.
\end{corollary}
From Corollary \ref{Consecutive CSC Corollary}, we get the following additional constraints for consecutive stages admitting only $CSC$ paths
%
\begin{equation}
    \beta^i_4 = \beta^{i + 1}_1
\end{equation}
\begin{equation}
    \beta^i_5 = \beta^{i + 1}_2
\end{equation}
\begin{equation}
    \xi^i_4 = \xi^{i + 1}_1
\end{equation}
\begin{equation}
    \xi^i_5 = \xi^{i + 1}_2
\end{equation}
\begin{equation} \label{Consecutive CSC Final Equation}
    \xi^{i + 1}_1 + \xi^{i + 1}_2 \leq \pi \rho 
\end{equation}
Constraints \eqref{CSC First Turn}-\eqref{Consecutive CSC Final Equation} can be added to the formulation for the MDPPP by a pre-processing step before solving. We note that constraint \eqref{Consecutive CSC Final Equation} could have also been written as
\begin{equation}
    \xi^i_4 + \xi^i_5 \leq \pi \rho
\end{equation}
by the first statement in Corollary \ref{Consecutive CSC Corollary}.


%% file: Sections/results_updated.tex
\label{Computational Results Section}

In this section, we present computational results related to the MDPPP described in Section \ref{Motivating Example Section}.

\subsection{Problem Generation}

Each instance consists of $n$ points, with $n$ ranging from 2 to 8, placed within a $10 \times 10$ grid. For each $n$, ten instances were generated by randomly placing $n$ points with the constraint that consecutive points are separated by a Euclidean distance of at least $4\rho$, where $\rho = 1$ for all instances. This requirement reduces overall computation time by enabling the use of additional constraints discussed in Section \ref{Additional Constraints Subsection}. The initial and final heading angles for each instance were randomly selected from the interval $[0, 2\pi]$.

\subsection{Implementation Details}

All computations were performed on a 64-bit Windows machine with an AMD Ryzen 7 2700 processor (3.2 GHz) and 16 GB RAM. All MILP models were implemented in Julia \cite{bezanson2017julia} using the JuMP modeling package \cite{DunningHuchetteLubin2017}, and solved using CPLEX 20.1 with default settings. The implementation of the algorithms presented in this article is open-sourced and made available at \url{https://github.com/yduaskme/NLP-Dubins}. A solution was considered optimal when the relative gap was under 1\%. Each instance had a one-hour time limit. Warm-starting was employed to improve computational efficiency: all variables were warm-started except binary and non-negative continuous variables involved in polyhedral relaxations, due to bookkeeping complexity.

Feasibility-based bounds tightening (FBBT) \cite{belotti2012feasibility} was used to narrow the initial domain and reduce the number of required partition points. Although FBBT can be applied repeatedly, for this MINLP a single application was sufficient to reach a fixed point. While this does not guarantee the tightest bounds, further tightening was not possible using feasibility alone.

Angle variables $\theta$ were initially partitioned by placing points at multiples of $\pm \frac{\pi}{2}$ from zero, which ensures breakpoints of both sine and cosine functions are captured and shared. For arc-length variables $\xi$, the tightened domain was used as the base partition.

As described in Section \ref{Overview Section}, at each iteration a feasible solution is constructed from the MILP relaxation. This is done by extracting the initial heading angle for each stage ($\theta^i_0$ for $i = 1, \ldots, N$) and computing the shortest Dubins path using these angles. This feasible solution helps update the best-known upper bound. Note that this procedure does not guarantee an improved bound at every iteration.

The final results compare the complete refinement strategy with the $k$-worst refinement strategy for $k$ values of 25\% and 50\%. The case of $k = 100\%$ corresponds to the complete refinement strategy.

\subsection{Results Overview}

Results are grouped into three main comparisons:

\begin{enumerate}
    \item \textbf{Original Formulation vs. Principal Domain Reformulation} — to demonstrate the benefits of reformulating the MDPPP.
    \item \textbf{Comparison of Variable Refinement Schemes} — to identify the most effective partitioning approach.
    \item \textbf{Complete vs. $k$-Worst Refinement Strategy} — to evaluate performance under partial refinement.
\end{enumerate}

\subsubsection{Original vs. Principal Domain Reformulation}

We first compare the original MDPPP formulation \eqref{MDPPP NLP} with the reformulation using principal domains, as described in Section \ref{Principal Domain Reformulation Subsection}. For simplicity, each angle is constrained to a principal domain of $[0, 2\pi]$. Tables \ref{tab:full original} and \ref{tab:full reformulation} in the Appendix \ref{sec:appendix} present the detailed results with additional instance level discussion and observations.

\begin{figure}
    \centering
    \subfloat[]{\includegraphics[width=0.45\textwidth]{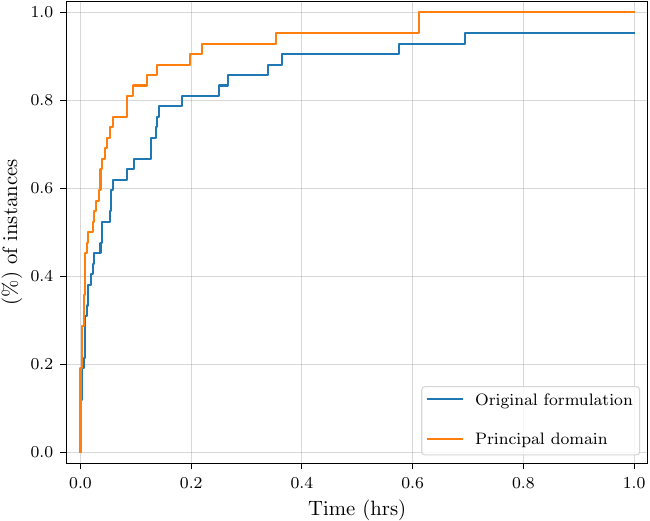}}\hfill 
    \subfloat[]{\includegraphics[width=0.45\textwidth]{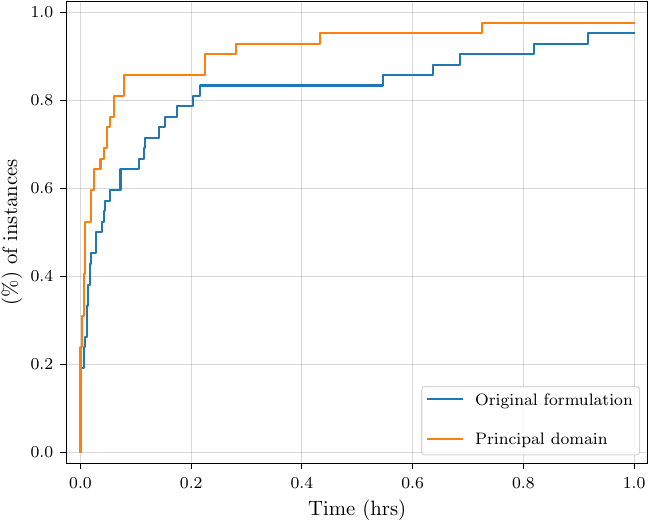}} \\
    \subfloat[]{\includegraphics[width=0.45\textwidth]{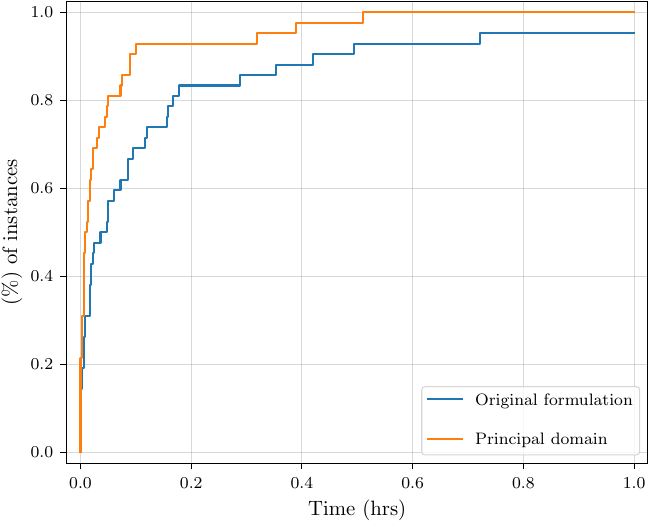}} \hfill
    \subfloat[]{\includegraphics[width=0.45\textwidth]{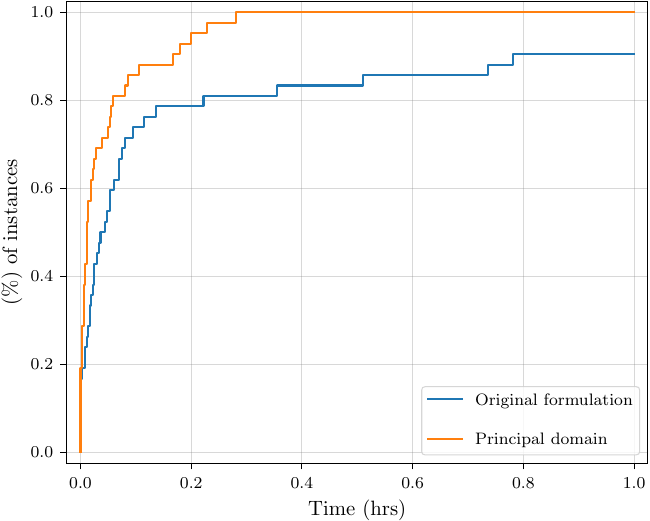}}
    \caption{Performance comparison of original vs. principal domain reformulation for: (a) Bisection, (b) Direct, (c) NU2, and (d) NU3 schemes.}
    \label{fig:performance_plot_orig_vs_principal}
\end{figure}

Figure \ref{fig:performance_plot_orig_vs_principal} shows performance plots for the percentage of instances solved to optimality over time. Each plot corresponds to a refinement scheme: bisection, direct, non-uniform two-point (NU2), and non-uniform three-point (NU3). It is evident that the principal domain reformulation significantly outperforms the original across all schemes.
Complete runtime and iteration data are available in the Appendix \ref{sec:appendix}.

\subsubsection{Comparison of Variable Refinement Schemes} \label{sec:refinement-schemes}

This section compares the performance of the four refinement schemes (bisection, direct, NU2, NU3) using only the principal domain reformulation, as it consistently outperforms the original formulation.
\begin{figure}
    \centering
    \includegraphics[width=0.5\linewidth]{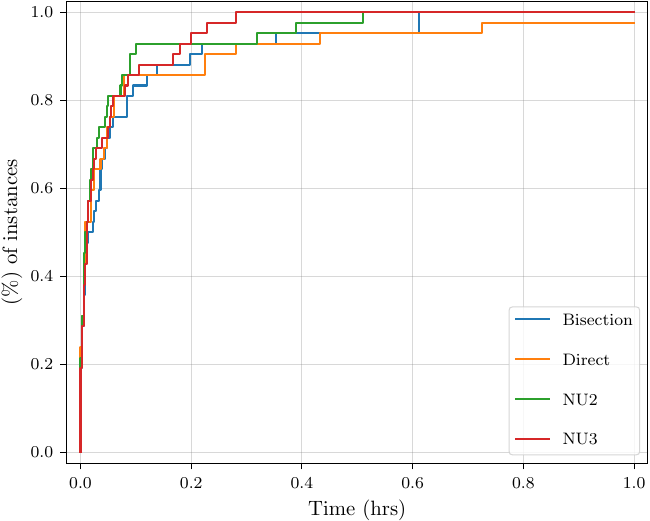}
    \caption{Performance plots for different variable refinement schemes with principal domain reformulation.}
    \label{fig:refinementschemes}
\end{figure}
As shown in Figure \ref{fig:refinementschemes}, NU2 and NU3 significantly outperform bisection and direct schemes. While bisection and direct methods are common in the literature, these results highlight the superiority of non-uniform approaches for MINLPs involving trigonometric functions. Between NU2 and NU3, neither is consistently superior, so we recommend testing both in practice.

\subsubsection{Complete Refinement vs. $k$-Worst Refinement}

Tables \ref{tab:k worst bisection}–\ref{tab:k worst nu3} (Appendix \ref{sec:appendix}) provide complete data for comparing full vs. partial ($k$-worst) refinement. Here, we summarize the findings using performance plots for NU2 and NU3 only, given their superior performance.

\begin{figure}
    \centering
    \subfloat[]{\includegraphics[width=0.45\textwidth]{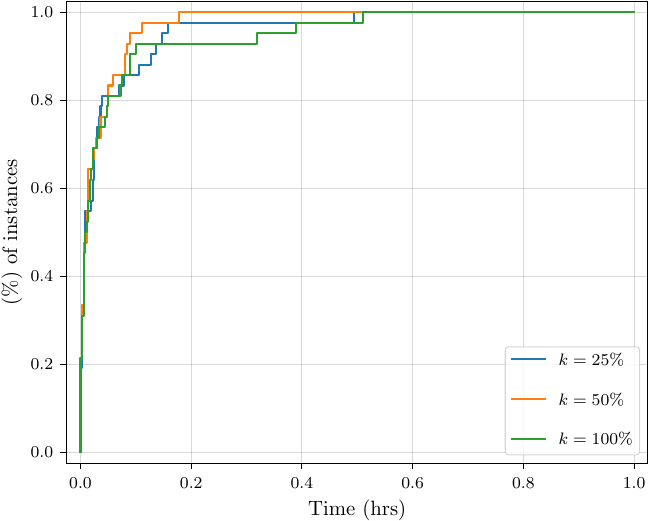}}\hfill 
    \subfloat[]{\includegraphics[width=0.45\textwidth]{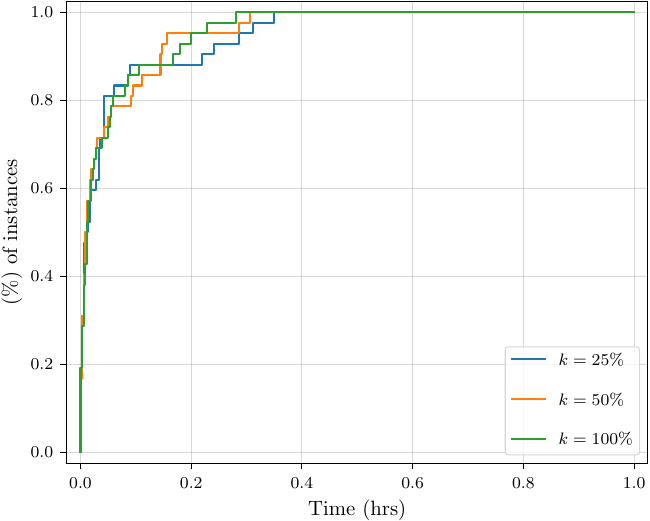}} 
    \caption{Performance comparison of $k$-worst refinement strategies for: (a) NU2 and (b) NU3 with $k = 25\%, 50\%, 100\%$.}
    \label{fig:nu2_nu3_k}
\end{figure}

\begin{figure}
    \centering
    \includegraphics[width=0.5\linewidth]{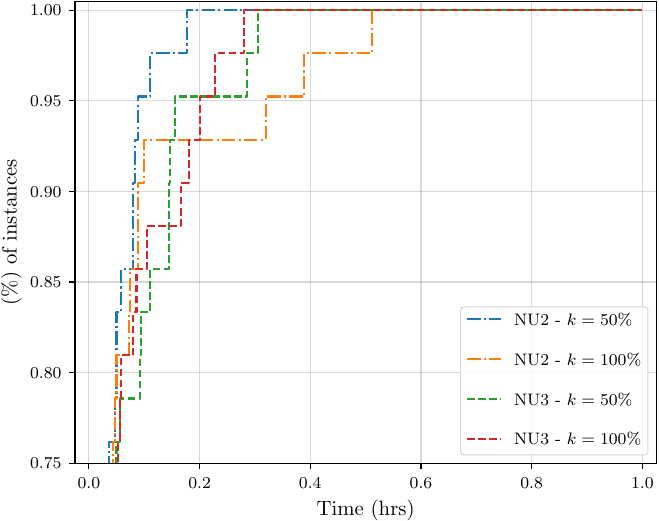}
    \caption{Refined performance comparison highlighting top configurations from Figure \ref{fig:nu2_nu3_k}.}
    \label{fig:nu23}
\end{figure}

Figure \ref{fig:nu2_nu3_k} shows that $k=50\%$ and $k=100\%$ outperform $k=25\%$ in both refinement schemes. Figure \ref{fig:nu23} focuses on these top-performing strategies, showing that NU2 with $k = 50\%$ achieves the best performance for most instances. Hence, we adopt this configuration as our default algorithm.

\subsubsection{Comparison with Other Solvers}

Finally, we compare our method with state-of-the-art MINLP solvers: Gurobi (commercial) and SCIP (open-source). All solvers were limited to 10 threads. We restrict this comparison to instances with 4 and 5 targets. Just in these small instances, SCIP times out on all of them, and Gurobi times out on 3 out of the 10 instances with 5 targets. Hence, we did not see the need to perform any runs for SCIP and Gurobi on the larger instances. The values in Table \ref{tab:solvers} clearly indicate the merit of the methods proposed in this article. 

\begin{table}[h]
    \centering
    \setlength\tabcolsep{3pt}
    \begin{tabular}{|c|c|c|c|}
    \hline
    Instance & Gurobi & SCIP & NU2 ($k=50\%$) \\
    \hline
    4--1 & 4 & ** & 4 \\
    4--2 & $<$1 & ** & 4 \\
    4--3 & 57 & ** & 1 \\
    4--4 & 123 & ** & 3 \\
    4--5 & 61 & ** & 17 \\
    4--6 & 29 & ** & 4 \\
    4--7 & 310 & ** & 4 \\
    4--8 & 40 & ** & 4 \\
    4--9 & 18 & ** & 8 \\
    4--10 & 20 & ** & 3 \\
    - & - & - & - \\
    5--1 & 1529 & ** & 14 \\
    5--2 & 339 & ** & 21 \\
    5--3 & 1840 & ** & 17 \\
    5--4 & 311 & ** & 22 \\
    5--5 & ** & ** & 14 \\
    5--6 & ** & ** & 22 \\
    5--7 & ** & ** & 18 \\
    5--8 & ** & ** & 24 \\
    5--9 & 1968 & ** & 28 \\
    5--10 & 2262 & ** & 21 \\
    \hline
    \end{tabular}
    \caption{Runtime (in seconds) comparison with Gurobi and SCIP. `**' indicates the run exceeded the 1-hour time limit.}
    \label{tab:solvers}
\end{table}

%% file: Sections/conclusion.tex
\label{Conclusion} This article develops an MILP-based iterative algorithm to solve an MINLP with periodic functions. In particular, we consider MINLPs with trigonometric functions that arise in a path planning problem for a single fixed-wing aerial vehicle. Domain partitioning approaches that partition variables that are shared between different nonlinear terms, exploration of different partition refinement schemes and strategies, and reformulations that leverage the principle domain information of the periodic functions are presented. Extensive computational results that corroborate the effectiveness of the proposed approach and study the pros and cons of the different partition refinement schemes and the strategies are presented. Future work would focus on showing the effectiveness of the approach on general MINLPs that not only have nonlinearities that arise from periodic functions but also contain other types of nonlinearities.  Another potential future direction for research is to consider non-smooth periodic functions that cannot be well-approximated by Fourier series.

%% file: Sections/appendix_new.tex
\label{sec:appendix}
\subsection{Original Formulation vs. Principal Domain Reformulation}
\normalsize

All tables in the Appendix follow a consistent format. The first column, `\textbf{instance}', lists entries of the form `$n$--$q$', where $n$ represents the number of points and $q$ the instance number for that value of $n$. Column `\textbf{t}' denotes the time (in seconds) taken by the algorithm to terminate, either due to achieving a relative gap below 1\% or reaching the one-hour time limit. Entries marked with `\textbf{**}' indicate instances that exceeded the time limit. Column `\textbf{iter.}' reports the number of refinement iterations performed, while `\textbf{bin}' shows the number of binary variables added to the formulation due to refinement (equivalently, the number of partition points introduced). The header row in each table indicates the refinement scheme used—bisection, direct, non-uniform two-point (NU2), or non-uniform three-point (NU3).

Both Table \ref{tab:full original} and Table \ref{tab:full reformulation} report results using the complete refinement strategy, meaning all partitions were refined at each iteration of the algorithm.

In Table \ref{tab:full original}, most instances are solved within the one-hour time limit, with the exception of instances 7-2 and 7-7 under the NU3 scheme. Notably, instances with $n = 7$ take significantly longer to solve than those with $n = 6$, regardless of the refinement method. However, both groups require a similar number of refinement iterations to meet the convergence criterion, and the number of added binary variables is also comparable. This suggests that the longer runtime for $n = 7$ is not due to refinement overhead but rather the difficulty in solving the MILP relaxation of the original formulation. Specifically, the failure to account for the periodic nature of trigonometric functions likely introduces ambiguity that increases with $n$, resulting in a substantial rise in computational effort. Among the refinement methods, NU2 shows the most consistent performance as $n$ increases.

In contrast, Table \ref{tab:full reformulation} shows that all instances up to $n = 7$ are solved within the time limit. The increase in computational time from $n = 6$ to $n = 7$ is much less dramatic compared to the original formulation, particularly under the NU2 and NU3 schemes. As before, the number of refinement iterations and added binary variables are similar between $n = 6$ and $n = 7$, suggesting that the reformulation using principal domains is primarily responsible for mitigating the rise in computational effort. By constraining each angle to a fixed principal domain, we simplify the relationships between stages and reduce ambiguity. Note that these results depend on the choice of principal domains; other selections may lead to further improvements. Nonetheless, Table \ref{tab:full reformulation} demonstrates the effectiveness of using principal domains over broader angle ranges. For $n = 8$, however, all refinement schemes encounter instances that exceed the time limit, indicating a sharp increase in difficulty. Even so, NU2 performs best among the four schemes at this size.

\input{Standalone_figures/full_original}
\input{Standalone_figures/full_reformulation}

\subsection{Complete Refinement Strategy vs. $k$-Worst Refinement Strategy}

The results for the bisection, direct, NU2, and NU3 refinement schemes under the complete and $k$-worst refinement strategies are shown in Tables \ref{tab:k worst bisection}--\ref{tab:k worst nu3}. We restrict attention to instances with $n \geq 6$, where refinement strategy choice has a more pronounced effect.

In Table \ref{tab:k worst bisection}, the three strategies show comparable performance under the bisection scheme. However, $k = 25\%$ performs slightly worse on average, often requiring more refinement iterations than the other two strategies.

In Table \ref{tab:k worst direct}, the $k$-worst strategy with $k = 50\%$ outperforms the others, sometimes by a substantial margin. The $k = 25\%$ strategy again requires more iterations, contributing to longer solve times. For most instances, $k = 50\%$ and the complete refinement strategy yield similar iteration counts. In cases where the complete strategy performs better, it typically requires fewer or equal iterations compared to $k = 50\%$.

Table \ref{tab:k worst nu2} shows similar trends under the NU2 scheme. The $k = 50\%$ strategy consistently outperforms the others, particularly in instances 8--2, 8--3, 8--5, 8--8, 8--9, and 8--10. As before, $k = 25\%$ tends to require more iterations. For $n = 8$, the performance of $k = 25\%$ is notably worse, suggesting that the extra time spent on additional iterations outweighs the benefit of fewer binary variables. Interestingly, the number of added binary variables for $k = 25\%$ is not significantly lower than for $k = 50\%$ in the $n = 8$ cases. However, for smaller $n$, $k = 25\%$ is able to achieve a 1\% gap with significantly fewer binary variables—especially evident in the $n = 7$ cases.

In Table \ref{tab:k worst nu3}, $k = 50\%$ again emerges as the best-performing strategy overall. The $k = 25\%$ strategy requires more iterations on average. Notably, the complete refinement strategy performs worse here, with three $n = 8$ instances exceeding the time limit. Despite this, the number of binary variables added is similar across all strategies (except possibly for instance 8--2), reinforcing that runtime bottlenecks are not primarily due to binary variable count.

\input{Standalone_figures/k_worst_bisection}
\input{Standalone_figures/k_worst_direct}
\input{Standalone_figures/k_worst_nu2}
\input{Standalone_figures/k_worst_nu3}

%% file: Standalone_figures/full_original.tex



\setlength\tabcolsep{2.5pt}
\renewcommand{\arraystretch}{0.95}  
\begin{filecontents*}{Data/FullOriginal.csv}
2--1,0,0.11,4,17,0,0.3,4,15,0,0.75,2,16,0,0.75,2,23,
2--2,1,0.44,4,20,0,0.13,4,18,0,0.9,2,20,0,0.09,3,43,
2--3,1,0.68,4,16,0,0.93,2,6,0,0.62,2,16,0,0.35,2,22,
2--4,0,0.13,5,11,0,0.39,3,7,0,0.77,3,14,0,0.3,3,23,
2--5,1,0.31,6,29,1,0.23,4,20,0,0.95,3,30,0,0.67,3,45,
2--6,1,0.43,4,16,0,0.36,4,16,0,0.32,3,24,0,0.53,3,36,
2--7,0,0.65,1,4,0,0.28,1,4,0,0.92,1,8,0,0.28,1,12,
2--8,1,0.93,5,20,0,0.75,3,11,0,0.19,3,24,0,0.19,3,34,
2--9,1,0.89,4,16,0,0.82,3,12,0,0.56,3,23,0,0.3,3,34,
2--10,0,0.64,2,8,0,0.07,2,8,0,0.69,1,8,0,0.65,1,12,
--,--,--,--,--,--,--,--,--,--,--,--,--,--,--,--,--,
3--1,2,0.78,4,41,2,0.44,4,37,2,0.64,3,60,2,0.43,3,88,
3--2,2,0.87,5,49,2,0.37,4,34,2,0.39,4,79,2,0.32,3,86,
3--3,3,0.29,5,50,2,0.53,4,34,2,0.73,3,60,3,0.23,3,88,
3--4,3,0.81,7,72,2,0.31,5,47,2,0.89,4,80,2,0.63,4,120,
3--5,4,0.66,6,61,3,0.45,5,45,2,0.76,3,58,3,0.5,3,94,
3--6,2,0.65,5,50,2,0.59,4,34,2,0.39,3,60,2,0.28,3,82,
3--7,1,0.56,4,40,1,0.28,3,27,1,0.21,3,60,1,0.27,2,60,
3--8,2,0.57,5,55,1,0.65,3,27,2,0.47,3,66,2,0.27,3,93,
3--9,2,0.57,4,40,2,0.18,4,37,1,0.49,2,40,1,0.43,2,58,
3--10,2,0.63,4,44,2,0.65,4,40,2,0.21,3,66,2,0.28,3,97,
--,--,--,--,--,--,--,--,--,--,--,--,--,--,--,--,--,
4--1,7,0.77,4,61,8,0.66,4,61,8,0.25,3,92,5,0.83,2,96,
4--2,9,0.61,5,80,5,0.82,4,57,6,0.94,3,95,7,0.51,3,141,
4--3,7,0.8,4,64,4,0.85,3,44,6,0.29,3,96,4,0.81,2,96,
4--4,6,0.5,5,80,6,0.49,4,62,6,0.72,3,96,7,0.69,3,144,
4--5,31,0.76,7,116,23,0.96,6,95,21,0.91,4,127,33,0.32,4,191,
4--6,17,0.58,5,82,8,0.35,3,48,13,0.35,3,96,7,0.44,2,93,
4--7,15,0.79,7,109,9,0.31,5,65,11,0.61,4,122,9,0.66,3,133,
4--8,12,0.66,5,85,7,0.6,3,50,10,0.32,3,102,14,0.22,3,152,
4--9,27,0.47,7,113,27,0.56,4,62,33,0.62,4,132,31,0.48,4,189,
4--10,9,0.3,5,71,8,0.29,4,54,7,0.32,3,84,5,0.97,2,78,
--,--,--,--,--,--,--,--,--,--,--,--,--,--,--,--,--,
5--1,39,0.72,4,89,46,0.33,3,67,25,0.82,2,92,54,0.88,2,137,
5--2,85,0.81,5,115,77,0.75,4,86,140,0.36,4,183,94,0.45,3,206,
5--3,35,0.81,4,91,42,0.26,4,85,62,0.2,3,138,60,0.1,3,203,
5--4,140,0.37,6,134,104,0.89,5,102,85,0.79,4,173,169,0.42,4,254,
5--5,57,0.84,5,110,61,0.77,4,85,66,0.36,3,132,78,0.21,3,195,
5--6,98,0.71,6,110,67,0.65,4,75,74,0.71,3,116,114,0.36,3,172,
5--7,39,0.87,5,108,39,0.37,4,86,30,0.7,3,132,66,0.45,3,198,
5--8,51,0.9,4,89,44,0.82,4,86,37,0.55,3,135,87,0.32,3,200,
5--9,45,0.86,4,91,58,0.63,3,64,69,0.23,3,134,49,0.11,3,198,
5--10,77,0.78,5,107,51,0.45,4,78,77,0.49,3,120,99,0.45,3,188,
--,--,--,--,--,--,--,--,--,--,--,--,--,--,--,--,--,
6--1,135,0.54,5,142,101,0.28,4,104,98,0.26,3,168,136,0.54,3,244,
6--2,214,0.76,5,141,196,0.69,4,103,181,0.67,3,165,253,0.4,3,248,
6--3,208,0.77,4,114,143,0.96,3,81,179,0.28,3,173,257,0.04,3,249,
6--4,468,0.73,6,173,424,0.94,5,132,425,0.93,3,174,417,0.61,3,256,
6--5,492,0.45,7,197,415,0.33,5,136,268,0.74,3,168,191,0.67,3,243,
6--6,304,0.75,6,160,384,0.47,6,141,347,0.4,4,212,179,0.72,3,227,
6--7,147,0.84,5,131,161,0.78,4,94,187,0.88,3,152,191,0.39,3,223,
6--8,461,0.82,6,174,736,0.29,5,140,313,0.41,4,232,220,0.6,3,260,
6--9,354,0.94,6,166,262,0.58,5,133,435,0.38,4,219,272,0.28,4,326,
6--10,193,0.58,5,135,263,0.33,4,95,227,0.48,3,162,128,0.91,2,150,
--,--,--,--,--,--,--,--,--,--,--,--,--,--,--,--,--,
7--1,205,0.88,4,137,551,0.74,4,133,318,0.43,3,204,345,0.82,2,201,
7--2,669,0.72,5,166,1975,0.64,5,157,573,0.44,3,200,**,5.34,1,102,
7--3,516,0.72,4,140,781,0.3,4,135,1782,0.85,3,210,804,0.34,3,309,
7--4,964,0.88,4,134,159,0.88,4,130,1047,0.76,3,199,1849,0.33,3,292,
7--5,2072,0.87,4,136,2477,0.5,4,136,1270,0.64,3,204,2816,0.34,3,305,
7--6,1229,0.77,6,205,2954,0.36,5,164,1519,0.53,3,204,2656,0.43,3,302,
7--7,2501,0.61,5,173,3307,0.64,5,162,2605,0.5,4,276,**,3.99,2,205,
7--8,506,0.87,5,172,638,0.52,4,129,566,0.91,3,210,1283,0.65,3,311,
7--9,1318,0.56,6,200,515,0.64,4,130,645,0.43,3,188,491,0.21,3,291,
7--10,905,0.85,5,171,2297,0.44,5,149,607,0.52,3,202,290,0.95,2,197,
--,--,--,--,--,--,--,--,--,--,--,--,--,--,--,--,--,
8--1,**,1.54,4,158,**,2.9,3,123,**,1.63,2,164,**,6.69,1,123,
8--2,**,1.02,5,200,**,1.78,3,118,3527,0.95,3,238,**,8.05,1,114,
8--3,**,1.09,4,155,**,1.52,3,118,2696,0.43,3,240,**,1.93,2,237,
8--4,820,0.65,5,198,**,1.52,3,111,1156,0.58,3,230,959,0.95,2,225,
8--5,**,2.59,4,161,**,2.95,3,117,**,3.7,2,158,**,3.25,2,235,
8--6,**,4.05,3,123,**,8.48,2,82,**,4.85,2,164,**,10.83,1,123,
8--7,1905,0.99,4,159,**,3.85,2,79,**,2.21,2,160,**,4.18,1,117,
8--8,**,1.39,4,152,**,2.66,3,112,**,2.49,2,150,**,1.95,2,216,
8--9,1720,0.51,5,195,**,4.29,2,73,**,1.08,2,154,**,6.2,1,111,
8--10,**,1.89,4,154,**,1.5,4,145,**,1.32,2,152,**,1.33,2,228,
,,,,,,,,,,,,,,,,,
\end{filecontents*}
\csvreader[
column count = 17,
    no head,
    before reading=\small\sisetup{round-mode=places, round-precision=2},
    after reading=\normalsize,
    longtable=|c|rrr|rrr|rrr|rrr|,
    filter expr = {test{\ifnumgreater{\thecsvinputline}{22}}},
    table head = \caption{Original formulation results.\label{tab:full original}}\\ \hline
     & \multicolumn{3}{c|}{Bisection} &\multicolumn{3}{c|}{Direct} &\multicolumn{3}{c|}{NU2} & \multicolumn{3}{c|}{NU3} \\ 
    \hline 
    \bfseries instance & \bfseries t & 
    \bfseries iter & 
    \bfseries bin &
    \bfseries t & 
    \bfseries iter. & 
    \bfseries bin &
    \bfseries t & 
    \bfseries iter. & 
    \bfseries bin &
    \bfseries t & 
    \bfseries iter. & 
    \bfseries bin
     \\ 
    \hline
    \endfirsthead
    \caption{Original formulation results (continued).}\\
    \hline 
    \bfseries instance & \bfseries t & 
    \bfseries iter & 
    \bfseries bin &
    \bfseries t & 
    \bfseries iter. & 
    \bfseries bin &
    \bfseries t & 
    \bfseries iter. & 
    \bfseries bin &
    \bfseries t & 
    \bfseries iter. & 
    \bfseries bin
     \\ 
    \hline\endhead
    \hline\endfoot,
    late after line=\\, 
    ]{Data/FullOriginal.csv}{
        1=\instance, 
        2=\tbisect, 
        3=\gapbisect, 
        4=\iterbisect,
        5=\binarybisect,
        6=\tdirect, 
        7=\gapdirect,
        8=\iterdirect,
        9=\binarydirect,
        10=\tnutwo,
        11=\gapnutwo,
        12=\iternutwo,
        13=\binarynutwo,
        14=\tnuthree,
        15=\gapnuthree,
        16=\iternuthree,
        17=\binarynuthree}%
    {\instance & \tbisect & \iterbisect & \binarybisect & \tdirect & \iterdirect & \binarydirect & \tnutwo & \iternutwo & \binarynutwo & \tnuthree & \iternuthree & \binarynuthree}


%% file: Standalone_figures/full_reformulation.tex



\setlength\tabcolsep{2.5pt}
\begin{filecontents*}{Data/FullReformulation.csv}
2--1,0,0.51,3,15,0,0.3,3,11,0,0.26,2,20,0,0.22,2,26
2--2,0,0.45,4,20,0,0.14,4,18,0,0.95,2,20,0,0.09,3,43
2--3,0,0.68,4,20,0,0.63,3,9,0,0.4,2,20,0,0.39,2,25
2--4,0,0.13,5,15,0,0.39,3,8,0,0.79,3,18,0,0.31,3,24
2--5,1,0.27,6,30,0,0.25,4,19,0,0.87,3,30,0,0.63,3,43
2--6,0,0.43,4,20,0,0.77,3,14,0,0.44,3,30,0,0.27,3,42
2--7,0,0.28,2,10,0,0.62,1,5,0,0.01,2,20,0,0.62,1,15
2--8,0,0.93,5,25,0,0.78,3,12,0,0.59,2,20,0,0.49,2,28
2--9,0,0.89,4,20,0,0.28,4,16,0,0.19,3,29,0,0.1,3,37
2--10,0,0.65,2,10,0,0.07,2,9,0,0.97,1,10,0,0.65,1,15
--,--,--,--,--,--,--,--,--,--,--,--,--,--,--,--,--
3--1,1,0.77,4,44,0,0.76,2,22,2,0.22,4,88,1,0.59,2,66
3--2,1,0.88,5,55,1,0.32,4,37,2,0.39,4,85,1,0.08,3,89
3--3,2,0.28,5,55,1,0.52,3,29,1,0.44,3,66,1,0.55,2,63
3--4,2,0.81,5,55,1,0.7,4,43,1,0.78,3,66,1,0.32,3,97
3--5,2,0.61,6,63,2,0.38,5,47,1,0.96,3,60,2,0.78,3,88
3--6,2,0.93,5,55,1,0.84,4,32,1,0.74,3,66,1,0.54,3,85
3--7,1,0.61,4,44,0,0.84,2,21,1,0.32,3,66,0,0.48,2,63
3--8,1,0.56,5,52,1,0.37,4,34,1,0.72,3,66,1,0.44,3,93
3--9,1,0.57,4,44,1,0.17,4,37,1,0.57,2,44,1,0.52,2,61
3--10,2,0.61,5,55,1,0.8,4,40,1,0.22,3,66,1,0.06,3,97
--,--,--,--,--,--,--,--,--,--,--,--,--,--,--,--,--
4--1,4,0.73,4,68,2,0.47,3,46,5,0.26,3,102,2,0.88,2,99
4--2,5,0.59,5,85,3,0.76,4,59,3,0.85,3,102,4,0.51,3,143
4--3,2,0.91,3,51,1,0.89,2,34,2,0.47,2,68,2,0.36,2,99
4--4,3,0.27,5,85,3,0.35,4,54,2,0.72,3,100,3,0.57,3,143
4--5,12,0.9,6,99,12,0.64,6,92,11,0.54,4,130,22,0.25,4,191
4--6,9,0.58,5,82,3,0.47,3,46,9,0.41,3,96,3,0.58,2,93
4--7,7,0.72,7,111,3,0.79,4,60,4,0.28,4,129,6,0.13,4,174
4--8,5,0.66,5,85,3,0.54,3,50,4,0.52,3,102,14,0.19,3,152
4--9,12,0.74,5,82,10,0.73,4,62,10,0.48,3,99,6,0.94,2,99
4--10,8,0.46,4,66,2,0.61,3,42,3,0.43,3,98,4,0.22,3,137
--,--,--,--,--,--,--,--,--,--,--,--,--,--,--,--,--
5--1,21,0.61,4,89,9,0.86,2,45,11,0.5,2,92,16,0.37,2,137
5--2,40,0.79,5,115,20,0.95,3,69,21,0.53,3,138,43,0.45,3,206
5--3,16,0.85,4,92,31,0.46,4,84,18,0.22,3,138,28,0.06,3,203
5--4,44,0.52,5,110,33,0.79,4,82,23,0.89,3,134,43,0.3,3,184
5--5,33,0.84,5,115,24,0.72,4,87,33,0.56,3,138,27,0.3,3,199
5--6,40,0.69,6,128,33,0.55,5,89,23,0.73,3,121,54,0.18,3,183
5--7,16,0.93,4,89,11,0.51,3,66,18,0.21,3,138,15,0.77,2,135
5--8,22,0.77,4,92,30,0.54,4,88,26,0.49,3,138,41,0.48,3,199
5--9,20,0.68,4,87,32,0.32,4,88,26,0.3,3,132,31,0.14,3,191
5--10,32,0.53,5,114,14,0.69,3,67,28,0.39,3,137,17,0.67,2,138
--,--,--,--,--,--,--,--,--,--,--,--,--,--,--,--,--
6--1,56,0.54,5,142,28,0.77,3,80,31,0.27,3,163,29,0.57,2,161
6--2,169,0.68,5,141,77,0.74,3,83,75,0.36,3,166,47,0.86,2,162
6--3,102,0.84,4,116,77,0.39,3,85,63,0.23,3,174,53,0.42,2,174
6--4,126,0.94,5,145,156,0.48,4,114,81,0.99,3,174,182,0.23,3,260
6--5,215,0.44,7,194,134,0.47,5,120,118,0.52,4,207,90,0.79,3,232
6--6,99,0.44,6,171,95,0.86,5,126,49,0.83,3,168,85,0.47,3,245
6--7,82,0.75,5,143,33,0.99,3,81,53,0.27,3,174,40,0.66,2,171
6--8,132,0.93,5,145,174,0.35,5,130,172,0.22,4,232,211,0.12,4,339
6--9,176,0.62,6,174,76,0.81,4,106,66,0.65,3,174,204,0.35,3,253
6--10,145,0.57,5,145,96,0.42,4,111,58,0.4,3,174,80,0.28,3,253
--,--,--,--,--,--,--,--,--,--,--,--,--,--,--,--,--
7--1,2205,0.48,5,202,223,0.38,4,133,267,0.46,3,206,73,0.95,2,204
7--2,345,0.68,5,175,284,0.46,4,122,271,0.2,3,208,386,0.15,3,302
7--3,308,0.47,4,140,225,0.31,4,134,87,0.99,2,140,101,0.35,2,210
7--4,438,0.72,4,140,179,0.89,3,102,367,0.51,3,210,194,0.59,2,207
7--5,1280,0.53,4,140,1568,0.35,4,131,322,0.51,2,140,293,0.46,2,210
7--6,719,0.77,6,207,1013,0.54,4,129,1156,0.47,3,198,610,0.87,2,201
7--7,508,0.61,4,140,814,0.92,4,132,323,0.43,3,209,824,0.29,3,306
7--8,135,0.84,5,175,286,0.7,4,126,122,0.44,3,207,313,0.23,3,299
7--9,198,0.98,4,139,198,0.38,4,133,185,0.19,3,206,143,0.57,2,203
7--10,303,0.85,5,169,812,0.28,5,138,169,0.52,3,198,727,0.21,3,289
--,--,--,--,--,--,--,--,--,--,--,--,--,--,--,--,--
8--1,2207,0.48,5,202,**,1.09,3,115,1405,0.33,3,240,1011,0.2,3,351
8--2,806,0.83,5,205,933,0.61,3,116,796,0.51,3,246,**,0.29 (1.54),3,364
8--3,665,0.71,5,202,1209,0.48,4,143,649,0.45,3,240,1026,0.31,3,350
8--4,978,0.65,5,199,443,0.42,3,113,895,0.15,3,232,228,0.64,2,228
8--5,2666,0.86,6,241,1998,0.96,4,151,3531,0.27,4,318,**,1.86,2,231
8--6,**,3.75,3,122,**,3.13,3,118,**,1.15,3,244,**,2.26,2,243
8--7,964,0.83,4,162,1205,0.68,3,120,929,0.47,3,242,630,0.89,2,240
8--8,**,0.49,6,245,3114,0.79,4,152,2418,0.32,4,324,2285,0.53,3,355
8--9,2966,0.48,5,199,2305,0.23,4,155,1157,0.14,3,246,858,0.55,2,243
8--10,796,0.88,5,202,2614,0.76,4,143,1849,0.28,3,236,657,0.18,3,333
\end{filecontents*}
\csvreader[
column count = 17,
    no head,
    before reading=\small\sisetup{round-mode=places, round-precision=2},
    after reading=\normalsize,
    filter expr = {test{\ifnumgreater{\thecsvinputline}{22}}},
    longtable=|c|rrr|rrr|rrr|rrr|,
    table head = \caption{Principal domain results.\label{tab:full reformulation}}\\ \hline
     & \multicolumn{3}{c|}{Bisection} &\multicolumn{3}{c|}{Direct} &\multicolumn{3}{c|}{NU2} & \multicolumn{3}{c|}{NU3} \\ 
    \hline 
    \bfseries instance & \bfseries t & 
    \bfseries iter & 
    \bfseries bin &
    \bfseries t & 
    \bfseries iter. & 
    \bfseries bin &
    \bfseries t & 
    \bfseries iter. & 
    \bfseries bin &
    \bfseries t & 
    \bfseries iter. & 
    \bfseries bin
     \\ 
    \hline
    \endfirsthead
    \caption{Principal domain results (continued).}\\
    \hline 
    \bfseries instance & \bfseries t & 
    \bfseries iter & 
    \bfseries bin &
    \bfseries t & 
    \bfseries iter. & 
    \bfseries bin &
    \bfseries t & 
    \bfseries iter. & 
    \bfseries bin &
    \bfseries t & 
    \bfseries iter. & 
    \bfseries bin
     \\ 
    \hline\endhead
    \hline\endfoot,
    late after line=\\, 
    ]{Data/FullReformulation.csv}{
        1=\instance, 
        2=\tbisect, 
        3=\gapbisect, 
        4=\iterbisect,
        5=\binarybisect,
        6=\tdirect, 
        7=\gapdirect,
        8=\iterdirect,
        9=\binarydirect,
        10=\tnutwo,
        11=\gapnutwo,
        12=\iternutwo,
        13=\binarynutwo,
        14=\tnuthree,
        15=\gapnuthree,
        16=\iternuthree,
        17=\binarynuthree}%
    {\instance & \tbisect & \iterbisect & \binarybisect & \tdirect & \iterdirect & \binarydirect & \tnutwo & \iternutwo & \binarynutwo & \tnuthree & \iternuthree & \binarynuthree}


%% file: Standalone_figures/k_worst_bisection.tex



\setlength\tabcolsep{3pt}
\begin{filecontents*}{Data/kWorstBisection.csv}
2--1,0,0.25,4,12,0,0.54,3,12,0,0.51,3,15
2--2,0,0.13,6,18,0,0.42,5,20,0,0.45,4,20
2--3,0,0.8,5,15,0,0.75,4,16,0,0.68,4,20
2--4,0,0.13,5,15,0,0.13,5,15,0,0.13,5,15
2--5,0,0.95,6,18,0,0.74,6,24,1,0.27,6,30
2--6,0,0.43,5,15,0,0.42,5,20,0,0.43,4,20
2--7,0,0.28,2,6,0,0.28,2,8,0,0.28,2,10
2--8,0,0.66,6,18,0,0.42,6,24,0,0.93,5,25
2--9,0,0.15,7,21,0,0.89,5,20,0,0.89,4,20
2--10,0,0.99,3,9,0,0.8,2,8,0,0.65,2,10
--,--,--,--,--,--,--,--,--,--,--,--,--
3--1,1,0.76,6,36,1,0.78,4,36,1,0.77,4,44
3--2,2,0.23,7,42,1,0.95,5,45,1,0.88,5,55
3--3,1,0.98,5,30,2,0.34,5,45,2,0.28,5,55
3--4,2,0.71,8,48,1,0.99,5,45,2,0.81,5,55
3--5,2,0.63,7,42,2,0.62,6,53,2,0.61,6,63
3--6,2,0.83,6,36,2,0.97,5,45,2,0.93,5,55
3--7,1,0.66,6,36,1,0.76,4,36,1,0.61,4,44
3--8,2,0.25,7,42,1,0.56,5,44,1,0.56,5,52
3--9,1,0.81,4,24,1,0.64,4,36,1,0.57,4,44
3--10,2,0.85,6,36,2,0.63,5,45,2,0.61,5,55
--,--,--,--,--,--,--,--,--,--,--,--,--
4--1,4,0.82,5,45,3,0.82,4,51,4,0.73,4,68
4--2,5,0.79,7,63,3,0.7,5,65,5,0.59,5,85
4--3,3,0.55,5,45,3,0.85,4,52,2,0.91,3,51
4--4,3,0.63,6,54,3,0.63,5,65,3,0.27,5,85
4--5,19,0.81,8,72,8,0.78,6,78,12,0.9,6,99
4--6,5,0.97,6,54,5,0.59,5,65,9,0.58,5,82
4--7,5,0.89,7,63,6,0.72,7,90,7,0.72,7,111
4--8,6,0.8,6,54,6,0.67,5,65,5,0.66,5,85
4--9,15,0.47,6,54,12,0.78,5,65,12,0.74,5,82
4--10,4,0.51,6,54,3,0.39,5,65,8,0.46,4,66
--,--,--,--,--,--,--,--,--,--,--,--,--
5--1,19,0.81,4,48,22,0.63,4,72,21,0.61,4,89
5--2,36,0.85,6,72,30,0.84,5,90,40,0.79,5,115
5--3,27,0.58,6,72,25,0.78,5,90,16,0.85,4,92
5--4,40,0.85,6,72,34,0.63,5,90,44,0.52,5,110
5--5,20,0.97,6,72,21,0.9,5,90,33,0.84,5,115
5--6,44,0.55,7,84,34,0.78,6,107,40,0.69,6,128
5--7,36,0.51,7,84,17,0.97,4,72,16,0.93,4,89
5--8,24,0.96,5,60,28,0.69,5,90,22,0.77,4,92
5--9,33,0.76,6,72,16,0.99,4,72,20,0.68,4,87
5--10,30,0.88,6,71,25,0.62,5,90,32,0.53,5,114
--,--,--,--,--,--,--,--,--,--,--,--,--
6--1,37,0.76,5,75,55,0.6,5,110,56,0.54,5,142
6--2,169,0.6,6,90,155,0.73,5,110,169,0.68,5,141
6--3,105,0.69,6,90,84,0.87,4,88,102,0.84,4,116
6--4,138,0.65,7,105,123,0.68,6,132,126,0.94,5,145
6--5,138,0.72,8,119,187,0.52,7,152,215,0.44,7,194
6--6,107,0.59,7,105,105,0.52,6,132,99,0.44,6,171
6--7,62,0.6,6,90,49,0.82,5,110,82,0.75,5,143
6--8,164,0.62,7,105,88,0.96,5,110,132,0.93,5,145
6--9,167,0.93,7,104,154,0.66,6,132,176,0.62,6,174
6--10,86,0.68,6,90,77,0.78,5,110,145,0.57,5,145
--,--,--,--,--,--,--,--,--,--,--,--,--
7--1,489,0.87,7,126,157,0.76,5,135,2205,0.48,5,202
7--2,327,0.75,6,108,286,0.88,5,135,345,0.68,5,175
7--3,274,0.6,5,90,104,0.56,4,108,308,0.47,4,140
7--4,366,0.52,6,108,364,0.85,4,108,438,0.72,4,140
7--5,524,0.73,5,90,703,0.67,4,108,1280,0.53,4,140
7--6,726,0.55,7,126,616,0.78,6,162,719,0.77,6,207
7--7,895,0.73,7,126,560,0.63,5,135,508,0.61,4,140
7--8,237,0.83,6,108,183,0.89,5,135,135,0.84,5,175
7--9,176,0.54,6,108,209,0.63,5,135,198,0.98,4,139
7--10,878,0.63,6,108,164,0.94,5,135,303,0.85,5,169
--,--,--,--,--,--,--,--,--,--,--,--,--
8--1,1361,0.9,6,126,1143,0.72,5,155,2207,0.48,5,202
8--2,3363,0.49,8,168,1111,0.96,5,155,806,0.83,5,205
8--3,1128,0.7,7,147,452,0.81,5,155,665,0.71,5,202
8--4,500,0.78,6,126,1042,0.73,5,155,978,0.65,5,199
8--5,2789,1,7,147,**,1.06,6,186,2666,0.86,6,241
8--6,**,1.93,5,105,**,2.6,4,124,**,3.75,3,122
8--7,1312,0.86,6,126,643,0.97,4,124,964,0.83,4,162
8--8,**,1.88,6,126,2464,0.62,6,185,**,0.49,6,245
8--9,**,1.25,5,105,2232,0.54,5,155,2966,0.48,5,199
8--10,1862,0.5,7,147,1219,0.62,6,186,796,0.88,5,202
\end{filecontents*}
\csvreader[
column count = 13,
    no head,
    before reading=\small\sisetup{round-mode=places, round-precision=2},
    after reading=\normalsize,
    filter expr = {test{\ifnumgreater{\thecsvinputline}{33}}}, 
    longtable=|c|rrr|rrr|rrr|,
    table head = \caption{k-Worst Results - Bisection.\label{tab:k worst bisection}}\\ \hline
     & \multicolumn{3}{c|}{k = 25 \%} &\multicolumn{3}{c|}{k = 50 \%} &\multicolumn{3}{c|}{k = 100 \%} \\ 
    \hline 
    \bfseries instance & \bfseries t & 
    \bfseries iter. & 
    \bfseries bin &
    \bfseries t & 
    \bfseries iter. & 
    \bfseries bin &
    \bfseries t & 
    \bfseries iter. & 
    \bfseries bin
     \\ 
    \hline
    \endfirsthead
    \caption{k-Worst Results - Bisection (continued).}\\
    \hline 
    \bfseries instance & \bfseries t & 
    \bfseries iter. & 
    \bfseries bin &
    \bfseries t & 
    \bfseries iter. & 
    \bfseries bin &
    \bfseries t & 
    \bfseries iter. & 
    \bfseries bin
     \\ 
    \hline\endhead
    \hline\endfoot,
    late after line=\\, 
    ]{Data/kWorstBisection.csv}{
        1=\instance, 
        2=\tone, 
        3=\gapone, 
        4=\iterone,
        5=\binaryone,
        6=\ttwo, 
        7=\gaptwo,
        8=\itertwo,
        9=\binarytwo,
        10=\tthree,
        11=\gapthree,
        12=\iterthree,
        13=\binarythree}%
    {\instance & \tone & \iterone & \binaryone & \ttwo & \itertwo & \binarytwo & \tthree & \iterthree & \binarythree}


%% file: Standalone_figures/k_worst_direct.tex



\setlength\tabcolsep{3pt}
\begin{filecontents*}{Data/kWorstDirect.csv}
2--1,0,0.43,3,9,0,0.3,3,11,0,0.3,3,11
2--2,1,0.84,4,12,1,0.72,4,16,0,0.14,4,18
2--3,0,0.62,3,9,1,0.93,2,7,0,0.63,3,9
2--4,0,0.39,3,8,0,0.39,3,8,0,0.39,3,8
2--5,0,0.44,6,18,1,0.55,4,16,0,0.25,4,19
2--6,0,0.32,4,12,1,0.3,4,13,0,0.77,3,14
2--7,0,0.66,1,3,0,0.62,1,4,0,0.62,1,5
2--8,0,0.94,3,9,0,0.78,3,12,0,0.78,3,12
2--9,0,0.47,4,12,0,0.34,4,16,0,0.28,4,16
2--10,0,0.34,3,9,0,0.5,2,8,0,0.07,2,9
--,--,--,--,--,--,--,--,--,--,--,--,--
3--1,1,0.23,4,24,1,0.99,2,18,0,0.76,2,22
3--2,1,0.47,5,30,1,0.33,4,35,1,0.32,4,37
3--3,1,0.5,4,24,1,0.55,3,24,1,0.52,3,29
3--4,1,0.34,6,36,1,0.73,5,45,1,0.7,4,43
3--5,1,0.65,5,30,2,0.38,5,43,2,0.38,5,47
3--6,1,0.71,5,30,1,0.83,4,32,1,0.84,4,32
3--7,0,0.93,3,18,0,0.45,3,27,0,0.84,2,21
3--8,1,0.7,4,24,1,0.39,4,35,1,0.37,4,34
3--9,1,0.17,5,30,1,0.17,4,36,1,0.17,4,37
3--10,2,0.58,5,29,2,0.11,5,43,1,0.8,4,40
--,--,--,--,--,--,--,--,--,--,--,--,--
4--1,4,0.74,5,44,2,0.51,3,39,2,0.47,3,46
4--2,3,0.97,5,45,3,0.27,5,64,3,0.76,4,59
4--3,2,0.83,3,27,1,0.78,2,26,1,0.89,2,34
4--4,3,0.77,5,45,3,0.48,4,48,3,0.35,4,54
4--5,10,0.83,8,72,9,0.68,6,77,12,0.64,6,92
4--6,2,0.73,3,27,3,0.58,3,38,3,0.47,3,46
4--7,5,0.5,7,63,4,0.72,5,65,3,0.79,4,60
4--8,5,0.38,5,45,3,0.7,3,39,3,0.54,3,50
4--9,11,0.51,5,45,7,0.78,4,52,10,0.73,4,62
4--10,3,0.3,4,35,2,0.53,3,39,2,0.61,3,42
--,--,--,--,--,--,--,--,--,--,--,--,--
5--1,13,0.62,3,36,9,0.99,2,36,9,0.86,2,45
5--2,29,0.61,5,59,19,0.67,4,72,20,0.95,3,69
5--3,36,0.43,5,59,24,0.47,4,72,31,0.46,4,84
5--4,47,0.41,6,72,37,0.68,5,89,33,0.79,4,82
5--5,19,0.77,5,60,18,0.88,4,72,24,0.72,4,87
5--6,23,0.67,4,48,38,0.56,5,86,33,0.55,5,89
5--7,14,0.84,4,47,10,0.71,3,54,11,0.51,3,66
5--8,26,0.57,5,58,20,0.73,4,70,30,0.54,4,88
5--9,27,0.45,5,60,24,0.46,4,72,32,0.32,4,88
5--10,16,0.87,4,46,13,0.87,3,54,14,0.69,3,67
--,--,--,--,--,--,--,--,--,--,--,--,--
6--1,24,0.76,4,59,20,0.81,3,66,28,0.77,3,80
6--2,69,0.78,4,60,44,0.78,3,66,77,0.74,3,83
6--3,60,0.62,4,60,61,0.54,3,66,77,0.39,3,85
6--4,144,0.59,6,90,122,0.46,4,88,156,0.48,4,114
6--5,117,0.58,6,90,135,0.55,5,105,134,0.47,5,120
6--6,119,0.38,6,90,82,0.92,5,107,95,0.86,5,126
6--7,45,0.92,4,60,64,0.34,4,87,33,0.99,3,81
6--8,175,0.87,6,90,138,0.38,5,110,174,0.35,5,130
6--9,121,0.74,5,75,80,0.9,4,88,76,0.81,4,106
6--10,65,0.96,4,60,70,0.37,4,88,96,0.42,4,111
--,--,--,--,--,--,--,--,--,--,--,--,--
7--1,286,0.71,5,89,177,0.43,4,108,223,0.38,4,133
7--2,319,0.39,5,88,451,0.5,4,107,284,0.46,4,122
7--3,170,0.45,5,88,127,0.56,4,108,225,0.31,4,134
7--4,395,0.9,5,90,240,0.86,4,108,179,0.89,3,102
7--5,444,0.79,4,72,817,0.49,4,108,1568,0.35,4,131
7--6,322,0.45,5,90,577,0.56,4,108,1013,0.54,4,129
7--7,308,0.96,5,90,1247,0.33,5,134,814,0.92,4,132
7--8,589,0.27,5,90,216,0.84,4,108,286,0.7,4,126
7--9,257,0.51,5,90,221,0.57,4,108,198,0.38,4,133
7--10,591,0.37,5,88,617,0.31,4,105,812,0.28,5,138
--,--,--,--,--,--,--,--,--,--,--,--,--
8--1,2599,0.85,5,104,1941,0.67,4,124,**,1.09,3,115
8--2,2041,0.76,5,105,959,0.61,4,124,933,0.61,3,116
8--3,648,0.56,5,104,477,0.69,4,124,1209,0.48,4,143
8--4,517,0.81,4,84,290,0.8,3,93,443,0.42,3,113
8--5,2533,0.64,6,125,3599,0.48,5,155,1998,0.96,4,151
8--6,**,2.95,4,84,3599,0.83,4,122,**,3.13,3,118
8--7,1343,0.33,5,102,928,0.8,3,93,1205,0.68,3,120
8--8,3530,0.65,6,125,1886,0.97,4,124,3114,0.79,4,152
8--9,2041,0.39,5,104,2279,0.33,4,122,2305,0.23,4,155
8--10,954,0.87,5,104,933,0.83,4,124,2614,0.76,4,143
\end{filecontents*}
\csvreader[
column count = 13,
    no head,
    before reading=\small\sisetup{round-mode=places, round-precision=2},
    after reading=\normalsize,
    filter expr = {test{\ifnumgreater{\thecsvinputline}{33}}}, 
    longtable=|c|rrr|rrr|rrr|,
    table head = \caption{k-Worst Results - Direct.\label{tab:k worst direct}}\\ \hline
     & \multicolumn{3}{c|}{k = 25 \%} &\multicolumn{3}{c|}{k = 50 \%} &\multicolumn{3}{c|}{k = 100 \%} \\ 
    \hline 
    \bfseries instance & \bfseries t & 
    \bfseries iter & 
    \bfseries bin &
    \bfseries t & 
    \bfseries iter. & 
    \bfseries bin &
    \bfseries t & 
    \bfseries iter. & 
    \bfseries bin
     \\ 
    \hline
    \endfirsthead
    \caption{k-Worst Results - Direct (continued).}\\
    \hline 
    \bfseries instance & \bfseries t & 
    \bfseries iter & 
    \bfseries bin &
    \bfseries t & 
    \bfseries iter. & 
    \bfseries bin &
    \bfseries t & 
    \bfseries iter. & 
    \bfseries bin
     \\ 
    \hline\endhead
    \hline\endfoot,
    late after line=\\, 
    ]{Data/kWorstDirect.csv}{
        1=\instance, 
        2=\tone, 
        3=\gapone, 
        4=\iterone,
        5=\binaryone,
        6=\ttwo, 
        7=\gaptwo,
        8=\itertwo,
        9=\binarytwo,
        10=\tthree,
        11=\gapthree,
        12=\iterthree,
        13=\binarythree}%
    {\instance & \tone & \iterone & \binaryone & \ttwo & \itertwo & \binarytwo & \tthree & \iterthree & \binarythree}


%% file: Standalone_figures/k_worst_nu2.tex



\setlength\tabcolsep{3pt}
\begin{filecontents*}{Data/kWorstNu2.csv}
2--1,0,0.44,3,18,0,0.38,2,16,0,0.26,2,20
2--2,0,0.78,3,18,0,0.63,3,24,0,0.95,2,20
2--3,0,0.42,3,18,0,0.46,2,16,0,0.4,2,20
2--4,0,0.79,3,18,0,0.79,3,18,0,0.79,3,18
2--5,0,0.4,4,24,0,0.44,3,24,0,0.87,3,30
2--6,0,0.55,3,18,0,0.52,3,24,0,0.44,3,30
2--7,0,0.06,2,12,0,0.01,2,16,0,0.01,2,20
2--8,0,0.92,3,18,0,0.85,2,16,0,0.59,2,20
2--9,0,0.49,4,24,0,0.26,3,24,0,0.19,3,29
2--10,0,0.93,2,12,0,0.39,2,16,0,0.97,1,10
--,--,--,--,--,--,--,--,--,--,--,--,--
3--1,1,0.28,4,48,1,0.24,4,72,2,0.22,4,88
3--2,1,0.54,4,48,2,0.39,4,71,2,0.39,4,85
3--3,1,0.85,3,36,1,0.52,3,54,1,0.44,3,66
3--4,1,0.59,5,60,1,0.4,4,72,1,0.78,3,66
3--5,1,0.98,3,36,1,0.99,3,52,1,0.96,3,60
3--6,1,0.97,4,48,1,0.67,3,54,1,0.74,3,66
3--7,0,0.82,3,36,1,0.35,3,52,1,0.32,3,66
3--8,1,0.2,5,60,1,0.86,3,54,1,0.72,3,66
3--9,1,0.62,3,36,1,0.62,2,36,1,0.57,2,44
3--10,2,0.34,5,60,1,0.39,3,54,1,0.22,3,66
--,--,--,--,--,--,--,--,--,--,--,--,--
4--1,4,0.26,4,72,4,0.3,3,78,5,0.26,3,102
4--2,3,0.98,4,72,4,0.33,4,104,3,0.85,3,102
4--3,2,0.64,3,52,1,0.57,2,52,2,0.47,2,68
4--4,2,0.82,3,54,3,0.79,3,78,2,0.72,3,100
4--5,26,0.33,8,144,17,0.3,5,130,11,0.54,4,130
4--6,3,0.77,3,54,4,0.46,3,78,9,0.41,3,96
4--7,4,0.87,5,90,4,0.36,4,104,4,0.28,4,129
4--8,4,0.86,3,54,4,0.42,3,78,4,0.52,3,102
4--9,12,0.65,4,72,8,0.44,3,78,10,0.48,3,99
4--10,4,0.23,4,72,3,0.42,3,78,3,0.43,3,98
--,--,--,--,--,--,--,--,--,--,--,--,--
5--1,13,0.62,2,48,14,0.54,2,72,11,0.5,2,92
5--2,31,0.65,4,96,21,0.59,3,108,21,0.53,3,138
5--3,15,0.52,3,72,17,0.29,3,108,18,0.22,3,138
5--4,26,0.67,4,96,22,0.95,3,108,23,0.89,3,134
5--5,22,0.44,4,96,14,0.46,3,108,33,0.56,3,138
5--6,28,0.4,4,96,22,0.77,3,106,23,0.73,3,121
5--7,12,0.81,3,72,18,0.34,3,108,18,0.21,3,138
5--8,19,0.92,3,72,24,0.75,3,108,26,0.49,3,138
5--9,13,0.78,3,72,28,0.31,3,108,26,0.3,3,132
5--10,21,0.49,4,95,21,0.35,3,108,28,0.39,3,137
--,--,--,--,--,--,--,--,--,--,--,--,--
6--1,28,0.82,3,90,44,0.32,3,132,31,0.27,3,163
6--2,95,0.34,4,120,54,0.29,3,132,75,0.36,3,166
6--3,38,0.93,3,90,54,0.28,3,132,63,0.23,3,174
6--4,71,0.82,4,120,54,0.88,3,132,81,0.99,3,174
6--5,114,0.42,5,150,88,0.59,4,172,118,0.52,4,207
6--6,96,0.29,5,150,41,0.79,3,132,49,0.83,3,168
6--7,40,0.71,3,90,52,0.33,3,132,53,0.27,3,174
6--8,139,0.4,5,150,137,0.32,4,176,172,0.22,4,232
6--9,97,0.51,4,120,104,0.82,3,132,66,0.65,3,174
6--10,86,0.39,4,120,47,0.45,3,132,58,0.4,3,174
--,--,--,--,--,--,--,--,--,--,--,--,--
7--1,112,0.71,4,144,98,0.54,3,162,267,0.46,3,206
7--2,285,0.77,3,108,215,0.34,3,162,271,0.2,3,208
7--3,140,0.64,3,108,134,0.29,3,162,87,0.99,2,140
7--4,387,0.55,4,144,170,0.48,3,162,367,0.51,3,210
7--5,493,0.57,3,108,304,0.91,2,108,322,0.51,2,140
7--6,574,0.86,3,108,641,0.46,3,162,1156,0.47,3,198
7--7,540,0.96,4,144,406,0.47,3,162,323,0.43,3,209
7--8,80,0.94,3,108,329,0.55,3,162,122,0.44,3,207
7--9,252,0.32,3,108,181,0.2,3,162,185,0.19,3,206
7--10,125,0.99,3,108,181,0.55,3,162,169,0.52,3,198
--,--,--,--,--,--,--,--,--,--,--,--,--
8--1,1789,0.33,4,168,297,0.54,3,186,1405,0.33,3,240
8--2,2986,0.85,4,168,1249,0.71,3,186,796,0.51,3,246
8--3,571,0.48,4,168,364,0.54,3,186,649,0.45,3,240
8--4,743,0.7,3,126,392,0.32,3,186,895,0.15,3,232
8--5,1599,0.51,5,210,1209,0.84,3,186,3531,0.27,4,318
8--6,**,1.63,4,168,**,1.61,3,186,**,1.15,3,244
8--7,1248,0.64,4,168,990,0.58,3,186,929,0.47,3,242
8--8,2555,0.53,5,208,1799,0.83,3,184,2418,0.32,4,324
8--9,2306,0.52,3,126,835,0.26,3,186,1157,0.14,3,246
8--10,468,0.47,4,168,297,0.75,3,186,1849,0.28,3,236
\end{filecontents*}
\csvreader[
column count = 13,
    no head,
    before reading=\small\sisetup{round-mode=places, round-precision=2},
    after reading=\normalsize,
    filter expr = {test{\ifnumgreater{\thecsvinputline}{33}}}, 
    longtable=|c|rrr|rrr|rrr|,
    table head = \caption{$k$-Worst Results - Non-Uniform Two-Point.\label{tab:k worst nu2}}\\ \hline
     & \multicolumn{3}{c|}{k = 25 \%} &\multicolumn{3}{c|}{k = 50 \%} &\multicolumn{3}{c|}{k = 100 \%} \\ 
    \hline 
    \bfseries instance & \bfseries t & 
    \bfseries iter & 
    \bfseries bin &
    \bfseries t & 
    \bfseries iter. & 
    \bfseries bin &
    \bfseries t & 
    \bfseries iter. & 
    \bfseries bin
     \\ 
    \hline
    \endfirsthead
    \caption{k-Worst Results - Non-Uniform Two-Point (continued).}\\
    \hline 
    \bfseries instance & \bfseries t & 
    \bfseries iter & 
    \bfseries bin &
    \bfseries t & 
    \bfseries iter. & 
    \bfseries bin &
    \bfseries t & 
    \bfseries iter. & 
    \bfseries bin
     \\ 
    \hline\endhead
    \hline\endfoot,
    late after line=\\, 
    ]{Data/kWorstNu2.csv}{
        1=\instance, 
        2=\tone, 
        3=\gapone, 
        4=\iterone,
        5=\binaryone,
        6=\ttwo, 
        7=\gaptwo,
        8=\itertwo,
        9=\binarytwo,
        10=\tthree,
        11=\gapthree,
        12=\iterthree,
        13=\binarythree}%
    {\instance & \tone & \iterone & \binaryone & \ttwo & \itertwo & \binarytwo & \tthree & \iterthree & \binarythree}


%% file: Standalone_figures/k_worst_nu3.tex



\setlength\tabcolsep{3pt}
\begin{filecontents*}{Data/kWorstNu3.csv}
2--1,0,0.16,3,27,0,0.22,2,24,0,0.22,2,26
2--2,0,0.59,3,27,0,0.47,3,36,0,0.09,3,43
2--3,0,0.49,2,18,0,0.39,2,21,0,0.39,2,25
2--4,0,0.31,3,24,0,0.31,3,24,0,0.31,3,24
2--5,0,0.74,3,27,0,0.19,3,34,0,0.63,3,43
2--6,0,0.37,3,27,0,0.33,3,36,0,0.27,3,42
2--7,0,0.66,1,9,0,0.62,1,12,0,0.62,1,15
2--8,0,0.63,2,18,0,0.49,2,24,0,0.49,2,28
2--9,0,0.5,3,27,0,0.19,3,36,0,0.1,3,37
2--10,0,0.65,2,18,0,0.83,1,12,0,0.65,1,15
--,--,--,--,--,--,--,--,--,--,--,--,--
3--1,1,0.42,3,53,1,0.69,2,54,1,0.59,2,66
3--2,2,0.29,4,72,1,0.11,3,81,1,0.08,3,89
3--3,1,0.43,3,54,1,0.82,2,54,1,0.55,2,63
3--4,2,0.46,5,90,1,0.95,3,81,1,0.32,3,97
3--5,2,0.37,4,72,1,0.78,3,78,2,0.78,3,88
3--6,1,0.47,3,54,1,0.54,3,81,1,0.54,3,85
3--7,1,0.29,3,54,0,0.7,2,54,0,0.48,2,63
3--8,1,0.23,4,72,1,0.5,3,81,1,0.44,3,93
3--9,1,0.51,3,51,1,0.61,2,54,1,0.52,2,61
3--10,2,0.18,5,90,1,0.74,3,81,1,0.06,3,97
--,--,--,--,--,--,--,--,--,--,--,--,--
4--1,3,0.68,3,81,3,0.9,2,78,2,0.88,2,99
4--2,2,0.9,3,81,3,0.74,3,117,4,0.51,3,143
4--3,3,0.42,3,81,2,0.51,2,78,2,0.36,2,99
4--4,3,0.79,3,81,3,0.56,3,117,3,0.57,3,143
4--5,14,0.79,5,132,15,0.29,4,156,22,0.25,4,191
4--6,5,0.39,3,81,3,0.68,2,78,3,0.58,2,93
4--7,5,0.45,5,134,6,0.14,4,156,6,0.13,4,174
4--8,7,0.64,3,81,10,0.33,3,117,14,0.19,3,152
4--9,10,0.99,3,80,13,0.2,3,117,6,0.94,2,99
4--10,3,0.52,3,81,4,0.25,3,116,4,0.22,3,137
--,--,--,--,--,--,--,--,--,--,--,--,--
5--1,11,0.59,2,72,14,0.45,2,108,16,0.37,2,137
5--2,47,0.53,4,144,38,0.48,3,162,43,0.45,3,206
5--3,25,0.33,3,108,32,0.11,3,162,28,0.06,3,203
5--4,29,0.69,3,108,38,0.31,3,160,43,0.3,3,184
5--5,22,0.58,3,108,24,0.34,3,162,27,0.3,3,199
5--6,26,0.41,3,108,29,0.16,3,159,54,0.18,3,183
5--7,19,0.59,3,108,12,0.89,2,108,15,0.77,2,135
5--8,27,0.73,3,108,34,0.51,3,162,41,0.48,3,199
5--9,29,0.43,3,108,29,0.15,3,160,31,0.14,3,191
5--10,18,0.89,3,108,15,0.9,2,108,17,0.67,2,138
--,--,--,--,--,--,--,--,--,--,--,--,--
6--1,29,0.66,3,135,22,0.71,2,132,29,0.57,2,161
6--2,151,0.82,3,135,49,0.83,2,132,47,0.86,2,162
6--3,57,0.82,3,135,42,0.48,2,132,53,0.42,2,174
6--4,64,0.76,3,132,95,0.74,3,198,182,0.23,3,260
6--5,122,0.42,4,180,74,0.81,3,192,90,0.79,3,232
6--6,105,0.24,4,180,61,0.49,3,195,85,0.47,3,245
6--7,66,0.54,3,135,41,0.85,2,132,40,0.66,2,171
6--8,128,0.52,4,180,159,0.14,4,264,211,0.12,4,339
6--9,160,0.35,4,180,111,0.28,3,198,204,0.35,3,253
6--10,67,0.8,3,135,62,0.34,3,198,80,0.28,3,253
--,--,--,--,--,--,--,--,--,--,--,--,--
7--1,137,0.74,3,162,205,0.26,3,243,73,0.95,2,204
7--2,313,0.59,3,160,334,0.18,3,243,386,0.15,3,302
7--3,151,0.31,3,162,186,0.2,3,243,101,0.35,2,210
7--4,220,0.59,3,162,562,0.22,3,243,194,0.59,2,207
7--5,1032,0.32,4,216,527,0.95,2,162,293,0.46,2,210
7--6,1265,0.65,3,162,410,0.98,2,162,610,0.87,2,201
7--7,876,0.43,4,216,533,0.33,3,243,824,0.29,3,306
7--8,151,0.6,3,162,340,0.24,3,243,313,0.23,3,299
7--9,124,0.72,2,108,107,0.63,2,162,143,0.57,2,203
7--10,325,0.38,3,162,521,0.25,3,243,727,0.21,3,289
--,--,--,--,--,--,--,--,--,--,--,--,--
8--1,1122,0.79,3,189,1101,0.32,3,278,1011,0.2,3,351
8--2,2005,0.26,4,252,1237,0.34,3,278,**,0.29 (1.54),3,364
8--3,1012,0.87,3,189,397,0.3,3,279,1026,0.31,3,350
8--4,292,0.33,3,189,110,0.98,2,186,228,0.64,2,228
8--5,2862,0.43,4,252,1584,0.57,3,279,**,1.86,2,231
8--6,3434,0.88,4,252,**,2.72,2,186,**,2.26,2,243
8--7,1188,0.65,3,189,386,0.93,2,186,630,0.89,2,240
8--8,**,2.55,3,189,1981,0.58,3,279,2285,0.53,3,355
8--9,2090,0.32,3,189,2433,0.79,2,186,858,0.55,2,243
8--10,799,0.23,4,252,1036,0.2,3,279,657,0.18,3,333
\end{filecontents*}
\csvreader[
column count = 13,
    no head,
    before reading=\small\sisetup{round-mode=places, round-precision=2},
    after reading=\normalsize,
    filter expr = {test{\ifnumgreater{\thecsvinputline}{33}}}, 
    longtable=|c|rrr|rrr|rrr|,
    table head = \caption{k-Worst Results - Non-Uniform Three-Point.\label{tab:k worst nu3}}\\ \hline
     & \multicolumn{3}{c|}{k = 25 \%} &\multicolumn{3}{c|}{k = 50 \%} &\multicolumn{3}{c|}{k = 100 \%} \\ 
    \hline 
    \bfseries instance & \bfseries t & 
    \bfseries iter & 
    \bfseries bin &
    \bfseries t & 
    \bfseries iter. & 
    \bfseries bin &
    \bfseries t & 
    \bfseries iter. & 
    \bfseries bin
     \\ 
    \hline
    \endfirsthead
    \caption{k-Worst Results - Non-Uniform Three-Point (continued).}\\
    \hline 
    \bfseries instance & \bfseries t & 
    \bfseries iter & 
    \bfseries bin &
    \bfseries t & 
    \bfseries iter. & 
    \bfseries bin &
    \bfseries t & 
    \bfseries iter. & 
    \bfseries bin
     \\ 
    \hline\endhead
    \hline\endfoot,
    late after line=\\, 
    ]{Data/kWorstNu3.csv}{
        1=\instance, 
        2=\tone, 
        3=\gapone, 
        4=\iterone,
        5=\binaryone,
        6=\ttwo, 
        7=\gaptwo,
        8=\itertwo,
        9=\binarytwo,
        10=\tthree,
        11=\gapthree,
        12=\iterthree,
        13=\binarythree}%
    {\instance & \tone & \iterone & \binaryone & \ttwo & \itertwo & \binarytwo & \tthree & \iterthree & \binarythree}
